\pdfoutput=1
\documentclass{amsart}
\usepackage[utf8]{inputenc}

\usepackage[margin=1in]{geometry}
\usepackage{microtype}

\usepackage{amsmath, amsfonts, amssymb, amsthm}
\usepackage{eucal}

\usepackage{tikz}
\usetikzlibrary{diagrams, decorations.pathreplacing}

\usepackage[pdfusetitle,unicode]{hyperref}

\let\itemref=\ref
\makeatletter
\def\theenumi{(\@arabic\c@enumi)}

\def\theenumii{(\@alph\c@enumii)}

\def\theenumiii{(\@roman\c@enumiii)}

\def\theenumiv{(\@Alph\c@enumiv)}

\makeatother

\def\resp{{\sfcode`\.1000 resp.}}
\def\ie{{\sfcode`\.1000 i.e.}}

\def\cf{{\sfcode`\.1000 cf.}}

\def\<{\discretionary{}{}{}}

\let\scr=\mathcal
\let\bb=\mathbf
\let\der=\mathbf
\let\cat=\mathcal
\let\smcat=\mathrm

\def\Z{\bb Z}
\def\Q{\bb Q}
\def\A{\bb A}
\def\P{\bb P}
\def\G{\bb G}

\def\L{\der L}
\def\R{\der R}

\def\ph{\mathord-}
\def\id{\mathrm{id}}
\def\tr{\mathrm{tr}}
\def\th{\vartheta}
\def\1{\mathbf 1}

\DeclareMathOperator{\Hom}{Hom{}}
\DeclareMathOperator{\Map}{Map{}}
\DeclareMathOperator{\cor}{Cor{}}
\DeclareMathOperator{\Tor}{Tor{}}
\DeclareMathOperator{\Spec}{Spec{}}
\DeclareMathOperator{\Th}{Th{}}
\DeclareMathOperator{\Pic}{Pic{}}
\DeclareMathOperator{\Char}{char{}}
\DeclareMathOperator{\essdim}{ess\,dim{}}
\DeclareMathOperator{\Sp}{Sp}

\let\phi=\varphi

\let\into=\hookrightarrow

\let\tens=\otimes
\let\pt=\ast
\let\minus=\smallsetminus
\def\cotens{\mathbin\Box}
\def\suchthat{\:\vert\:}
\def\abs#1{\lvert #1\rvert}

\def\Gr{\mathrm{Gr}{}}
\def\MU{\mathrm{MU}{}}
\def\BGL{\mathrm{BGL}{}}
\def\MGL{\mathrm{MGL}{}}
\def\BP{\mathrm{BP}{}}
\def\MH{\Lambda}

\def\Spc{\cat S\smcat{pc}{}}
\def\Spt{\cat S\smcat{pt}{}}
\def\DM{\cat{DM}{}}
\def\Ab{\cat A\smcat{b}{}}
\def\H{\cat H}
\def\SH{\cat{SH}{}}
\def\D{\cat D}
\def\Mod{\cat M\smcat{od}{}}

\def\Sm{\cat S\smcat{m}{}}
\def\Cor{\cat C\smcat{or}{}}

\def\s{\Delta^\op}
\def\op{\mathrm{op}}
\def\eff{\mathrm{eff}}
\def\Nis{\mathrm{Nis}}
\def\et{\mathrm{\acute et}}

\def\aast{{\ast\ast}}
\def\spi{\underline{\pi}}

\let\lim=\relax
\DeclareMathOperator*{\lim}{lim}
\DeclareMathOperator*{\Rlim}{lim^1}
\DeclareMathOperator*{\colim}{colim}
\DeclareMathOperator*{\hocolim}{hocolim}
\DeclareMathOperator*{\holim}{holim}

\let\Im=\relax
\DeclareMathOperator{\Im}{Im}

\numberwithin{equation}{section}

\theoremstyle{plain}
\newtheorem{theorem}[equation]{Theorem}
\newtheorem{proposition}[equation]{Proposition}
\newtheorem{lemma}[equation]{Lemma}
\newtheorem{corollary}[equation]{Corollary}

\theoremstyle{definition}
\newtheorem{definition}[equation]{Definition}
\newtheorem{example}[equation]{Example}

\theoremstyle{remark}
\newtheorem{remark}[equation]{Remark}

\title{From algebraic cobordism to motivic cohomology}
\author{Marc Hoyois}
\date{\today}
\address{Department of Mathematics, Northwestern University, 2033 Sheridan Road, Evanston, IL 60208, USA}
\email{hoyois@math.northwestern.edu}
\urladdr{\url{http://math.northwestern.edu/~hoyois/}}


\begin{document}

\begin{abstract}
		Let $S$ be an essentially smooth scheme over a field of characteristic exponent $c$.
		We prove that there is a canonical equivalence of motivic spectra over $S$
	\[\MGL/(a_1,a_2,\dotsc)[1/c]\simeq H\Z[1/c],\]
	where $H\Z$ is the motivic cohomology spectrum, $\MGL$ is the algebraic cobordism spectrum, and the elements $a_n$ are generators of the Lazard ring.
	We discuss several applications including the computation of the slices of $\Z[1/c]$-local Landweber exact motivic spectra and the convergence of the associated slice spectral sequences.
\end{abstract}

\maketitle
\tableofcontents
\newpage
\section{Introduction}

Complex cobordism plays a central role in stable homotopy theory as the universal complex-oriented cohomology theory, and one of the most fruitful advances in the field was Quillen's identification of the complex cobordism of a point with the Lazard ring $L\cong\Z[a_1,a_2,\dotsc]$. Quillen's theorem can be rephrased purely in terms of spectra as a canonical equivalence $\MU/(a_1,a_2,\dotsc)\simeq H\Z$, where $\MU$ is the spectrum representing complex cobordism and $H\Z$ the spectrum representing ordinary cohomology. This formulation immediately suggests a motivic version of Quillen's theorem: in the motivic world, the universal oriented cohomology theory is represented by the algebraic cobordism spectrum $\MGL$, and the theory of Chern classes in motivic cohomology provides a canonical map
\begin{equation*}\label{eqn:intro}\tag{\textasteriskcentered}\MGL/(a_1,a_2,\dotsc)\to H\Z\end{equation*}
where $H\Z$ is now the spectrum representing motivic cohomology. By analogy with the topological situation, one is tempted to conjecture that~\eqref{eqn:intro} is an equivalence of motivic spectra. In this paper we prove this conjecture over fields of characteristic zero (and, more generally, over essentially smooth schemes over such fields). This result was previously announced by Hopkins and Morel, but their proof was never published. Our proof is essentially reverse-engineered from a talk given by Hopkins at Harvard in the fall of 2004 (recounted in \cite{Hopkins:2004}). For another account of the work of Hopkins and Morel, see \cite{Ayoub:2005}.

Several applications of this equivalence are already known. In \cite{Spitzweck:2010,Spitzweck:2012}, Spitzweck computes the slices of Landweber exact spectra: if $E$ is the motivic spectrum associated with a Landweber exact $L$-module $M_\ast$, then its $q$th slice $s_qE$ is the shifted Eilenberg–Mac Lane spectrum $\Sigma^{2q,q}HM_{q}$ (for $E=\MGL$, this should be taken as the motivic analogue of Quillen's computation of the homotopy groups of $\MU$). This shows that one can approach $E$-cohomology from motivic cohomology with coefficients in $M_\ast$ by means of a spectral sequence, generalizing the classical spectral sequence for algebraic $K$-theory. In \cite{Levine:2013}, Levine computes the slices of the motivic sphere spectrum in terms of the $\G_m$-stack of strict formal groups. In \cite{Levine:2009}, the Levine–Morel algebraic cobordism $\Omega^\ast(\ph)$ is identified with $\MGL^{(2,1)\ast}(\ph)$ on smooth schemes; in particular, the formal group law on $\MGL_{(2,1)\ast}$ is universal. Another interesting fact which follows at once from the equivalence~\eqref{eqn:intro} is that $H\Z$ is a cellular spectrum (\ie, an iterated homotopy colimit of stable motivic spheres). We will review all these applications at the end of the paper.

Assume now that the base field has characteristic $p>0$ (or, more generally, that the base scheme is essentially smooth over such a field). We will then show that~\eqref{eqn:intro} induces an equivalence
\[\MGL/(a_1,a_2,\dotsc)[1/p]\simeq H\Z[1/p].\]
This (partial) extension of the Hopkins–Morel equivalence is made possible by the recent computation of the motivic Steenrod algebra in positive characteristic (\cite{HKO:2013}).
It is not clear whether~\eqref{eqn:intro} itself is an equivalence in this case; it would suffice to show that the induced map
\[H\Z/p\wedge \MGL/(a_1,a_2,\dotsc)\to H\Z/p\wedge  H\Z\]
is an equivalence over the finite field $\bb F_p$. The left-hand side can easily be computed as an $H\Z/p$-module, but the existing methods to compute $H\Z/\ell\wedge H\Z$ for primes $\ell\neq p$ (namely, representing groups of algebraic cycles by symmetric powers and studying the motives of the latter using resolutions of singularities) all fail when $\ell=p$. In particular, it remains unknown whether $H\Z$ is a cellular spectrum over fields of positive characteristic.

\begingroup
\def\nocontent#1#2#3{}
\let\addcontentsline=\nocontent

\subsection*{Outline of the proof}
Assume for simplicity that the base is a field of characteristic zero and let $f\colon \MGL/(a_1,a_2,\dotsc)\to H\Z$ be the map to be proved an equivalence. Then
\vskip1em
\[
\begingroup
\def\vp{\vphantom{\displaystyle\sum}}
\hskip-1em\left.
\begin{array}{r}
	\smash{
	\left.
	\begin{array}{r}
		\text{\itemref{itm:outline:1} $H\Q\wedge f$ is an equivalence} \\
		\vp\\
		\text{\itemref{itm:outline:2} $H\Z/\ell\wedge f$ is an equivalence} \\
	\end{array}
	\right\rbrace\stackrel{\textstyle\itemref{itm:outline:3}\,\vp}{\Longrightarrow} \text{$H\Z\wedge f$ is an equivalence}
	} \\
	\vp\\
	\text{\itemref{itm:outline:4} $\MGL_{\leq 0}\simeq H\Z_{\leq 0}$} \\
\end{array}
\right\rbrace\stackrel{\textstyle\itemref{itm:outline:5}\,\vp}{\Longrightarrow} \text{$f$ is an equivalence,}
\endgroup
\]
\vskip.5em
\noindent
where $\ell$ is any prime number and $(\ph)_{\leq 0}$ is the truncation for Morel's homotopy $t$-structure. Here follows a summary of each key step (references are given in the main text).
\begin{enumerate}
	\item\label{itm:outline:1} This is a straightforward consequence of the work of Naumann, Spitzweck, and Østvær on motivic Landweber exactness, more specifically of the fact that $H\Q$ is the Landweber exact spectrum associated with the additive formal group over $\Q$.
	\item\label{itm:outline:2} $H\Z/\ell\wedge H\Z$ can be computed using Voevodsky's work on the motivic Steenrod algebra and motivic Eilenberg–Mac Lane spaces: it is a cellular $H\Z/\ell$-module and its homotopy groups are the kernel of the Bockstein acting on the dual motivic Steenrod algebra. To apply Voevodsky's results we also need the fact proved by Röndigs and Østvær that ``motivic spectra with $\Z/\ell$-transfers'' are equivalent to $H\Z/\ell$-modules. We compute the homotopy groups of $H\Z/\ell\wedge \MGL/(a_1,a_2,\dotsc)$ by elementary means, and direct comparison then shows that $H\Z/\ell\wedge f$ is an isomorphism on homotopy groups, whence an equivalence by cellularity.
	\item\label{itm:outline:3} This is a simple algebraic result.
	\item\label{itm:outline:4} By a theorem of Morel it suffices to show that $\MGL_{\leq 0}\to H\Z_{\leq 0}$ induces isomorphisms on the stalks of the homotopy sheaves at field extensions $L$ of $k$. For $H\Z_{\leq 0}$ these stalks are given by the motivic cohomology groups $H^{n,n}(\Spec L,\Z)$, which have been classically identified with the Milnor $K$-theory groups of the field $L$. We identify the homotopy sheaves of $\MGL_{\leq 0}$ with the cokernel of the Hopf element acting on the homotopy sheaves of the truncated sphere spectrum $\1_{\leq 0}$; it then follows from Morel's explicit computation of the latter that the stalks of the homotopy sheaves of $\MGL_{\leq 0}$ over $\Spec L$ are also the Milnor $K$-theory groups of $L$.
	\item\label{itm:outline:5} To complete the proof we show that $\MGL/(a_1,a_2,\dotsc)$ is $H\Z$-local. Since the homotopy $t$-structure is left complete, it suffices to show that the truncations of $\MGL/(a_1,a_2,\dotsc)$ are $H\Z_{\leq 0}$-local. As proved by Gutiérrez, Röndigs, Spitzweck, and Østvær, these truncations are modules over $\MGL_{\leq 0}$ and hence are $\MGL_{\leq 0}$-local. By \itemref{itm:outline:4}, they are also $H\Z_{\leq 0}$-local.
\end{enumerate}

\subsection*{Acknowledgements}
I am very grateful to Paul Arne Østvær and Markus Spitzweck for pointing out a mistake in the proof of the main theorem in a preliminary version of this text, and to Paul Goerss for many helpful discussions. I thank Lukas Brantner for pointing out a mistake in the proof of Proposition~\ref{prop:piMGL} in the published version of this article.

\endgroup

\section{Preliminaries}

Let $S$ be a Noetherian scheme of finite Krull dimension; we will call such a scheme a \emph{base scheme}. Let $\Sm/S$ be the category of separated smooth schemes of finite type over $S$, or \emph{smooth schemes} for short. 
We denote by $\Spc(S)$ and $\Spc_\pt(S)$ the categories of simplicial presheaves and of pointed simplicial presheaves on $\Sm/S$. The motivic spheres $S^{p,q}\in\Spc_\pt(S)$ are defined by
\[S^{p,q}=(S^1)^{\wedge (p-q)}\wedge\tilde\G_m^{\wedge q}\]
where $\tilde\G_m=\Delta^1\vee\G_m$ is pointed away from $\G_m$ (the reason we do not use $\G_m$ itself is that we need $S^{2,1}$ to be projectively cofibrant for some of the model structures considered in \S\ref{sub:transfers} to exist; the reader can safely ignore this point).
We denote by $\Spt(S)$ the category of symmetric $S^{2,1}$-spectra in $\Spt_\pt(S)$, by $\Sigma^\infty\colon \Spc_\pt(S)\to\Spt(S)$ the stabilization functor, and by $\1\in\Spt(S)$ the sphere spectrum $\Sigma^\infty S_+$. We endow each of the categories $\Spc(S)$, $\Spc_\pt(S)$, and $\Spt(S)$ with its usual class of equivalences (often called motivic weak equivalences), and we denote by $\H(S)$, $\H_\pt(S)$, and $\SH(S)$ their respective homotopy categories.

Standard facts that we will occasionally use are that $\Spc(S)$ has a left proper model structure in which monomorphisms are cofibrations and that $\Spc(S)$ and $\Spt(S)$ admit model structures which are left Bousfield localizations of objectwise and levelwise model structures, respectively. The latter implies that an objectwise (\resp{} levelwise) homotopy colimit in $\Spc(S)$ (\resp{} in $\Spt(S)$), by which we mean a homotopy colimit with respect to objectwise (\resp{} levelwise) equivalences, is a homotopy colimit. In particular, filtered colimits are always homotopy colimits in these categories since this is true for simplicial sets.

In Lemma~\ref{lem:combinatorial} we will recall that $\Spt(S)$ admits a model structure which is symmetric monoidal, simplicial, combinatorial, and which satisfies the monoid axiom of \cite[Definition 3.3]{Schwede:2000}. By \cite[Theorem 4.1]{Schwede:2000}, if $E$ is a commutative monoid in $\Spt(S)$, the category $\Mod_E$ of $E$-modules inherits a symmetric monoidal model structure. In particular, we can consider homotopy colimits of $E$-modules (which are also homotopy colimits in the underlying category $\Spt(S)$). We denote by $\D(E)$ the homotopy category of $\Mod_E$. To distinguish these highly structured modules from modules over commutative monoids in the monoidal category $\SH(S)$, we call the latter \emph{weak modules}.

We use the notation $\Map(X,Y)$ for derived mapping spaces, while $[X,Y]=\pi_0\Map(X,Y)$ is the set of morphisms in the homotopy category. As a general rule, when a functor between model categories preserves equivalences, we use the same symbol for the induced functor on homotopy categories. Otherwise we use the prefixes $\L$ and $\R$ to denote left and right derived functors, but here are some exceptions.
The smash product of spectra $E\wedge F$ always denotes the derived monoidal structure on $\SH(S)$. In particular, if $E$ is a commutative monoid in $\Spt(S)$, $E\wedge\ph$ is the left derived functor of the free $E$-module functor $\Spt(S)\to\Mod_E$. The bigraded suspension and loop functors $\Sigma^{p,q}$ and $\Omega^{p,q}$ are also always considered at the level of homotopy categories.

For $E\in\SH(S)$, $\spi_{p,q}(E)$ will denote the Nisnevich sheaf associated with the presheaf
\[X\mapsto [\Sigma^{p,q}\Sigma^\infty X_+,E]\]
on $\Sm/S$.
By \cite[Proposition~5.1.14]{Morel:2003}, the family of functors $\spi_{p,q}$, $p,q\in\Z$, detects equivalences in $\SH(S)$. We note that the functors $\spi_{p,q}$ preserve sums and filtered colimits, because the objects $\Sigma^{p,q}\Sigma^\infty X_+$ are compact (this follows by higher abstract nonsense from the observation that filtered homotopy colimits of simplicial presheaves on $\Sm/S$ preserve $\A^1$-local objects and Nisnevich-local objects, \cf{} Appendix~\ref{sec:appendix}; the proof in \cite[\S9]{Dugger:2005} also works over an arbitrary base scheme).

We do not give here the definition of the motivic Thom spectrum $\MGL$, but we recall that it can be constructed as a commutative monoid in $\Spt(S)$ (\cite[\S2.1]{Panin:2008}), and that it is a cellular spectrum (the unstable cellularity of Grassmannians over any base is proved in \cite[Proposition 3.7]{Wendt:2012}; the proof of \cite[Theorem 6.4]{Dugger:2005} also works without modifications over $\Spec\Z$, which implies the cellularity of $\MGL$ over an arbitrary base).

\subsection{The homotopy \texorpdfstring{$t$}{𝑡}-structure}\label{sub:tstructure}

Let $\SH(S)_{\geq d}$ denote the subcategory of $\SH(S)$ generated under homotopy colimits and extensions by \[\{\Sigma^{p,q}\Sigma^\infty X_+\suchthat X\in\Sm/S\text{ and }p-q\geq d\}.\]
Spectra in $\SH(S)_{\geq d}$ are called \emph{$d$-connective} (or simply \emph{connective} if $d=0$). Since $\Spt(S)$ is a combinatorial simplicial model category, $\SH(S)_{\geq 0}$ is the nonnegative part of a unique $t$-structure on $\SH(S)$ (combine \cite[Proposition A.3.7.6]{HTT} and \cite[Proposition 1.4.4.11]{HA}), called the \emph{homotopy $t$-structure}. The associated truncation functors are denoted by $E\mapsto E_{\geq d}$ and $E\mapsto E_{\leq d}$, so that we have cofiber sequences
\[E_{\geq d}\to E\to E_{\leq d-1}\to\Sigma^{1,0}E_{\geq d}.\]
Write $\kappa_dE$ for the cofiber of $E_{\geq d+1}\to E_{\geq d}$, or equivalently the fiber of $E_{\leq d}\to E_{\leq d-1}$.

The filtration of $\SH(S)$ by the subcategories $\SH(S)_{\geq d}$ adheres to the axiomatic framework of \cite[\S2.1]{Gutierrez:2012}. It follows from \cite[\S2.3]{Gutierrez:2012} that the full slice functor
\[\kappa_\ast\colon\SH(S)\to \SH(S)^{\Z}\]
has a lax symmetric monoidal structure, \ie, there are natural coherent maps
\[\kappa_mE\wedge \kappa_nF\to \kappa_{m+n}(E\wedge F).\]
In particular, $\kappa_\ast$ preserves monoids and modules.

\begin{lemma}\label{lem:filtered}
	The truncation functor $(\ph)_{\leq d}\colon \SH(S)\to\SH(S)$ preserves filtered homotopy colimits.
\end{lemma}

\begin{proof}
	It suffices to show that $\SH(S)_{\leq d}$ is closed under filtered homotopy colimits.	Since
	\[\SH(S)_{\leq d}=\{E\in\SH(S)\suchthat [F,E]=0\text{ for all }F\in\SH(S)_{\geq d+1}\},\]
	this follows from the fact that $\SH(S)_{\geq d+1}$ is generated under homotopy colimits and extensions by compact objects.
\end{proof}

\begin{lemma}\label{lem:basechangetrunc}
	Let $f\colon T\to S$ be an essentially smooth morphism of base schemes, and let $f^\ast\colon\SH(S)\to\SH(T)$ be the induced base change functor. Then, for any $E\in\SH(S)$ and any $d\in\Z$, $f^\ast(E_{\leq d})\simeq (f^\ast E)_{\leq d}$.
\end{lemma}

\begin{proof}
	Let $f_\ast$ be the right adjoint to $f^\ast$. With no assumption on $f$, we have \[f^\ast(\Sigma^{p,q}\Sigma^\infty X_+)\simeq \Sigma^{p,q}\Sigma^\infty(X\times_ST)_+\] for every $X\in\Sm/S$. Since $f^\ast$ preserves homotopy colimits, it follows that $f^\ast(\SH(S)_{\geq d+1})\subset \SH(T)_{\geq d+1}$ and hence, by adjunction, that $f_\ast(\SH(T)_{\leq d})\subset\SH(S)_{\leq d}$. This shows that $f_\ast$ is compatible with the inclusions $\SH_{\leq d}\subset\SH$. Taking left adjoints, we obtain a canonical equivalence $(f^\ast(E_{\leq d}))_{\leq d}\simeq (f^\ast E)_{\leq d}$. Thus, it remains to show that $f^\ast(\SH(S)_{\leq d})\subset\SH(T)_{\leq d}$. By Lemma~\ref{lem:stablecontinuity} \itemref{lem:stablecontinuity:1} we can assume that $f$ is smooth. Then $f^\ast$ admits a left adjoint $f_\sharp$ such that, for $Y\in \Sm/T$, \[f_\sharp(\Sigma^{p,q}\Sigma^\infty Y_+)\simeq \Sigma^{p,q}\Sigma^\infty Y_+.\] As before, this implies $f_\sharp(\SH(T)_{\geq d+1})\subset\SH(S)_{\geq d+1}$ whence the desired result by adjunction.
\end{proof}

Morel's connectivity theorem gives a more explicit description of the homotopy $t$-structure when $S$ is the spectrum of a field:

\begin{samepage}
\begin{theorem}\label{thm:t-structure}
    Let $k$ be a field and let $E\in\SH(k)$. Then
    \begin{enumerate}
        \item\label{thm:t-structure:1} $E\in\SH(k)_{\geq d}$ if and only if $\spi_{p,q}E=0$ for $p-q<d$;
        \item\label{thm:t-structure:2} $E\in\SH(k)_{\leq d}$ if and only if $\spi_{p,q}E=0$ for $p-q>d$.
    \end{enumerate}
\end{theorem}
\end{samepage}

\begin{proof}
	We first observe that the vanishing condition in \itemref{thm:t-structure:2} (say for $d=-1$) implies in fact the vanishing of the individual groups $[\Sigma^{p,q}\Sigma^\infty X_+,E]$ for $p-q\geq 0$ and $X\in\Sm/k$. By the standard adjunctions, this group is equal to the set of maps $\Sigma^{p-q}X_+\to L_{\A^1}\Omega^\infty\Sigma^{-q,-q}E$ in the homotopy category of pointed simplicial sheaves.
	Thus, the vanishing of the sheaves for all $p-q\geq 0$ implies that $L_{\A^1}\Omega^\infty\Sigma^{-q,-q}E$ is contractible, whence the vanishing of the presheaves.
	 
	By \cite[\S5.2]{Morel:2003}, the right-hand sides of \itemref{thm:t-structure:1} and \itemref{thm:t-structure:2} define a $t$-structure on $\SH(k)$;\footnote{This is only proved when $k$ is perfect in \cite{Morel:2003}, but that hypothesis is removed in \cite{Morel:2005}. In any case we can assume that $k$ is perfect in all our applications of this theorem.} call it $\scr T$. To show that this $t$-structure coincides with ours, it suffices to show the implications from left to right in \itemref{thm:t-structure:1} and \itemref{thm:t-structure:2}.
	For \itemref{thm:t-structure:1}, we have to show that if $F\in\SH(k)_{\geq 0}$, then $F$ is $\scr T$-nonnegative, or equivalently $[F,E]=0$ for every $\scr T$-negative $E$. Now $\scr T$-nonnegative spectra are easily seen to be closed under homotopy colimits and extensions, so we may assume that $F=\Sigma^{p,q}\Sigma^\infty X_+$ with $p-q\geq 0$. But then $[F,E]=0$ by our preliminary observation.
	For \itemref{thm:t-structure:2}, let $E\in\SH(k)_{\leq -1}$, \ie, $[F,E]=0$ for all $F\in\SH(k)_{\geq 0}$. Taking $F=\Sigma^{p,q}\Sigma^\infty X_+$ with $p-q\geq 0$, we deduce that $E$ is $\scr T$-negative.
\end{proof}

\begin{corollary}\label{cor:nondeg}
Let $k$ be a field, $X\in\Sm/k$, and $p,q\in\Z$. For every $E\in\SH(k)$ and $d>p-q+\dim X$, \[[\Sigma^{p,q}\Sigma^\infty X_+,E_{\geq d}]=0.\] In particular, the canonical map $E\to\holim_{n\to\infty} E_{\leq n}$ is an equivalence, \ie, the homotopy $t$-structure on $\SH(k)$ is left complete.
\end{corollary}

\begin{proof}
	By Theorem~\ref{thm:t-structure} \itemref{thm:t-structure:1}, the Nisnevich-local presheaf of spectra $\Omega^\infty_{\G_m}\Omega^{p,q}(E_{\geq d})$ is $(d-p+q)$-connective. Since the Nisnevich cohomological dimension of $X$ is at most $\dim X$, this implies the first claim. It follows that $\holim_{n\to\infty}E_{\geq n}=0$, whence the second claim.
\end{proof}

\begin{remark}
	It is clear that the homotopy $t$-structure is right complete over any base scheme.
\end{remark}

\begin{remark}
	Theorem~\ref{thm:t-structure} and Corollary~\ref{cor:nondeg} are true more generally over any base scheme satisfying the stable $\A^1$-connectivity property in the sense of \cite[Definition 1]{Morel:2005}. Our definition of the homotopy $t$-structure is thus a conservative extension of Morel's definition to all base schemes, and it allows us for example to state the results of \S\ref{sec:grassmannians} unconditionally over any base scheme. However, the main result of this paper only uses the homotopy $t$-structure when $S$ is the spectrum of a perfect field, so this generality is merely a convenience.
\end{remark}

\subsection{Strictly \texorpdfstring{$\A^1$}{𝐀¹}-invariant sheaves}
\label{sub:strictlyA1inv}

In this paragraph we recall the following fact, proved by Morel: if $k$ is a perfect field, equivalences in $\SH(k)$ are detected by the stalks of the sheaves $\spi_{p,q}$ at generic points of smooth schemes.

If $\scr F$ is a presheaf on $\Sm/S$ and $(X_\alpha)$ is a cofiltered diagram in $\Sm/S$ with affine transition maps, we can define the value of $\scr F$ at $X=\lim_\alpha X_\alpha$ by the usual formula
\[\scr F(X)=\colim_\alpha\scr F(X_\alpha).\]
This is well-defined: in fact, we have $\scr F(X)\cong\Hom(rX,\scr F)$ where $rX$ is the presheaf on $\Sm/S$ represented by $X$ (\cite[Proposition 8.14.2]{EGA4-3}).
In particular, if $S$ is the spectrum of a field $k$ and $L$ is a separably generated extension of $k$, $\scr F(\Spec L)$ is defined in this way. Note that if $\scr F'$ is the Nisnevich sheaf associated with $\scr F$, then $\scr F(\Spec L)\cong\scr F'(\Spec L)$.

A Nisnevich sheaf of abelian groups $\scr F$ on $\Sm/S$ is called \emph{strictly $\A^1$-invariant} if the map
\[H^i_\Nis(X,\scr F)\to H^i_\Nis(X\times\A^1,\scr F)\]
induced by the projection $X\times\A^1\to X$ is an isomorphism for every $X\in\Sm/S$ and $i\geq 0$.
When $S$ is the spectrum of a field (or, more generally, when the stable $\A^1$-connectivity property holds over $S$), the category of strictly $\A^1$-invariant sheaves is an exact abelian subcategory of the category of all abelian sheaves (\cite[Corollary~6.24]{Morel:2012}). In this case a typical example of a strictly $\A^1$-invariant sheaf is $\spi_{p,q}E$ for a spectrum $E\in\SH(S)$ and $p,q\in\Z$ (\cite[Remark~5.1.13]{Morel:2003}).

\begin{theorem}\label{thm:fields}
    Let $f\colon\scr F\to\scr G$ be a morphism of strictly $\A^1$-invariant Nisnevich sheaves on $\Sm/k$ where $k$ is a \emph{perfect} field. Then $f$ is an isomorphism if and only if, for every finitely generated field extension $k\subset L$, the map $\scr F(\Spec L)\to\scr G(\Spec L)$ induced by $f$ is an isomorphism.
\end{theorem}

\begin{proof}
	This follows from \cite[Example 2.3 and Proposition 2.8]{Morel:2012}.
\end{proof}

\section{The stable path components of \texorpdfstring{$\MGL$}{MGL}}
\label{sec:grassmannians}

In this section we compute $\kappa_0\MGL$ as a weak $\kappa_0\1$-module over an arbitrary base scheme $S$.
If $X\in\Sm/S$ and $E\to X$ is a vector bundle, we denote its Thom space by
\[\Th(E)=E/(E\minus X).\]

\begin{lemma}\label{lem:thomconnected}
    Let $E$ be a rank $d$ vector bundle on $X\in\Sm/S$. Then $\Sigma^\infty\Th(E)$ is $d$-connective.
\end{lemma}

\begin{proof}
    Let $\{U_\alpha\}$ be a trivializing Zariski cover of $X$. Then $\Th(E)$ is the homotopy colimit of the simplicial diagram
\begin{tikzmath}
	\def\colsep{.85em}
	\diagram{
	\dotsb & \bigvee_{\alpha,\beta}\Th(E|_{U_{\alpha\beta}}) & \bigvee_\alpha\Th(E|_{U_\alpha})\rlap. \\
	};
	\arrows
	(11-) edge[-top,vshift=2*\dbl] (-12) edge[-mid] (-12) edge[-bot,vshift=-2*\dbl] (-12)
	(12-) edge[-top,vshift=\dbl] (-13) edge[-bot,vshift=-\dbl] (-13);
\end{tikzmath}
Since $E$ is trivial on $U_{\alpha_1\dotso\alpha_n}$, $\Th(E|_{U_{\alpha_1\dotso\alpha_n}})$ is equivalent to $\Sigma^{2d,d}(U_{\alpha_1\dotso\alpha_n})_+$ which is stably $d$-connective.
    Thus, $\Sigma^\infty\Th(E)$ is $d$-connective as a homotopy colimit of $d$-connective spectra.
\end{proof}

\begin{lemma}\label{lem:connected}
    Let $X\in\Sm/S$ and let $Z\subset X$ be a smooth closed subscheme of codimension $d$. Then $\Sigma^\infty(X/(X\minus Z))$ is $d$-connective.
\end{lemma}

\begin{proof}
    This follows from Lemma~\ref{lem:thomconnected} and the equivalence $X/(X\minus Z)\simeq\Th(\scr N)$ of \cite[Theorem 3.2.23]{Morel:1999}, where $\scr N$ is the normal bundle of $Z\subset X$.
\end{proof}

Let $V$ be a vector bundle of finite rank over $S$. Denote by $\Gr(r,V)$ the Grassmann scheme of $r$-planes in $V$, and let $E(r,V)$ denote the tautological rank $r$ vector bundle over it. When $V=\A^n$ we will also write $\Gr(r,n)$ and $E(r,n)$. A subbundle $W\subset V$ induces a closed immersion
\[i_W\colon\Gr(r,W)\into\Gr(r,V),\]
and given vector bundles $V_1$, \dots, $V_t$, there is a closed immersion
\[j_{V_1,\dotsc,V_t}\colon\Gr(r_1,V_1)\times\dotsb\times\Gr(r_t,V_t)\into\Gr(r_1+\dotsb+r_t,V_1\times\dotsb\times V_t).\]
It is clear that these immersions are compatible with the tautological bundles as follows:
\begin{gather*}
    i_W^\ast E(r,V)\cong E(r,W)\text{ and}\\
    j_{V_1,\dotsc,V_t}^\ast E(r_1+\dotsb+r_t,V_1\times\dotsb\times V_t)\cong E(r_1,V_1)\times\dotsb\times E(r_t,V_t).
\end{gather*}

Suppose now that $W\subset V\cong\A^n$ is a hyperplane with complementary line $L\subset V$ (so that $\Gr(1,L)=S$). Then the closed immersions
\[i_W\colon\Gr(r,W)\into\Gr(r,V)\quad\text{and}\quad j_{L,W}\colon\Gr(r-1,W)\into\Gr(r,V)\]
have disjoint images and are complementary in the following sense. The inclusion \[i_W\colon\Gr(r,W)\into\Gr(r,V)\minus \Im(j_{L,W})\] is the zero section of a rank $r$ vector bundle \[p\colon \Gr(r,V)\minus \Im(j_{L,W})\to\Gr(r,W)\] whose fiber over an $S$-point $P\in\Gr(r,W)$ is the vector bundle of $r$-planes in $P\oplus L$ not containing $L$. Similarly, the inclusion \[j_{L,W}\colon\Gr(r-1,W)\into\Gr(r,V)\minus \Im(i_W)\] is the zero section of a rank $n-r$ vector bundle
\[q\colon\Gr(r,V)\minus \Im(i_W)\to \Gr(r-1,W)\]
whose fiber over $P\in\Gr(r-1,W)$ can be identified with the vector bundle of lines in a complement of $P$ that are not contained in $W$.

\begin{lemma}\label{lem:grass0}
    The immersion $i_W$ induces an equivalence \[\Th(E(r,W))\simeq\Th(E(r,V)|_{\Gr(r,V)\minus \Im(j_{L,W})}),\] and the immersion $j_{L,W}$ an equivalence \[\Th(E(1,L)\times E(r-1,W))\simeq\Th(E(r,V)|_{\Gr(r,V)\minus \Im(i_W)}).\]
\end{lemma}

\begin{proof}
    The vector bundles $E(r,W)$ and $E(r,V)|_{\Gr(r,V)\minus \Im(j_{L,W})}$ are pullbacks of one another along the immersion $i_W$ and its retraction $p$. It follows that they are strictly $\A^1$-homotopy equivalent in the category of vector bundles and fiberwise isomorphisms. In particular, the complements of their zero sections are $\A^1$-homotopy equivalent, and therefore their Thom spaces are equivalent. The second statement is proved in the same way.
\end{proof}

\begin{samepage}
\begin{lemma}\label{lem:grass}
    Let $W\subset V\cong\A^n$ be a hyperplane with complementary line $L$. Then
    \begin{enumerate}
        \item\label{lem:grass:1} the cofiber of $i_W\colon\Gr(r,W)\into\Gr(r,V)$ is stably $(n-r)$-connective;
        \item\label{lem:grass:2} the cofiber of $j_{L,W}\colon\Gr(r-1,W)\into\Gr(r,V)$ is stably $r$-connective;
        \item\label{lem:grass:3} the cofiber of $\Th(E(r,W))\into\Th(E(r,V))$ is stably $n$-connective;
        \item\label{lem:grass:4} the cofiber of $\Th(E(1,L)\times E(r-1,W))\into\Th(E(r,V))$ is stably $2r$-connective.
    \end{enumerate}
\end{lemma}
\end{samepage}

\begin{proof}
    By the preceding discussion, the cofiber of $i_W$ is equivalent to
    \[\Gr(r,V)/(\Gr(r,V)\minus \Im(j_{L,W}))\]
    which is stably $(n-r)$-connective by Lemma~\ref{lem:connected}. The proof of \itemref{lem:grass:2} is identical.
    
    By Lemma~\ref{lem:grass0}, the cofiber in \itemref{lem:grass:3} is equivalent to \[\Th(E(r,V))/\Th(E(r,V)|_{\Gr(r,V)\minus \Im(j_{L,W})}).\] This quotient is isomorphic to \[ E(r,V)/(E(r,V)\minus \Im(j_{L,W}))\] which is stably $n$-connective by Lemma~\ref{lem:connected}. The proof of \itemref{lem:grass:4} is identical.
\end{proof}

The Hopf map is the projection $h\colon\A^2\minus \{0\}\to\P^1$; let $C(h)$ be its cofiber.
The commutative diagram
\begin{tikzmath}
	\diagram{
	\A^2\minus \{0\} & \P^1 \\
	E(1,2)\minus \P^1 & E(1,2) \\
	};
	\arrows
	(11-) edge node[above]{$h$} (-12)
	(21-) edge[c->] (-22)
	(21) edge node[left]{$\cong$} (11)
	(22) edge node[right]{$\simeq$} (12);
\end{tikzmath}
(where $E(1,2)$ is the tautological bundle on $\P^1$) shows that $C(h)\simeq\Th(E(1,2))$. Thus, we have a canonical map $C(h)\to\MGL_1=\colim_{n\to\infty}\Th(E(1,n))$. Using the bonding maps of the spectrum $\MGL$, we obtain maps
\begin{equation}\label{eqn:Ch}
\Sigma^{2r-2,r-1}C(h)\to \MGL_r
\end{equation}
for every $r\geq 1$, and in the colimit we obtain a map
\begin{equation}\label{eqn:Chstable}
\Sigma^{-2,-1}\Sigma^\infty C(h)\to \hocolim_{r\to\infty}\Sigma^{-2r,-r}\Sigma^\infty\MGL_r\simeq \MGL.
\end{equation}

\begin{lemma}\label{lem:pikkMGL}
	 The map~\eqref{eqn:Chstable} induces an equivalence
	 $(\Sigma^{-2,-1}\Sigma^\infty C(h))_{\leq 0}\simeq\MGL_{\leq 0}$.
\end{lemma}

\begin{proof}
	By definition of~\eqref{eqn:Chstable} and Lemma~\ref{lem:filtered}, it suffices to show that~\eqref{eqn:Ch} induces an equivalence
 \[(\Sigma^\infty\Sigma^{2r-2,r-1}C(h))_{\leq r}\simeq(\Sigma^\infty\MGL_r)_{\leq r},\]
 for every $r\geq 1$. Consider in $\H_\pt(S)$ the commutative diagram
 \begin{tikzmath}
	\def\colsep{2em}
	\def\rowsep{2em}
	\def\vd{\raisebox{0pt}[1.1em][.2em]{\vdots}}
 	\diagram{
	\Sigma^{2r-2,r-1}\Th(E(1,2)) & \Sigma^{2r-4,r-2}\Th(E(2,3)) & \dotsb & \Th(E(r,r+1)) \\
	\Sigma^{2r-2,r-1}\Th(E(1,3)) & \Sigma^{2r-4,r-2}\Th(E(2,4)) & \dotsb & \Th(E(r,r+2)) \\
	\vd & \vd & \ddots & \vd \\
	\Sigma^{2r-2,r-1}\MGL_1 & \Sigma^{2r-4,r-2}\MGL_2 & \dotsb & \MGL_r \\
	};
	\arrows (11-) edge (-12) (12-) edge (-13) (13-) edge (-14)
	(21-) edge (-22) (22-) edge (-23) (23-) edge (-24)
	(41-) edge (-42) (42-) edge (-43) (43-) edge (-44)
	(11) edge (21) (21) edge (31) (31) edge (41)
	(12) edge (22) (22) edge (32) (32) edge (42)
	(14) edge (24) (24) edge (34) (34) edge (44);
 \end{tikzmath}
 with $r$ columns and $\omega+1$ rows. The vertical maps are induced by the closed immersions $i$ and the horizontal maps are induced by the closed immersions $j$; the last row is the colimit of the previous rows and the maps in the last row are the bonding maps defining the spectrum $\MGL$. Now the map~\eqref{eqn:Ch} is obtained by traveling down the left column and then along the bottom row, and hence it is also the composition of the top row followed by the right column. Explicitly, the top row is composed of the maps
	\[
		\Th(E(1,1)^{\times (r-s+1)}\times E(s-1,{s}))\to \Th(E(1,1)^{\times (r-s)}\times E(s,{s+1}))
	\]
	induced by $j_{\A^1,\A^s}$	for $2\leq s\leq r$, and the right column is composed of the maps
	\[
		\Th(E(r,{n}))\to\Th(E(r,{n+1}))
	\]
	induced by $i_{\A^n}$ for $n\geq r+1$. The former have stably $(r+s)$-connective cofiber by Lemma~\ref{lem:grass} \itemref{lem:grass:4}, and the latter have stably $(n+1)$-connective cofiber by Lemma~\ref{lem:grass} \itemref{lem:grass:3}. Since $r+s\geq r+1$ and $n+1\geq r+1$, all those maps become equivalences in $\SH(S)_{\leq r}$. By Lemma~\ref{lem:filtered}, their composition also becomes an equivalence in $\SH(S)_{\leq r}$, as was to be shown.
\end{proof}

Recall that there are canonical equivalences $\A^2\minus\{0\}\simeq S^{3,2}$ and $\P^1\simeq S^{2,1}$ in $\H_\pt(S)$, and hence $h$ stabilizes to a map 
\[\eta\colon\Sigma^{1,1}\1\to\1.\]

\begin{theorem}\label{thm:pikkMGL}
    The unit $\1\to\MGL$ induces an equivalence $(\1/\eta)_{\leq 0}\simeq\MGL_{\leq 0}$.
\end{theorem}

\begin{proof}
	Follows from Lemma~\ref{lem:pikkMGL} and the easy fact that the composition \[\1\to\1/\eta=\Sigma^{-2,-1}\Sigma^\infty C(h)\to\MGL\] is the unit of $\MGL$.
\end{proof}

The following corollary is also proved using different arguments in \cite[Theorem 5.7]{Spitzweck2:2010}.

\begin{corollary}\label{cor:MGLeffective}
	The spectrum $\MGL$ is connective.
\end{corollary}

\begin{proof}
	Follows from Theorem~\ref{thm:pikkMGL} since $\1/\eta$ is obviously connective.
\end{proof}

\begin{remark}
	Suppose that $S$ is essentially smooth over a field. Combined with the computation of $\spi_{n,n}(\1)$ over perfect fields from \cite[Remark~6.42]{Morel:2012}, Theorem~\ref{thm:t-structure}, and Lemma~\ref{lem:stablecontinuity} \itemref{lem:stablecontinuity:1}, Theorem~\ref{thm:pikkMGL} shows that for any $X\in\Sm/S$, $\spi_{-n,-n}(\MGL)(X)$ is the $n$th unramified Milnor $K$-theory group of $X$.
\end{remark}

\section{Complements on motivic cohomology}

\subsection{Spaces and spectra with transfers}
\label{sub:transfers}

Let $S$ be a base scheme and let $R$ be a commutative ring. We begin by recalling the existence of a commutative diagram of symmetric monoidal Quillen adjunctions:
\pagebreak[0]
\begin{tikzequation}\label{eqn:transfers}
	\def\colsep{4em}
	\diagram{
	\Sm/S &[between origins] &[between origins] \Cor(S,R) \\
	\Spc_\pt(S) & & \Spc_\tr(S,R) \\
	\Spt(S) & & \Spt_\tr(S,R) \\
	& \Mod_{HR}\rlap. & \\
	};
	\arrows
	(11-) edge node[above]{$\Gamma$} (-13)
	(11) edge (21) (13) edge (23)
	(21) edge[vshift=\dbl] node[above=\dbl]{$R_\tr$} (23) edge[<-,vshift=-\dbl] node[below=\dbl]{$u_\tr$} (23)
	(31) edge[vshift=\dbl] node[above=\dbl]{$R_\tr$} (33) edge[<-,vshift=-\dbl] node[below=\dbl]{$u_\tr$} (33)
	(21) edge[vshift=-\dbl] node[left=\dbl]{$\Sigma^\infty$} (31) edge[<-,vshift=\dbl] node[right=\dbl]{$\Omega^\infty$} (31)
	(23) edge[vshift=-\dbl] node[left=\dbl]{$\Sigma^\infty_\tr$} (33) edge[<-,vshift=\dbl] node[right=\dbl]{$\Omega^\infty_\tr$} (33)
	(42) edge[<-,vshift=\dbl] node[below left]{$HR\wedge\ph$} (31) edge[vshift=-\dbl] (31) edge[<-,vshift=-\dbl] node[below right]{$\Psi$} (33) edge[vshift=\dbl] node[above left]{$\Phi$} (33);
\end{tikzequation}
For a more detailed description and for proofs we refer to \cite[\S2]{Rondigs:2008} (where $\Z$ can harmlessly be replaced with $R$). The homotopy categories of $\Spc_\tr(S,R)$ and $\Spt_\tr(S,R)$ will be denoted by $\H_\tr(S,R)$ and $\SH_\tr(S,R)$, respectively. For the purpose of this diagram, the model structure on $\Spc_\pt(S)$ is the projective $\A^1$-Nisnevich-local structure, \ie, it is the left Bousfield localization of the projective objectwise model structure on simplicial presheaves. Similarly, $\Spt(S)$ is endowed with the projective stable model structure, \ie, the left Bousfield localization of the projective levelwise model structure. These model structures are all simplicial. We record the following result for which we could not find a complete reference.

\begin{lemma}\label{lem:combinatorial}
	The projective stable model structure on $\Spt(S)$ is combinatorial and it satisfies the monoid axiom of \cite[Definition 3.3]{Schwede:2000}.
\end{lemma}

\begin{proof}
	The stable projective model structure is the left Bousfield localization at a small set of maps (described in \cite[Definition 8.7]{Hovey:2001}) of the projective levelwise model structure which is left proper and cellular by \cite[Theorem 8.2]{Hovey:2001}. It follows from \cite[Theorem 4.1.1]{Hirschhorn:2009} that the projective stable model structure is cellular. By definition, symmetric $S^{2,1}$-spectra are algebras over a colimit-preserving monad on the category of symmetric sequences of pointed simplicial presheaves on $\Sm/S$. The latter category is merely a category of presheaves of pointed sets on an essentially small category and hence is locally presentable. We deduce from \cite[Remark 2.75]{Adamek:1994} that $\Spt(S)$ is locally presentable and therefore combinatorial.
	
	It remains to check the monoid axiom. Let $i$ be an acyclic cofibration and let $X\in\Spt(S)$. It is obvious that the identity functor is a left Quillen equivalence from our model structure to the Jardine model structure of \cite[Theorem 4.15]{Jardine:2000}, so that $i$ is also an acyclic cofibration for the Jardine model structure. By \cite[Proposition 4.19]{Jardine:2000}, $i\wedge X$ is an equivalence and a levelwise monomorphism. We conclude by showing that the class of maps that are simultaneously equivalences and levelwise monomorphisms is stable under cobase change and transfinite composition. Any pushout along a levelwise monomorphism is in fact a levelwise homotopy pushout since there exists a left proper levelwise model structure on $\Spt(S)$ in which levelwise monomorphisms are cofibrations (\cite[Theorem 4.2]{Jardine:2000}), and hence it is a homotopy pushout. The stability under transfinite composition follows from the fact that filtered colimits in $\Spt(S)$ are homotopy colimits.
\end{proof}

The category $\Cor(S,R)$ of finite correspondences with coefficients in $R$ has as objects the separated smooth schemes of finite type over $S$. The set of morphisms from $X$ to $Y$ will be denoted by $\cor_R(X,Y)$; in the notation and terminology of \cite[\S 9.1.1]{Cisinski:2012}, it is the $R$-module $c_0(X\times_SY/X,\Z)\tens R$, where $c_0(X/S,\Z)$ is the group of finite $\Z$-universal $S$-cycles with domain $X$. Then $\Cor(S,R)$ is an additive category, with direct sum given by disjoint union of schemes. It also has a monoidal structure given by direct product on objects such that the graph functor $\Gamma\colon\Sm/S\to\Cor(S,R)$ is symmetric monoidal.

$\Spc_\tr(S,R)$ is the category of additive simplicial presheaves on $\Cor(S,R)$, endowed with the projective $\A^1$-Nisnevich-local model structure. The functor
\[R_\tr\colon\Spc_\pt(S)\to\Spc_\tr(S,R)\]
is the unique colimit-preserving simplicial functor such that, for all $X\in \Sm/S$, $R_\tr X_+$ is the presheaf on $\Cor(S,R)$ represented by $X$. Its right adjoint $u_\tr$ is restriction along $\Gamma$ and it detects equivalences. The tensor product $\tens_R$ on $\Spc_\tr(S,R)$ is the composition of the ``external'' cartesian product with the additive left Kan extension of the monoidal product on $\Cor(S,R)$ (see the proof of Lemma~\ref{lem:balanced} for a formulaic version of this definition). Since the cartesian product on $\Spc(S)$ is obtained from that of $\Sm/S$ in the same way, $R_\tr$ has a symmetric monoidal structure. 

$\Spt_\tr(S,R)$ is the category of symmetric $R_\tr S^{2,1}$-spectra in $\Spc_\tr(S,R)$ with the projective stable model structure, and the stable functors $R_\tr$ and $u_\tr$ are defined levelwise. The Eilenberg–Mac Lane spectrum $HR$ is by definition the monoid $u_\tr R_\tr\1$. This immediately yields the monoidal adjunction $(\Phi,\Psi)$ between $\Spt_\tr(S,R)$ and $\Mod_{HR}$, which completes our description of diagram~\eqref{eqn:transfers}.

We make the following observation which is lacking from \cite[\S2]{Rondigs:2008}.

\begin{lemma}\label{lem:utrstable}
	The functor $u_\tr\colon\Spt_\tr(S,R)\to\Spt(S)$ detects equivalences.
\end{lemma}

\begin{proof}
	It detects levelwise equivalences since $u_\tr\colon\Spc_\tr(S,R)\to\Spc_\pt(S)$ detects equivalences. Define a functor $Q\colon\Spt(S)\to\Spt(S)$ by $(QE)_n=\Hom(S^{2,1}, E_{n+1})$ (with action of $\Sigma_n$ induced by that of $\Sigma_{n+1}$), and let $Q^\infty E=\colim_{n\to\infty}Q^nE$. Similarly, let $Q_\tr\colon\Spt_\tr(S,R)\to\Spt_\tr(S,R)$ be given by $(Q_\tr E)_n=\Hom(R_\tr S^{2,1}, E_{n+1})$. Then a morphism $f$ in $\Spt(S)$ (\resp{} in $\Spt_\tr(S,R)$) is a stable equivalence if and only if $Q^\infty(f)$ (\resp{} $Q^\infty_\tr(f)$) is a levelwise equivalence. The proof is completed by noting that $u_\tr Q^\infty_\tr\cong Q^\infty u_\tr$.
\end{proof}

An $HR$-module is called \emph{cellular} if it is an iterated homotopy colimit of $HR$-modules of the form $\Sigma^{p,q}HR$ with $p,q\in\Z$. Similarly, an object in $\Spt_\tr(S,R)$ is cellular if it is an iterated homotopy colimit of objects of the form $R_\tr\Sigma^{p,q}\1$ with $p,q\in\Z$.

\begin{lemma}\label{lem:HZmod}
	The derived adjunction $(\L\Phi,\R\Psi)$ between $\D(HR)$ and $\SH_\tr(S,R)$ restricts to an equivalence between the full subcategories of cellular objects.
\end{lemma}

\begin{proof}
	The proof of \cite[Corollary 62]{Rondigs:2008} works with any ring $R$ instead of $\Z$.
\end{proof}

\subsection{Eilenberg–Mac Lane spaces and spectra}
\label{sub:EML}

Denote by $\s\Mod_R$ the category of simplicial $R$-modules with its usual model structure. Then there is a symmetric monoidal Quillen adjunction
\begin{tikzmath}
	\diagram{\s\Mod_R & \Spc_\tr(S,R) \\};
	\arrows (11-) edge[vshift=\dbl] node[above=\dbl]{$c_\tr$} (-12) (-12) edge[vshift=\dbl] (11-);
\end{tikzmath}
where $c_\tr A$ is the ``constant additive presheaf'' with value $A$. It is clear that $c_\tr$ preserves equivalences. It is easy to show that this adjunction extends to a symmetric monoidal Quillen adjunction
\begin{tikzmath}
	\diagram{\Sp(\s\Mod_R) & \Spt_\tr(S,R) \\};
	\arrows (11-) edge[vshift=\dbl] node[above=\dbl]{$c_\tr$} (-12) (-12) edge[vshift=\dbl] (11-);
\end{tikzmath}
where $\Sp(\s\Mod_R)$ is the category of symmetric spectra in $\s\Mod_R$ with the projective stable model structure. Explicitly, the stable functor $c_\tr$ is given by the formula
\[(A_0,A_1,A_2,\dotsc)\mapsto(c_\tr A_0,R_\tr \tilde\G_m\tens_R c_\tr A_1, R_\tr(\tilde\G_m^{\wedge 2})\tens_R c_\tr A_2,\dotsc).\]
Of course, $\Sp(\s\Mod_R)$ is Quillen equivalent to the model category of unbounded chain complexes of $R$-modules.

If $p\geq q\geq 0$ and $A\in\s\Mod_R$, the \emph{motivic Eilenberg–Mac Lane space} of degree $p$ and of weight $q$ with coefficients in $A$ is defined by
\[K(A(q),p)=u_\tr(R_\tr S^{p,q}\tens_R c_\tr A).\]
Note that $HR=u_\tr c_\tr(\Sigma^\infty R)$, where $R$ is viewed as a constant simplicial $R$-module. More generally, for any $A\in\Sp(\s\Mod_R)$, we define the \emph{motivic Eilenberg–Mac Lane spectrum} with coefficients in $A$ by
\[HA=u_\tr c_\tr A.\]
Since $u_\tr$ is lax monoidal, $HA$ is a module over the monoid $HR$. If $A\in\s\Mod_R$, the symmetric spectrum $H(\Sigma^\infty A)$ is given by the sequence of motivic spaces $K(A(n),2n)$ for $n\geq 0$. It is clear that the space $K(A(q),p)$ and the spectrum $HA$ do not depend on the ring $R$.

Note that the functors $R_\tr$ and $\tens_R$ do \emph{not} preserve equivalences and that we did not derive them in our definitions of $K(A(q),p)$ and $HA$. We will now justify these definitions by showing that the canonical maps
\begin{gather}
	\label{eqn:EMLderived}
	u_\tr(\L R_\tr S^{p,q}\tens_R^\L c_\tr A)\to K(A(q),p)\text{ and}\\
	\label{eqn:stableEMLderived}
	u_\tr\L c_\tr A\to HA
\end{gather}
in $\H_\pt(S)$ and $\SH(S)$, respectively, are equivalences. In particular, if $A$ is an abelian group, our definition of $K(A(q),p)$ agrees with the definition in \cite[\S3.2]{MEMS} (which in our notations is $u_\tr(\L\Z_\tr S^{p,q}\tens_\Z^\L c_\tr A)$).

\begin{remark}
Our definitions directly realize $HR$ and $HA$ as commutative monoids and modules in the symmetric monoidal model category $\Spt(S)$. For the purposes of this paper, however, all we need to know is that $HR$ is an $\mathbb E_\infty$-algebra in the underlying symmetric monoidal $(\infty,1)$-category and that $HA$ is a module over it; the reader who is familiar with these notions can take the left-hand sides of~\eqref{eqn:EMLderived} and~\eqref{eqn:stableEMLderived} as definitions of $K(A(q),p)$ and $HA$ and skip the proof of the strictification result (which ends with Proposition~\ref{prop:technical}).
\end{remark}

\begin{lemma}\label{lem:freeaction}
    Let $G$ be a group object acting freely on an object $X$ in the category $\Spc_\pt(S)$ or $\Spc_\tr(S,R)$. Then the quotient $X/G$ is the homotopy colimit of the simplicial diagram
\begin{tikzmath}
	\def\colsep{.85em}
	\diagram{
	\dotsb & G\times G\times X & G\times X & X\rlap. \\
	};
	\arrows
	(11-) edge[-top,vshift=3*\dbl] (-12) edge[-mid,vshift=\dbl] (-12) edge[-mid,vshift=-\dbl] (-12) edge[-bot,vshift=-3*\dbl] (-12)
	(12-) edge[-top,vshift=2*\dbl] (-13) edge[-mid] (-13) edge[-bot,vshift=-2*\dbl] (-13)
	(13-) edge[-top,vshift=\dbl] (-14) edge[-bot,vshift=-\dbl] (-14);
\end{tikzmath}
\end{lemma}

\begin{proof}
    This is true for simplicial sets and hence is true objectwise. Objectwise homotopy colimits are homotopy colimits since there exist model structures on $\Spc_\pt(S)$ and $\Spc_\tr(S,R)$ which are left Bousfield localizations of objectwise model structures.
\end{proof}

\begin{lemma}\label{lem:ow9e0}
    Let $Y\subset X$ be an inclusion of objects in $\Spc_\tr(S,R)$. Then
	 \begin{tikzmath}
		\def\colsep{2.2em}
		\diagram{Y & X \\ 0 & X/Y \\};
		\arrows (11-) edge[c->] (-12) (21-) edge (-22) (11) edge (21) (12) edge (22);
	 \end{tikzmath}
	  is a homotopy pushout square.
\end{lemma}

\begin{proof}
    Let $\mathbf 2$ denote the one-arrow category and let $Q\colon\Spc_\tr(S,R)^{\mathbf 2}\to\Spc_\tr(S,R)$ denote the functor $(Y\to X)\mapsto X/Y$. This functor preserves equivalences between cofibrant objects for the projective structure on $\Spc_\tr(S,R)^{\mathbf 2}$, and our claim is that the canonical map $\L Q(Y\to X)\to Q(Y\to X)$ is an equivalence when $Y\to X$ is an inclusion. Factor $Q$ as
\begin{tikzmath}
	\diagram{
	\Spc_\tr(S,R)^{\mathbf 2} & \Spc_\tr(S,R)^{\Delta^\op} & \Spc_\tr(S,R)\rlap, \\
	};
	\arrows (11-) edge node[above]{$B$} (-12) (12-) edge node[above]{$\colim$} (-13);
\end{tikzmath}
    where $B$ sends $Y\to X$ to the simplicial diagram
\begin{tikzmath}
	\def\colsep{.85em}
	\diagram{
	\dotsb & Y\oplus Y\oplus X & Y\oplus X & X\rlap. \\
	};
	\arrows
	(11-) edge[-top,vshift=3*\dbl] (-12) edge[-mid,vshift=\dbl] (-12) edge[-mid,vshift=-\dbl] (-12) edge[-bot,vshift=-3*\dbl] (-12)
	(12-) edge[-top,vshift=2*\dbl] (-13) edge[-mid] (-13) edge[-bot,vshift=-2*\dbl] (-13)
	(13-) edge[-top,vshift=\dbl] (-14) edge[-bot,vshift=-\dbl] (-14);
\end{tikzmath}
    Since $\Spc_\tr(S,R)$ is left proper, $B$ preserves equivalences.
	 Assume now that $Y\to X$ is an inclusion.
	 Then by Lemma~\ref{lem:freeaction}, the canonical map
	 \[\hocolim B(Y\to X)\to \colim B(Y\to X)=Q(Y\to X)\]
	 is an equivalence. Let $\tilde Y\to \tilde X$ be a cofibrant replacement of $Y\to X$, so that $\L Q(Y\to X)\simeq Q(\tilde Y\to\tilde X)$. Then $\tilde Y\to\tilde X$ is a cofibration in $\Spc_\tr(S,R)$ and hence an inclusion. In the commutative square
	 \begin{tikzmath}
	 	\diagram{\hocolim B(\tilde Y\to\tilde X) & Q(\tilde Y\to\tilde X) \\
		\hocolim B(Y\to X) & Q(Y\to X)\rlap, \\};
		\arrows (11-) edge node[above]{$\simeq$} (-12) (21-) edge node[above]{$\simeq$} (-22) (11) edge node[right]{$\simeq$} (21) (12) edge (22);
	 \end{tikzmath}
	 the top, bottom, and left arrows are equivalences, and hence the right arrow is an equivalence as well, as was to be shown.
\end{proof}

An object in $\Spc_\pt(S)$ (\resp{} $\Spc_\tr(S,R)$) will be called \emph{flat} if in each degree it is a filtered colimit of finite sums of objects of the form $X_+$ (\resp{} $R_\tr X_+$) for $X\in\Sm/S$. Cofibrant replacements for the projective model structures can always be chosen to be flat.
Note that the functors
\begin{gather*}
    R_\tr\colon \Spc_\pt(S)\to\Spc_\tr(S,R)\text{ and}\\
    \mathord\tens_R\colon\Spc_\tr(S,R)\times\Spc_\tr(S,R)\to\Spc_\tr(S,R)
\end{gather*}
preserve flat objects. Moreover, by \cite[Theorem 4.8]{Voevodsky:2010b}, $R_\tr$ preserves objectwise equivalences between flat objects, and so does $T\tens_R\ph$ for any $T\in\Spc_\tr(S,R)$.

For every $F\in\Spc_\tr(S,R)$, we define a functorial resolution $\epsilon\colon L_\ast F\to F$ where $L_\ast F$ is flat. Let $\cat C$ be a set of representatives of isomorphism classes of objects in $\Cor(S,R)$.
The inclusion $\cat C\into\Cor(S,R)$ induces an adjunction between families of sets indexed by $\cat C$ and additive presheaves on $\Cor(S,R)$. The associated comonad $L$ has the form
\[LF=\bigoplus_{X\in\cat C}R_\tr X_+\tens F(X).\]
Here and in what follows, the unadorned tensor product is the right action of the category of \emph{sets}; that is, $R_\tr X_+\tens F(X)=\bigoplus_{F(X)}R_\tr X_+$. The right adjoint evaluates an additive presheaf $F$ on each $U\in\cat C$. Thus, the augmented simplicial object $\epsilon \colon L_\ast F\to F$ induced by this comonad is a simplicial homotopy equivalence when evaluated on any $U\in\cat C$, and in particular is an objectwise equivalence. If $F\in\Spc_\tr(S,R)$, $L_\ast F$ is defined by applying the previous construction levelwise and taking the diagonal.

\begin{lemma}\label{lem:balanced}
	If $T\in\Spc_\tr(S,R)$ is flat, then $T\tens_R\ph$ preserves objectwise equivalences.
\end{lemma}

\begin{proof}
	Since it preserves objectwise equivalences between flat objects, it suffices to show that, for any $F\in\Spc_\tr(S,R)$, $T\tens_R\epsilon\colon T\tens_R L_\ast F\to T\tens_R F$ is an objectwise equivalence. By the definition of $L_\ast$, we can obviously assume that $T$ and $F$ are functors $\Cor(S,R)^\op\to\Ab$, and since filtered colimits preserve equivalences, we can further assume that $T$ is $R_\tr X_+$ for some $X\in\Sm/S$. Given $U\in\cat C$, we will then define a candidate for an extra degeneracy operator
        \[
            s_F\colon (R_\tr X_+\tens_R F)(U)\to (R_\tr X_+\tens_R LF)(U).
        \]
		Given a finite correspondence $\psi\colon U\to Y$, define the correspondence $\psi_U\colon U\to Y\times U$ by
		\begin{equation}\label{eqn:extra}
			\psi_U=(\psi\times\id_U)\circ\Delta_U.
		\end{equation}
    By definition of $\tens_R$, we have
    \[(R_\tr X_+\tens_R F)(U)=\int^{C,D\in\Cor(S,R)}(\cor_R(C,X)\times F(D))\tens\cor_R(U,C\times D),\]
    and since $R_\tr$ is monoidal,
    \[(R_\tr X_+\tens_R LF)(U)=\bigoplus_{Y\in\cat C}\cor_R(U,X\times Y)\tens F(Y).\]
    Given $(\phi, x)\tens\psi\in(\cor_R(C,X)\times F(D))\tens\cor_R(U,C\times D)$, let
    \[(s_F)_{C,D}((\phi, x)\tens\psi)=(\phi\circ\psi_1)_U\tens\psi_2^\ast(x)\]
    in the summand indexed by $U$. One checks easily that this is an extranatural transformation and hence induces the map $s_F$. A straightforward computation shows that $L^{n+1}(\epsilon)s_{L^{n+1}F}=s_{L^nF}L^n(\epsilon)$ for all $n\geq 0$, which takes care of all the identities for a contraction of the augmented simplicial object except the identity $\epsilon s_F=\id$ which is slightly more involved.
    
    Start with $(\phi, x)\tens\psi\in(\cor_R(C,X)\times F(D))\tens\cor_R(U,C\times D)$, representing the element $[(\phi, x)\tens\psi]\in (R_\tr X_+\tens_R F)(U)$. The identity $\psi=(p_1\times\psi_2)\circ\psi_U$, which follows at once from~\eqref{eqn:extra}, shows that the element $(\phi, x)\tens\psi$ is the pushforward of $(\phi, x)\tens\psi_U$ under the pair of correpondences $p_1\colon C\times D\to C$ and $\psi_2\colon U\to D$. In the coend, it is therefore identified with the pullback of that element, which is
    \[(\phi\circ p_1, \psi_2^\ast(x))\tens\psi_U.\]
    Let
\begin{tikzmath}
	\diagram[text height=1.6ex,text depth=.5ex]{
	\bigoplus_{Y\in\cat C}(\cor_R(C,X)\times \cor_R(D,Y))\tens\cor_R(U,C\times D)\tens F(Y) \\
	(\cor_R(C,X)\times F(D))\tens\cor_R(U,C\times D) \\
	};
	\arrows
	(11) edge node[left]{\normalsize $\epsilon_{C,D}$} (21);
\end{tikzmath}
    be the family of maps inducing $\epsilon$ in the coends. By~\eqref{eqn:extra}, we have $(\phi\circ p_1\times \id_U)\circ\psi_U=(\phi\circ\psi_1)_U$. This shows that the element $(\phi\circ\psi_1)_U\tens\psi_2^\ast(x)$ is represented by
    \[(\phi\circ p_1,\id_U)\tens\psi_U\tens\psi_2^\ast(x)\in(\cor_R(C\times D,X)\times \cor_R(U,U))\tens\cor_R(U,C\times D\times U)\tens F(U),\]
    and we have
    \[\epsilon_{C\times D, U}((\phi\circ p_1,\id_U)\tens\psi_U\tens\psi_2^\ast(x))=(\phi\circ p_1, \psi_2^\ast(x))\tens\psi_U\equiv(\phi, x)\tens\psi,\]
    \ie, $\epsilon s_F([(\phi, x)\tens \psi])=[(\phi, x)\tens\psi]$.
\end{proof}

\begin{proposition}\label{prop:technical}
	Let $\cat E$ be the class of pointed spaces $X\in\Spc_\pt(S)$ with the property that, for all $F\in\Spc_\tr(S,R)$, the canonical map $\L R_\tr X\tens_R^\L F\to R_\tr X\tens_R F$ is an equivalence. Then $\cat E$ contains all flat objects and is closed under finite coproducts, filtered colimits, smash products, and quotients.
\end{proposition}

\begin{proof}
	 Let $X\in\Spc_\pt(S)$ be flat and let $\tilde X\to X$ be an objectwise equivalence where $\tilde X$ is projectively cofibrant and flat. Then $\L R_\tr X\tens_R^\L F\simeq R_\tr\tilde X\tens_RL_\ast F$ and we must show that the composition
    \[R_\tr \tilde X\tens_R L_\ast F\to R_\tr X\tens_R L_\ast F\to R_\tr X\tens_R F\]
    is an equivalence. The first map is an objectwise equivalence because $\tilde X$, $X$, and $L_\ast F$ are all flat. The second is also an objectwise equivalence by Lemma~\ref{lem:balanced}. It is clear that $\cat E$ is closed under finite coproducts, filtered colimits, and smash products. Finally, $\cat E$ is closed under quotients by Lemma~\ref{lem:ow9e0}.
\end{proof}

Any sphere $S^{p,q}$ clearly belongs to the class $\cat E$ considered in Proposition~\ref{prop:technical}, which immediately shows that~\eqref{eqn:EMLderived} is an equivalence, as promised. Furthermore, it shows that the functors $\Sigma^\infty_\tr\colon\Spc_\tr(S,R)\to\Spt_\tr(S,R)$ and $c_\tr\colon\Sp(\s\Mod_R)\to\Spt_\tr(S,R)$ preserve levelwise equivalences. An adjunction argument then shows that $c_\tr$ also preserves stable equivalences, so that~\eqref{eqn:stableEMLderived} is an equivalence.

\begin{samepage}
\begin{proposition}\label{prop:Kcolim}
	\leavevmode
	\begin{enumerate}
		\item\label{prop:Kcolim:1} For any $p\geq q\geq 0$, the functor $K(\ph(q),p)\colon\s\Mod_R\to\Spc_\pt(S)$ preserves sifted homotopy colimits and transforms finite homotopy coproducts into finite homotopy products.
		\item\label{prop:Kcolim:2} The functor $H\colon\Sp(\s\Mod_R)\to\Spt(S)$ preserves all homotopy colimits.
	\end{enumerate}
\end{proposition}
\end{samepage}

\begin{proof}
	Note that the functor $u_\tr\colon\Spc_\tr(S,R)\to\Spc_\pt(S)$ preserves sifted homotopy colimits: since homotopy colimits can be computed objectwise in both model categories this follows from the fact that the forgetful functor from simplicial abelian groups to simplicial sets preserves sifted homotopy colimits. It is clear that $u_\tr$ transforms finite sums into finite products, which proves \itemref{prop:Kcolim:1}. The stable functor $u_\tr\colon\Spt_\tr(S,R)\to\Spt(S)$ also preserves sifted homotopy colimits since they can be computed levelwise, but in addition it preserves finite homotopy colimits since it is a right Quillen functor between stable model categories. It follows that it preserves all homotopy colimits, whence \itemref{prop:Kcolim:2}.
\end{proof}

Note that any $R$-module $A$ is a sifted homotopy colimit in $\s\Mod_R$ of finitely generated free $R$-modules: if $A$ is flat then it is a filtered colimit of finitely generated free $R$-modules, and in general $A$ admits a projective simplicial resolution of which it is the homotopy colimit; furthermore, any simplicial $R$-module is the sifted homotopy colimit of itself as a diagram of $R$-modules. Thus, part \itemref{prop:Kcolim:1} of Proposition~\ref{prop:Kcolim} gives a recipe to build $K(A(q),p)$ from copies of $K(R(q),p)$ using finite homotopy products and sifted homotopy colimits.
From part \itemref{prop:Kcolim:2} of the proposition we obtain the Bockstein cofiber sequences
\begin{gather}
\label{eqn:bockstein1}
H\Z\stackrel{n}{\to} H\Z\to H(\Z/n)\text{ and}\\
\label{eqn:bockstein2}
H(\Z/n)\stackrel{n}{\to} H(\Z/n^2)\to H(\Z/n)
\end{gather}
for any $n\in\Z$. In particular, the canonical map $(H\Z)/n\to H(\Z/n)$ is an equivalence, which removes any ambiguity from the notation $H\Z/n$.

Let $f\colon T\to S$ be a morphism of base schemes, and denote by $f^\ast$ the induced derived base change functors. Note that $f^\ast$ commutes with the functor $c_\tr$. For any $p\geq q\geq 0$ and any $A\in\s\Mod_R$ there is a canonical map
\begin{equation}\label{eqn:EMLunstable}
	f^\ast K(A(q),p)_S\to K(A(q),p)_T,
\end{equation}
adjoint to the composition
\begin{multline*}
\L R_\tr f^\ast u_\tr (\L R_\tr S^{p,q}_S\tens_R^\L c_\tr A)\simeq f^\ast \L R_\tr u_\tr (\L R_\tr S^{p,q}_S\tens_R^\L c_\tr A)\\
\to f^\ast (\L R_\tr S^{p,q}_S\tens_R^\L c_\tr A) \simeq \L R_\tr f^\ast S^{p,q}_S\tens_R^\L f^\ast c_\tr A\simeq \L R_\tr S^{p,q}_T\tens_R^\L c_\tr A.
\end{multline*}
Similarly, for any $A\in\Sp(\s\Mod_R)$ there is a canonical map
\begin{equation}\label{eqn:EMLstable}
	f^\ast HA_S\to HA_T.
\end{equation}

\begin{theorem}\label{thm:EMLpb}
	Let $S$ and $T$ be essentially smooth schemes over a base scheme $U$, and let $f\colon T\to S$ be a $U$-morphism. Then \eqref{eqn:EMLunstable} and~\eqref{eqn:EMLstable} are equivalences.
\end{theorem}

\begin{proof}
	We may clearly assume that $S=U$ so that $f$ is essentially smooth. Let us consider the unstable case first.
	It suffices to show that the canonical map 
	\begin{equation}\label{eqn:futr}
		f^\ast u_\tr\to u_\tr f^\ast
	\end{equation}
	is an equivalence.
	If $f$ is smooth, then the functors $f^\ast$ have left adjoints $f_\sharp$ such that $f_\sharp\L R_\tr\simeq\L R_\tr f_\sharp$, and so \eqref{eqn:futr} is an equivalence by adjunction. In the general case, let $f$ be the cofiltered limit of the smooth maps $f_\alpha\colon T_\alpha\to S$. Let $\scr F\in\H_\tr(S,R)$.
	To show that $f^\ast u_\tr\scr F\to u_\tr f^\ast\scr F$ is an equivalence in $\H_\pt(T)$, it suffices to show that for any $X\in\Sm/T$, the induced map
	\[\Map(X_+,f^\ast u_\tr\scr F)\to\Map(X_+,u_\tr f^\ast\scr F)\simeq\Map(\L R_\tr X_+,f^\ast\scr F)\]
	is an equivalence. Write $X$ as a cofiltered limit of smooth $T_\alpha$-schemes $X_\alpha$. Then by Lemma~\ref{lem:continuity}, the above map is the homotopy colimit of the maps
	\[\Map((X_\alpha)_+,f_\alpha^\ast u_\tr\scr F)\to\Map(\L R_\tr(X_\alpha)_+,f_\alpha^\ast\scr F)\]
	which are equivalences since $f_\alpha$ is smooth. The proof that~\eqref{eqn:EMLstable} is an equivalence is similar, using Lemma~\ref{lem:stablecontinuity} instead of Lemma~\ref{lem:continuity}.
\end{proof}

\begin{remark}\label{rmk:slices}
	Using Lemma~\ref{lem:stablecontinuity} \itemref{lem:stablecontinuity:1}, one immediately verifies the hypothesis of \cite[Theorem 2.12]{Pelaez:2013} for an essentially smooth morphism $f\colon T\to S$. Thus, for every $q\in\Z$, we have a canonical equivalence $f^\ast s_q\simeq s_q f^\ast$. By Theorem~\ref{thm:EMLpb}, the equivalence $s_0\1\simeq H\Z$ proved over perfect fields in \cite[Theorems 9.0.3 and 10.5.1]{Levine:2008} holds in fact over any base scheme $S$ which is essentially smooth over a field. Since both $s_q$ and $H$ preserve homotopy colimits, we deduce that, for any $A\in\Sp(\s\Ab)$,
\[
	s_q HA\simeq\begin{cases} HA & \text{if $q=0$,} \\ 0 & \text{otherwise.} \end{cases}
\]
\end{remark}

\subsection{Representability of motivic cohomology}
\label{sub:motivic}

We prove that the Eilenberg–Mac Lane spaces and spectra represent motivic cohomology of essentially smooth schemes over fields. If $k$ is a field and $X$ is a smooth $k$-scheme, the motivic cohomology groups $H^{p,q}(X,A)$ are defined in \cite[Definition 3.4]{Mazza:2006} for any abelian group $A$. By \cite[Proposition 3.8]{Mazza:2006}, these groups do not depend on the choice of the field $k$ on which $X$ is smooth. More generally, if $X$ is an essentially smooth scheme over a field $k$, cofiltered limit of smooth $k$-schemes $X_\alpha$, we define
\[H^{p,q}(X,A)=\colim_\alpha H^{p,q}(X_\alpha,A).\]
This does not depend on the choice of the diagram $(X_\alpha)$ since it is unique as a pro-object in the category of smooth $k$-schemes (\cite[Corollaire 8.13.2]{EGA4-3}). Moreover, by \cite[Lemma 3.9]{Mazza:2006}, this definition is still independent of the choice of $k$.

Let now $k$ be a perfect field and let $\DM^{\eff,-}(k,R)$ be Voevodsky's triangulated category of effective motives: it is the homotopy category of bounded below $\A^1$-local chain complexes of Nisnevich sheaves of $R$-modules with transfers. Normalization defines a symmetric monoidal functor
\[N\colon \H_\tr(k,R)\to\DM^{\eff,-}(k,R)\]
(see \cite[\S1.2]{MEMS}) which is fully faithful by \cite[Theorem 1.15]{MEMS}. By \cite[Proposition 14.16]{Mazza:2006} (where $k$ is implicitly assumed to be perfect), for any $X\in\Sm/k$ and any $R$-module $A$ there is a canonical isomorphism
\begin{equation}\label{eqn:mot1}
	H^{p,q}(X,A)\cong [N\L R_\tr X_+,A(q)[p]]
\end{equation}
for all $p,q\in\Z$, where the chain complex $A(q)\in\DM^{\eff,-}(k,R)$ is defined as follows:
\[
A(q)=
\begin{cases}
	N(\L R_\tr S^{q,q}\tens_R^\L c_\tr A)[-q] & \text{if $q\geq 0$,} \\
	0 & \text{if $q<0$.}
\end{cases}
\]
Since $N$ is fully faithful, we have, for every $F\in\H_\pt(k)$ and every $p\geq q\geq 0$,
\begin{equation}\label{eqn:mot2}
	[N\L R_\tr F,A(q)[p]]\cong[F,K(A(q),p)].
\end{equation}
Moreover, by \cite[Theorem 2.4]{Voevodsky:2003}, we have bistability isomorphisms
\begin{equation}\label{eqn:mot3}
	[N\L R_\tr F,A(q)[p]]\cong [N\L R_\tr\Sigma^{r,s} F, A(q+s)[p+r]]
\end{equation}
for every $p,q\in\Z$, $r\geq s\geq 0$, and $F\in\H_\pt(k)$.

\begin{theorem}\label{thm:motivic}
    Assume that $S$ is essentially smooth over a field. Let $A$ be an $R$-module and $X\in\Sm/S$. For any $p\geq q\geq 0$ and $r\geq s\geq 0$, there is a natural isomorphism
    \[H^{p-r,q-s}(X,A)\cong[\Sigma^{r,s}X_+,K(A(q),p)].\]
    For any $p,q\in\Z$, there is a natural isomorphism
    \[H^{p,q}(X,A)\cong[\Sigma^\infty X_+,\Sigma^{p,q}HA].\]
\end{theorem}

\begin{proof}
	Suppose first that $S$ is the spectrum of a perfect field $k$. Then the first isomorphism is a combination of the isomorphisms \eqref{eqn:mot1}, \eqref{eqn:mot2}, and \eqref{eqn:mot3}. From the latter two it also follows that the canonical maps
	 \begin{equation}\label{eqn:Omegaspectrum}
		 K(A(q),p)_k\to\Omega^{r,s}K(A(q+s),p+r)_k
	\end{equation}
	 are equivalences, so the second isomorphism follows from the first one and the definition of $HA$.
	 
	 In general, choose an essentially smooth morphism $f\colon S\to\Spec k$ where $k$ is a perfect field, and let $f^\ast\colon\H_\pt(k)\to\H_\pt(S)$ be the corresponding base change functor. Let $f$ be the cofiltered limit of the smooth morphisms $f_\alpha\colon S_\alpha\to k$, and let $X$ be the cofiltered limit of the smooth $S_\alpha$-schemes $X_\alpha$. By Theorem~\ref{thm:EMLpb}, $f^\ast K(A(q),p)_{k}\simeq K(A(q),p)_S$.
 By Lemma~\ref{lem:continuity} \itemref{lem:continuity:1}, we therefore have
\[
	[\Sigma^{r,s}X_+,K(A(q),p)_S]\cong\colim_\alpha[\Sigma^{r,s}(X_\alpha)_+,f_\alpha^\ast K(A(q),p)_{k}].
\]
Using the left adjoint $(f_\alpha)_\sharp$ of $f_\alpha^\ast$, we obtain the isomorphisms
\[
	[\Sigma^{r,s}(X_\alpha)_+,f_\alpha^\ast K(A(q),p)_{k}]\cong[\Sigma^{r,s}(X_\alpha)_+, K(A(q),p)_{k}]\cong H^{p-r,q-s}(X_\alpha,A).
\]
Finally, the colimit of the right-hand side is $H^{p-r,q-s}(X,A)$ by definition. The second isomorphism can either be proved in the same way, using Lemma~\ref{lem:stablecontinuity}, or it can be deduced from~\eqref{eqn:Omegaspectrum} and the fact (observed in Appendix~\ref{sec:appendix}) that $f^\ast\Omega^{r,s}\simeq\Omega^{r,s}f^\ast$.
\end{proof}

The following corollary summarizes the standard vanishing results for motivic cohomology (which we use freely later on).

\begin{samepage}
\begin{corollary}\label{cor:vanish}
	Assume that $S$ is essentially smooth over a field. Let $X\in\Sm/S$ and $p,q\in\Z$ satisfy any of the following conditions:
	\begin{enumerate}
		\item $q<0$;
		\item $p>q+\essdim X$;
		\item $p>2q$;
	\end{enumerate}
	where $\essdim X$ is the least integer $d$ such that $X$ can be written as a cofiltered limit of smooth $d$-dimensional schemes over a field.
Then, for any $R$-module $A$, $[\Sigma^\infty X_+, \Sigma^{p,q}HA]=0$.
\end{corollary}
\end{samepage}

\begin{proof}
	By Theorem~\ref{thm:motivic}, we must show that $H^{p,q}(X,A)=0$. If $X$ is smooth over a field, we have $H^{p,q}(X,A)=0$ when $q<0$ or $p>q+\dim X$ by definition of motivic cohomology and \cite[Theorem 3.6]{Mazza:2006}, respectively. If $X$ is smooth over a perfect field and $p>2q$, then $H^{p,q}(X,A)=0$ by \cite[Theorem 19.3]{Mazza:2006}. For general $X$ the result follows by the definition of motivic cohomology of essentially smooth schemes over fields.
\end{proof}

If $S$ is smooth over a field, recall from \cite[Corollary 4.2]{Mazza:2006} that
\begin{equation}\label{eqn:picard}
	H^{2,1}(S,\Z)\cong\Pic(S).
\end{equation}
This generalizes immediately to schemes $S$ that are essentially smooth over a field. It follows from Theorem~\ref{thm:motivic} that the motivic ring spectrum $HR\in\SH(S)$ can be oriented as follows. Given $X\in\Sm/S$ and a line bundle $\scr L$ over $X$, define $c_1(\scr L)\in H^{2,1}(X,R)$ to be the image of the integral cohomology class corresponding to the isomorphism class of $\scr L$ in $\Pic(X)$. By the universality of $\MGL$ (\cite[Theorem 3.1]{NSO2}), this orientation is equivalently determined by a morphism of ring spectra
\[\th\colon\MGL\to HR.\]
Moreover, if $S$ and $T$ are both essentially smooth over a field and if $f\colon T\to S$ is any morphism, then the canonical equivalence $f^\ast HR_S\simeq HR_T$ of Theorem~\ref{thm:EMLpb} is compatible with the orientations; this follows at once from the naturality of the isomorphism~\eqref{eqn:picard}. In other words, through the identifications $f^\ast\MGL_S\simeq\MGL_T$ and $f^\ast HR_S\simeq HR_T$, we have $f^\ast(\th_S)=\th_T$.

\section{Operations and co-operations in motivic cohomology}
\label{sec:operations_in_motivic_cohomology}

In this section, the base scheme $S$ is essentially smooth over a field and $\ell\neq\Char S$ is a prime number. We recall the structures of the motivic Steenrod algebra over $S$ and its dual, and we compute the $H\Z/\ell$-module $H\Z/\ell\wedge H\Z$.

\subsection{Duality and Künneth formulas}\label{sub:duality}

In this paragraph we formulate a convenient finiteness condition on the homology of cellular spectra that ensures that their homology and cohomology are dual to one another and satisfy Künneth formulas. We fix a commutative ring $R$. Given an $HR$-module $M$, we denote by $M^\vee$ its dual in the symmetric monoidal category $\D(HR)$. Note that if $M=HR\wedge E$, then $\pi_\aast M$ is the motivic homology of $E$ and $\pi_{-\ast,-\ast} M^\vee$ is its motivic cohomology (with coefficients in $R$).

\begin{remark}
	We always consider bigraded abelian groups as a symmetric monoidal category where the symmetry is Koszul-signed with respect to the first grading, \cf{} \cite[\S3]{Naumann:2009}. For any oriented ring spectrum $E$, $E_\aast:=\pi_\aast E$ is a commutative monoid in this symmetric monoidal category.
\end{remark}

An $HR$-module will be called \emph{split} if it is equivalent to an $HR$-module of the form
\[\bigvee_\alpha\Sigma^{p_\alpha,q_\alpha}HR.\]
Note that if $M$ is split and $S$ is not empty, the family of bidegrees $(p_\alpha,q_\alpha)$ is uniquely determined: this follows for example from Remark~\ref{rmk:slices}: if $M\simeq\bigvee_{p,q\in\Z}\Sigma^{p,q}HV_{p,q}$ where $V_{p,q}$ is an $R$-module, then $V_{p,q}\cong\pi_{p,q}s_q j^\ast M$ where $j$ is the inclusion of a connected component of $S$. Split $HR$-modules are obviously cellular, but, unlike in topology, the converse does not hold even if $R$ is a field: the motivic spectrum constructed in \cite[Remark~2.1]{Voevodsky:2002} representing étale cohomology with coefficients in $\mu_\ell^{\tens\ast}$ is a cellular $H\Z/\ell$-module, yet it is not split since all its slices are zero.

\begin{lemma}\label{lem:kunneth}
	Let $M$ and $N$ be $HR$-modules. The canonical map
	\[\pi_\aast M\tens_{HR_\aast}\pi_\aast N\to\pi_\aast(M\wedge_{HR}N)\]
	is an isomorphism under either of the following conditions:
	\begin{enumerate}
		\item\label{lem:kunneth:1} $M$ is cellular and $\pi_\aast N$ is flat over $HR_\aast$;
		\item\label{lem:kunneth:2} $M$ is split.
	\end{enumerate}
\end{lemma}

\begin{proof}
	Assuming \itemref{lem:kunneth:1}, this is a natural map between homological functors of $M$, so we may assume \itemref{lem:kunneth:2}, in which case the result is obvious.
\end{proof}

\begin{samepage}
\begin{lemma}\label{lem:split}
	For an $HR$-module $M$, the following conditions are equivalent:
	\begin{enumerate}
		\item\label{lem:split:1} $M$ is split;
		\item\label{lem:split:2} $\spi_\aast M$ is free over $\spi_\aast HR$;
		\item\label{lem:split:3} $M$ is cellular and $\pi_\aast M$ is free over $\pi_\aast HR$.
	\end{enumerate}
\end{lemma}
\end{samepage}

\begin{proof}
	It is clear that \itemref{lem:split:1} implies \itemref{lem:split:2} and \itemref{lem:split:3}. Assuming \itemref{lem:split:2} or \itemref{lem:split:3}, use a basis of $\spi_\aast M$ (or $\pi_\aast M$) to define a morphism of $HR$-modules $\bigvee_\alpha\Sigma^{p_\alpha,q_\alpha}HR\to M$. This morphism is then a $\spi_\aast$-isomorphism (or a $\pi_\aast$-isomorphism between cellular $HR$-modules) and so it is an equivalence.
\end{proof}

\begin{definition}
	A split $HR$-module is called \emph{psf} (short for \emph{proper and slicewise finite}) if it is equivalent to $\bigvee_\alpha\Sigma^{p_\alpha,q_\alpha}HR$ where the bidegrees $(p_\alpha,q_\alpha)$ satisfy the following conditions: they are all contained in the cone $q\geq 0$, $p\geq 2q$, and for each $q$ there are only finitely many $\alpha$ such that $q_\alpha=q$.
\end{definition}

This condition is satisfied in many interesting cases. For example, if $E$ is the stabilization of any Grassmannian or Thom space of the tautological bundle thereof, or if $E=\MGL$, then $E$ is cellular and the calculus of oriented cohomology theories (see \S\ref{sub:hurewicz}) shows that $HR_\aast E$ is free over $HR_\aast$, with finitely many generators in each bidegree $(2n,n)$; it follows from Lemma~\ref{lem:split} that $HR\wedge E$ is psf. Later we will show that, for $\ell\neq \Char S$ a prime number, $H\Z/\ell\wedge H\Z/\ell$ and $H\Z/\ell\wedge H\Z$ are psf.

\pagebreak[2]

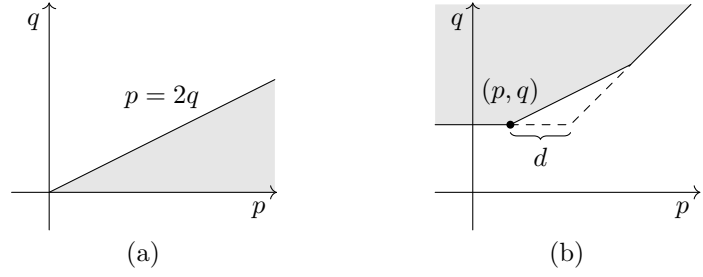
\begin{figure}[htb]
\centering
\begin{tikzpicture}
	\fill[color=gray!20] (.5cm,.5cm) -- (3.5cm, .5cm) -- (3.5cm,2cm) -- cycle;
	\draw (.5cm,.5cm) -- node[above=7]{$p=2q$} (3.5cm, 2cm);
	\draw[->] (.5cm,0) -- (.5cm, 3cm) node[below left]{$q$};
	\draw[->] (0,.5cm) -- node[below=15]{(a)} (3.5cm, .5cm) node[below left]{$p$};
\end{tikzpicture}
\hskip 2cm
\begin{tikzpicture}[decoration=brace]
	\fill[color=gray!20] (0,1.4cm) -- (1cm, 1.4cm) -- (2.6cm,2.2cm) -- (3.4cm,3cm) -- (0, 3cm) -- cycle;
	\draw (0,1.4cm) -- (1cm, 1.4cm) -- (2.6cm,2.2cm) -- (3.4cm,3cm);
	\draw[dashed] (1cm, 1.4cm) -- ( 1.8cm,1.4cm) -- (2.6cm,2.2cm);
	\draw[decorate,decoration={brace,raise=3pt}] ( 1.8cm,1.4cm) -- node[below=5pt]{$d$} (1cm, 1.4cm);
	\fill (1cm, 1.4cm) circle (1.5pt) node[above=4pt]{$(p,q)$};
	\draw[->] (.5cm,0) -- (.5cm, 3cm) node[below left]{$q$};
	\draw[->] (0,.5cm) -- node[below=15]{(b)} (3.5cm, .5cm) node[below left]{$p$};
\end{tikzpicture}
\caption{(a) The proper cone. (b) If $d=\essdim X$, the shaded area is the potentially nonzero locus of $H^{\ast-p,\ast-q}(X,R)$ according to Corollary~\ref{cor:vanish}. If for every $(p,q)$ and every $d$ there are only finitely many bidegrees $(p_\alpha,q_\alpha)$ in the shaded area, it follows that the canonical map $\bigvee_\alpha\Sigma^{p_\alpha,q_\alpha}HR\to\prod_\alpha\Sigma^{p_\alpha,q_\alpha}HR$ is a $\spi_\aast$-isomorphism.}
\label{fig:psf}
\end{figure}

\begin{samepage}
\begin{proposition}\label{prop:psf}
	Let $M$ and $N$ be psf $HR$-modules. Then
	\begin{enumerate}
		\item\label{prop:psf:1} $M\wedge_{HR}N$ is psf;
		\item\label{prop:psf:2} $M^\vee$ is split;
		\item\label{prop:psf:3} the pairing $M^\vee\wedge_{HR}M\to HR$ is perfect;
		\item\label{prop:psf:4} the canonical map $M^\vee\wedge_{HR}N^\vee\to(M\wedge_{HR}N)^\vee$ is an equivalence;
		\item\label{prop:psf:5} the pairing $\pi_\aast M^\vee\tens_{HR_\aast}\pi_\aast M\to HR_\aast$ is perfect;
		\item\label{prop:psf:6} the canonical map
		$\pi_\aast M\tens_{HR_\aast}\pi_\aast N\to \pi_\aast(M\wedge_{HR}N)$
		is an isomorphism;
		\item\label{prop:psf:7} the canonical map
		$\pi_\aast M^\vee\tens_{HR_\aast}\pi_\aast N^\vee\to \pi_\aast(M\wedge_{HR}N)^\vee$
		is an isomorphism.
	\end{enumerate}
\end{proposition}
\end{samepage}

\begin{proof}
	Assertion \itemref{prop:psf:1} is clear from the definition. Let $M=\bigvee_\alpha\Sigma^{p_\alpha,q_\alpha} HR$. Corollary~\ref{cor:vanish} and the psf condition imply that the canonical maps \[\bigvee_\alpha\Sigma^{p_\alpha,q_\alpha}HR\to\prod_\alpha\Sigma^{p_\alpha,q_\alpha}HR\quad\text{and}\quad\bigvee_\alpha\Sigma^{-p_\alpha,-q_\alpha}HR\to\prod_\alpha\Sigma^{-p_\alpha,-q_\alpha}HR\] are equivalences (compare Figure~\ref{fig:psf} (a) and (b)). This implies~\itemref{prop:psf:2}, \itemref{prop:psf:3}, and~\itemref{prop:psf:4}. In particular, the two inclusions $\bigoplus_\alpha\pi_\aast\Sigma^{\pm p_\alpha,\pm q_\alpha}HR\into\prod_\alpha\pi_\aast\Sigma^{\pm p_\alpha,\pm q_\alpha}HR$ are isomorphisms, which shows~\itemref{prop:psf:5}.
	Assertion \itemref{prop:psf:6} is a just special case of Lemma~\ref{lem:kunneth}, and assertion \itemref{prop:psf:7} follows from~\itemref{prop:psf:2}, \itemref{prop:psf:4}, and Lemma~\ref{lem:kunneth}.
\end{proof}

\subsection{The motivic Steenrod algebra}

For the rest of this section we fix a prime number $\ell\neq\Char S$, and we abbreviate $K(\Z/\ell(n),2n)$ to $K_n$ and $H\Z/\ell$ to $H$. We denote by $\scr A^\aast$ the motivic Steenrod algebra at $\ell$. By this we mean the algebra of all bistable natural transformations $\tilde H^\aast(\ph,\Z/\ell)\to\tilde H^\aast(\ph,\Z/\ell)$ (as functors on the pointed homotopy category $\H_\pt(S)$), that is,
\[\scr A^\aast=\lim_{n\to\infty}\tilde H^{\ast+2n,\ast+n}(K_n).\]

In \cite[\S9]{Voevodsky:2003}, the reduced power operations \[P^{i}\in \scr A^{2i(\ell-1),i(\ell-1)}\] are constructed for all $i\geq 0$, provided that $S$ be the spectrum of a perfect field. By inspection of their definitions, if $f\colon\Spec k'\to \Spec k$ is an extension of perfect fields, then, under the identifications $f^\ast K_n\simeq K_n$ and $f^\ast H\simeq H$, we have $f^\ast(P^i)=P^i$. If $S$ is essentially smooth over $k$, the reduced power operations over $k$ therefore induce reduced power operations over $S$ which are independent of the choice of $k$.

Given a sequence of integers $(\epsilon_0,i_1,\epsilon_1,\dotsc,i_r,\epsilon_r)$ satisfying $i_j>0$, $\epsilon_j\in\{0,1\}$, and $i_{j}\geq \ell i_{j+1}+\epsilon_j$, we can form the operation $\beta^{\epsilon_0}P^{i_1}\dotso P^{i_r}\beta^{\epsilon_r}$, where $\beta\colon H\to\Sigma^{1,0}H$ is the Bockstein morphism defined by the cofiber sequence~\eqref{eqn:bockstein2}; the analogous operations in topology form a $\Z/\ell$-basis of the topological mod $\ell$ Steenrod algebra $\scr A^\ast$, and so we obtain a map of left $H^{\ast\ast}$-modules
\begin{equation}\label{eqn:steenrod}
    H^{\ast\ast}\tens_{\Z/\ell}\scr A^\ast\to \scr A^\aast.
\end{equation}

\begin{lemma}\label{lem:steenrod}
    The map~\eqref{eqn:steenrod} is an isomorphism. In particular, the algebra $\scr A^\aast$ is generated by the reduced power operations $P^i$,  the Bockstein $\beta$, and the operations $u\mapsto au$ for $a\in H^\aast(S,\Z/\ell)$.
\end{lemma}

\begin{proof}
	If $S$ is the spectrum of a field of characteristic zero, this is \cite[Theorem 3.49]{MEMS}; the general case is proved in \cite{HKO:2013}.
\end{proof}

By a \emph{split proper Tate object of weight $\geq n$} we mean an object of $\H_\tr(S,\Z/\ell)$ which is a direct sum of objects of the form $\L\Z/\ell_\tr S^{p,q}$ with $p\geq 2q$ and $q\geq n$.

\begin{lemma}\label{lem:splittate}
	$\L \Z/\ell_\tr K_n$ is split proper Tate of weight $\geq n$.
\end{lemma}

\begin{proof}
If $S$ is the spectrum of a field admitting resolutions of singularities, this is proved in \cite[Corollary 3.33]{MEMS}; the general case is proved in \cite{HKO:2013}.
\end{proof}

\begin{lemma}\label{lem:steenrod2}
    The canonical map $H^\aast H\to \scr A^\aast$ is an isomorphism.
\end{lemma}

\begin{proof}
    This map fits in the exact sequence
    \[0\to\Rlim_{n\to\infty} \tilde H^{p-1+2n,q+n}(K_n)\to H^{p,q}H\to\lim_{n\to\infty} \tilde H^{p+2n,q+n}(K_n)\to 0,\]
    and we must show that the $\Rlim$ term vanishes. 
    By Lemma~\ref{lem:splittate}, $\L\Z/\ell_\tr K_n\simeq\Sigma^{2n,n}M_n$ where $M_n$ is split proper Tate of weight $\geq 0$. All functors should be derived in the following computations. Using the standard adjunctions, we get
    \begin{multline*}
        \tilde H^{p-1+2n,q+n}(K_n)\cong[\Sigma^\infty K_n,\Sigma^{p-1+2n,q+n}H\Z/\ell]\cong[\Sigma^\infty_\tr\Z/\ell_\tr K_n,\Sigma^{p-1+2n,q+n}\Z/\ell_\tr\1]\\
        \cong[\Sigma^{2n,n}\Sigma^\infty_\tr M_n,\Sigma^{p-1+2n,q+n}\Z/\ell_\tr\1]\cong[\Sigma^\infty_\tr M_n,\Sigma^{p-1,q}\Z/\ell_\tr\1].
    \end{multline*}
	 To show that $\Rlim [\Sigma^\infty_\tr M_n,\Sigma^{p-1,q}\Z/\ell_\tr\1]=0$, it remains to show that the cofiber sequence
    \[\bigoplus_{n\geq 0}\Sigma^\infty_\tr M_n\to\bigoplus_{n\geq 0}\Sigma^\infty_\tr M_n\to\hocolim_{n\to\infty}\Sigma^\infty_\tr M_n\]
	 splits in $\SH_\tr(S,\Z/\ell)$. If $S$ is the spectrum of a perfect field, this follows from \cite[Corollary 2.71]{MEMS}. In general, let $f\colon S\to\Spec k$ be essentially smooth where $k$ is a perfect field. Then $f^\ast M_n\simeq M_n$ by Theorem~\ref{thm:EMLpb}, so the above cofiber sequence splits.
\end{proof}

As a consequence of Lemma~\ref{lem:steenrod2}, $H^\aast E$ is a left module over $\scr A^\aast$ for every spectrum $E$. The proof of the following theorem was also given in \cite[\S6]{Dugger:2010} for fields of characteristic zero.

\begin{theorem}\label{thm:HsmashH}
	$H\wedge H$ is a psf $H$-module.
\end{theorem}

\begin{proof}
	By Lemma~\ref{lem:splittate}, $\L \Z/\ell_\tr K_n\simeq\Sigma^{2n,n}M_n$ where $M_n$ is split proper Tate of weight $\geq 0$. By \cite[Corollary 2.71]{MEMS} and Theorem~\ref{thm:EMLpb}, $\hocolim_{n\to\infty} M_n$ is again a split proper Tate object of weight $\geq 0$, \ie, can be written in the form
    \[\hocolim_{n\to\infty} M_n\simeq\bigoplus_\alpha \L \Z/\ell_\tr S^{p_\alpha,q_\alpha}\]
    with $p_\alpha\geq 2q_\alpha\geq 0$.
    In the following computations, all functors must be appropriately derived. We have
    \begin{multline*}\textstyle
         \Z/\ell_\tr H\simeq\Z/\ell_\tr\colim\Sigma^{-2n,-n}\Sigma^\infty K_n\simeq \colim\Sigma^{-2n,-n}\Sigma^\infty_\tr \Z/\ell_\tr K_n\\
        \simeq \colim\Sigma^{-2n,-n}\Sigma^\infty_\tr \Sigma^{2n,n} M_n\simeq \colim\Sigma^\infty_\tr M_n
        \simeq \Sigma^\infty_\tr \colim M_n\\
        \simeq \Sigma^\infty_\tr\bigoplus_\alpha \Z/\ell_\tr S^{p_\alpha,q_\alpha}
        \simeq \Z/\ell_\tr \bigvee_\alpha \Sigma^\infty S^{p_\alpha,q_\alpha};
    \end{multline*}
	 in particular, $\Z/\ell_\tr H =\Phi(H\wedge H)$ is cellular. By Lemma~\ref{lem:HZmod}, we obtain the equivalences
    \[H\wedge H\simeq H\wedge\bigvee_\alpha \Sigma^\infty S^{p_\alpha,q_\alpha}\simeq\bigvee_\alpha\Sigma^{p_\alpha,q_\alpha}H\]
    in $\D(H)$. It remains to identify the bidegrees $(p_\alpha,q_\alpha)$. By Theorem~\ref{thm:EMLpb} we may as well assume that $S$ is the spectrum of an algebraically closed field, so that $H^\aast(S,\Z/\ell)\cong\Z/\ell[\tau]$ with $\tau$ in degree $(0,1)$. By Lemma~\ref{lem:steenrod2} and the above decomposition, we have
   \[\scr A^{\ast\ast}\cong [H,\Sigma^{\ast\ast}H]\cong [H\wedge H,\Sigma^{\ast\ast}H]_H\cong\prod_\alpha[H,\Sigma^{\ast-p_\alpha,\ast-q_\alpha}H]_H\cong\prod_\alpha H^{\ast-p_\alpha,\ast-q_\alpha}.\]
   On the other hand, by Lemma~\ref{lem:steenrod}, $\scr A^{p,q}$ is finite for every $p,q\in\Z$. This implies that the above product is a direct sum, and hence that the bidegrees $(p_\alpha,q_\alpha)$ are the bidegrees of a basis of $\scr A^\aast$ over $H^\aast$; in particular, $H\wedge H$ is psf.
\end{proof}

By Theorem~\ref{thm:HsmashH} and Proposition~\ref{prop:psf} \itemref{prop:psf:7}, we have a canonical isomorphism
\[H^\aast(H\wedge H)\cong H^\aast H\tens_{H^\aast}H^\aast H\]
(where the tensor product $\tens_{H^\aast}$ uses the \emph{left} $H^\aast$-module structure on both sides). The multiplication $H\wedge H\to H$ therefore induces a coproduct
\[\Delta\colon \scr A^\aast\to \scr A^\aast\tens_{H^\aast}\scr A^\aast.\]
If $S$ is the spectrum of a perfect field, this coproduct coincides with the one studied in \cite[\S11]{Voevodsky:2003} by virtue of \cite[Lemma 11.6]{Voevodsky:2003}.

\subsection{The Milnor basis}
\label{sub:the_dual_motivic_steenrod_algebra}

Let $\scr A_\aast$ denote the dual of $\scr A^{-\ast,-\ast}$ in the symmetric monoidal category of left $H_\aast$-modules. By Theorem~\ref{thm:HsmashH}, the $H$-module $H\wedge H$ is psf, and by Lemma~\ref{lem:steenrod2}, $\scr A^{-\ast,-\ast}=\pi_\aast(H\wedge H)^\vee$. It follows from Proposition~\ref{prop:psf} \itemref{prop:psf:5} that $\scr A_\aast=\pi_\aast(H\wedge H)$ and that $\scr A^{-\ast,-\ast}$ is in turn the dual of $\scr A_\aast$. Moreover, by Lemma~\ref{lem:kunneth}, there is a canonical isomorphism
\begin{equation}\label{eqn:Akunneth}
	\scr A_\aast \tens_{H_\aast} \pi_\aast M\cong \pi_\aast(H\wedge M)
\end{equation}
for any $M\in\D(H)$. Applying~\eqref{eqn:Akunneth} with $M=H^{\wedge i}$, we deduce by induction on $i$ that
\[H_\aast(H^{\wedge i})\cong\scr A_\aast^{\tens i}\]
(the tensor product being over $H_\aast$), so that $(H_\aast,\scr A_\aast)$ is a Hopf algebroid. Applying~\eqref{eqn:Akunneth} with $M=H\wedge E$, we deduce further that $H_\aast E$ is a left comodule over $\scr A_\aast$, for any $E\in\SH(S)$.

Define a Hopf algebroid $(A, \Gamma)$ as follows. Let 
\begin{align*}
	A&=\Z/\ell[\rho,\tau],\\
	\Gamma &=A[\tau_0,\tau_1,\dotsc,\xi_1,\xi_2,\dotsc]/(\tau_i^2-\tau\xi_{i+1}-\rho\tau_{i+1}-\rho\tau_0\xi_{i+1}).
\end{align*}
The structure maps $\eta_L$, $\eta_R$, $\epsilon$, and $\Delta$ are given by the formulas
\begin{align*}
	\eta_L\colon A\to\Gamma,\quad &\eta_L(\rho)=\rho,\\
	&\eta_L(\tau)=\tau,\\
	\eta_R\colon A\to\Gamma,\quad &\eta_R(\rho)=\rho,\\
	&\eta_R(\tau)=\tau+\rho\tau_0,\\
	\epsilon\colon\Gamma\to A,\quad &\epsilon(\rho)=\rho,\\
	&\epsilon(\tau)=\tau,\\
	&\epsilon(\tau_r)=0,\\
	&\epsilon(\xi_{r})=0,\\
	\Delta\colon \Gamma\to\Gamma\tens_A\Gamma,\quad& \Delta(\rho)=\rho\tens 1,\\
	&\Delta(\tau)=\tau\tens 1,\\
	&\Delta(\tau_r)=\tau_r\tens 1+1\tens\tau_r+\sum_{i=0}^{r-1}\xi_{r-i}^{\ell^i}\tens \tau_i,\\
	&\Delta(\xi_{r})=\xi_r\tens 1+1\tens \xi_r+\sum_{i=1}^{r-1}\xi_{r-i}^{\ell^i}\tens\xi_i.
\end{align*}
The coinverse map $\iota\colon\Gamma\to\Gamma$ is determined by the identities it must satisfy. Namely, we have
\begin{align*}
 	&\iota(\rho)=\rho,\\
 	&\iota(\tau)=\tau+\rho\tau_0,\\
 	&\iota(\tau_r)=-\tau_r-\sum_{i=0}^{r-1}\xi_{r-i}^{\ell^i} \iota(\tau_i),\\
 	&\iota(\xi_r)=-\xi_r-\sum_{i=1}^{r-1}\xi_{r-i}^{\ell^i} \iota(\xi_i).
\end{align*}
We will not use this map.

We view $H_\aast$ as an $A$-algebra via the map $A\to H_\aast$ defined as follows: if $\ell$ is odd it sends both $\rho$ and $\tau$ to $0$, while if $\ell=2$ it sends $\rho$ to the image of $-1\in\G_m(S)$ in \[H_{-1,-1}=H^1_\et(S,\mu_2)\] and $\tau$ to the nonvanishing element of \[H_{0,-1}=\mu_2(S)\cong\Hom(\pi_0(S),\Z/2)\] (recall that $\Char S\neq 2$ if $\ell=2$). We will also denote by $\rho,\tau\in H_\aast$ the images of $\rho,\tau\in A$ under this map; so if $\ell\neq 2$, $\rho=\tau=0$ in $H_\aast$. All the arguments in this paper work regardless of what $\rho$ and $\tau$ are, and with this setup we will not have to worry about the parity of $\ell$ from now on.

\begin{theorem}\label{thm:A**}
	$\scr A_\aast$ is isomorphic to $\Gamma\tens_AH_\aast$ with
	\[\abs{\tau_r}=(2\ell^r-1,\ell^r-1)\quad\text{and}\quad\abs{\xi_r}=(2\ell^r-2,\ell^r-1).\]
	The map $H_\aast\to\scr A_\aast$ dual to the left action of $\scr A^\aast$ on $H^\aast$ is a left coaction of $(A,\Gamma)$ on the ring $H_\aast$, and the Hopf algebroid $(H_\aast,\scr A_\aast)$ is isomorphic to the twisted tensor product of $(A,\Gamma)$ with $H_\aast$.
\end{theorem}

This means that
\begin{itemize}
	\item $\scr A_\aast=\Gamma\tens_AH_\aast$;
	\item $\eta_L$ and $\epsilon$ are extended from $(A,\Gamma)$;
	\item $\eta_R\colon H_\aast\to\scr A_\aast$ is the coaction;
	\item $\Delta\colon\scr A_\aast\to \scr A_\aast\tens_{H_\aast}\scr A_\aast$ is induced by the comultiplication of $\Gamma$ and the map $\eta_R$ to the second factor;
	\item $\iota\colon \scr A_\aast\to\scr A_\aast$ is induced by the coinverse of $\Gamma$ and $\eta_R$.
\end{itemize}

\begin{proof}[Proof of Theorem~\ref{thm:A**}]
	If $S$ is the spectrum of a perfect field, this is proved in \cite[\S12]{Voevodsky:2003}. In general, choose an essentially smooth morphism $f\colon S\to\Spec k$ where $k$ is a perfect field. Note that the induced map $(H_k)_\aast\to (H_S)_\aast$ is a map of $A$-algebras. It remains to observe that the Hopf algebroid $\scr A_\aast$ is obtained from $(H_k)_\aast H_k$ by extending scalars from $(H_k)_\aast$ to $(H_S)_\aast$, which follows formally from the following facts: $f^\ast$ is a symmetric monoidal functor, $f^\ast(H_k)\simeq H_S$ as ring spectra (Theorem~\ref{thm:EMLpb}), and $H_k^{\wedge i}$ is a split $H_k$-module (Theorem~\ref{thm:HsmashH}).
\end{proof}

As usual, for a sequence $E=(\epsilon_0,\epsilon_1,\dotsc)$ with $\epsilon_i\in\{0,1\}$ and $\epsilon_i=0$ for almost all $i$, we set
\[\tau(E)=\tau_0^{\epsilon_0}\tau_1^{\epsilon_1}\dotso,\]
and for a sequence $R=(r_1,r_2,\dotsc)$ of nonnegative integers (almost all zero) we set
\[\xi(R)=\xi_1^{r_1}\xi_2^{r_2}\dotso.\]
Sequences can be added termwise, and we write $R'\subset R$ if there exists a sequence $R''$ such that $R'+R''=R$. We write $\emptyset$ for a sequence of zeros.

The products $\tau(E)\xi(R)$ form a basis of $\scr A_\aast$ as a left $H_\aast$-module. If $\rho(E,R)\in\scr A^\aast$ is dual to $\tau(E)\xi(R)$ with respect to that basis, then
\[\rho(E,R)=Q(E)P^R=Q_0^{\epsilon_0}Q_1^{\epsilon_1}\dotso P^R,\]
where $Q(E)$ is dual to $\tau(E)$, $Q_i$ to $\tau_i$, and $P^R$ to $\xi(R)$, and we have $P^n=P^{n,0,0,\dotsc}$ and $\beta=Q_0$. We set
\[q_i=P^{e_i}\]
where $e_0=\emptyset$ and, for $i\geq 1$, $e_i$ is the sequence with $1$ in the $i$th position and $0$ elsewhere. By dualizing the comultiplication of $\scr A_\aast$ we see at once that, for $i\geq 1$,
\[Q_i=q_i\beta-\beta q_i\quad\text{and}\quad q_i=P^{\ell^{i-1}}\dotso P^{\ell}P^1.\]

The following lemma completely describes the coproduct on the basis elements $\rho(E,R)$. It is proved by dualizing the product on $\scr A_\aast$. Explicit formulas for the products of elements $\rho(E,R)$ are more complicated, so we will not attempt to derive them.

\begin{samepage}
\begin{lemma}[Cartan formulas]\label{lem:cartan}
	\leavevmode
	\begin{itemize}
		\item $\Delta(P^R)=\sum_{E=(\epsilon_0,\epsilon_1,\dotsc)}\sum_{R_1+R_2=R-E}\tau^{\sum_{i\geq 0} \epsilon_i}Q(E)P^{R_1}\tens Q(E)P^{R_2}$;
		\item $\Delta(Q_i)=Q_i\tens 1+1\tens Q_i+\sum_{j=1}^i\sum_{E_1+E_2=[i-j+1,i-1]}\rho^jQ_{i-j}Q(E_1)\tens Q_{i-j}Q(E_2)$.
	\end{itemize}
\end{lemma}
\end{samepage}

It is also easy to prove that the subalgebra of $\scr A^\aast$ generated by $\rho$ and $Q_i$, $i\geq 0$, is an exterior algebra in the $Q_i$'s over $\Z/\ell[\rho]\subset H^\aast$ (but the algebra generated by $Q_i$ and $H^\aast$, which is the left $H^\aast$-submodule generated by the operations $Q(E)$, is not even commutative if $\rho\neq 0$).

A well-known result in topology states that the left and right ideals of $\scr A^\ast$ generated by $Q_i$, $i\geq 0$, coincide and are generated by $Q_0$ as two-sided ideals. This can fail altogether in the motivic Steenrod algebra: for example, if $S$ is a field and $\rho\neq 0$, $Q_0\tau$ and $\tau Q_0$ are the unique nonzero elements of degree $(1,1)$ in the right and left ideals and $Q_0\tau-\tau Q_0=\rho$, so neither ideal is included in the other and in particular neither ideal is a two-sided ideal. It is true that the $H^\aast$-bimodules generated by those various ideals coincide, but this is not very useful.

\begin{lemma}\label{lem:leftideal}
	The left ideal of $\scr A^\aast$ generated by $\{Q_i\suchthat i\geq 0\}$ is the left $H^\aast$-submodule generated by $\{\rho(E,R)\suchthat E\neq \emptyset\}$.
\end{lemma}

\begin{proof}
	Define a matrix $a$ by the rule
	\[P^RQ(E)=\sum_{E',R'}a^{E,R}_{E',R'}\rho(E',R').\]
	Then $a^{E,R}_{E',R'}$ is the coefficient of $\xi(R)\tens\tau(E)$ in $\Delta(\tau(E')\xi(R'))$. The only term in $\Delta(\xi(R'))$ that can be a factor of $\xi(R)\tens\tau(E)$ is $\xi(R')\tens 1$, so we must have $R'\subset R$ for $a^{E,R}_{E',R'}$ to be nonzero. If $R=R'$, we must further have a term $1\tens\tau(E)$ in $\Delta(\tau(E'))$, and it is easy to see that this cannot happen unless also $E=E'$, in which case $\xi(R)\tens\tau(E)$ appears with coefficient $1$ in $\Delta(\tau(E)\xi(R))$. It is also clear that $a^{E,R}_{\emptyset,R'}=0$ if $E\neq \emptyset$. Combining these three facts, we can write
	\[P^RQ(E)=\rho(E,R)+\sum_{\substack{ E'\neq \emptyset \\ R'\subsetneqq R}}a^{E,R}_{E',R'} \rho(E',R').\]
	We can then use induction on the $\subset$-order of $R$ to prove that \emph{for all $E\neq\emptyset$, $\rho(E,R)$ is an $H^\aast$-linear combination of elements of the form $P^{R'}Q(E')$ with $E'\neq\emptyset$}. In particular, $\rho(E,R)$ is in the left ideal if $E\neq\emptyset$, which proves one inclusion.
	
	Conversely, let $\rho(E,R)Q_i$ be in the left ideal. Given what was just proved this is an $H^\aast$-linear combination of elements of the form $P^{R'}Q(E')Q_i$; because the $Q_i$'s generate an exterior algebra, such an element is either $0$ or $\pm P^{R'}Q(E'')$ with $E''\neq\emptyset$. The above formula shows directly that this is in turn an $H^\aast$-linear combination of elements of the desired form.
\end{proof}

We will denote by $\scr P_\aast$ the left $H_\aast$-submodule of $\scr A_\aast$ generated by the elements $\xi(R)$; it is clearly a left $\scr A_\aast$-comodule algebra (but it is not a Hopf algebroid in general, since it may not even be a right $H_\aast$-module). As an $H_\aast$-algebra it is the polynomial ring $H_\aast[\xi_1,\xi_2,\dotsc]$.

\begin{corollary}\label{cor:P**dual}
	The inclusion $\scr P_\aast\into\scr A_\aast$ is dual to the projection $\scr A^\aast\to\scr A^\aast/\scr A^\aast(Q_0,Q_1,\dotsc)$.
\end{corollary}

\begin{proof}
	Follows at once from Lemma~\ref{lem:leftideal}.
\end{proof}

\subsection{The motive of \texorpdfstring{$H\Z$}{𝐻𝐙} with finite coefficients}

Denote by $\scr M$ the basis of $\scr A_\aast$ formed by the elements $\tau(E)\xi(R)$. Since $\scr A_\aast=\pi_\aast(H\wedge H)$, $\scr M$ defines a map of $H$-modules
\begin{equation}\label{eqn:HH}
	\bigvee_{\zeta\in\scr M}\Sigma^{\abs\zeta}H\to H\wedge H
\end{equation}
which is an equivalence as $H\wedge H$ is a cellular $H$-module (Theorem~\ref{thm:HsmashH}).

Let $B\colon H\to \Sigma^{1,0}H\Z$ be the Bockstein morphism defined by the cofiber sequence~\eqref{eqn:bockstein1}. This cofiber sequence induces the short exact sequence
\begin{tikzmath}
	\def\colsep{1.2em}
	\diagram{0 & {H_\aast H\Z} & {H_\aast H} & {H_\aast \Sigma^{1,0}H\Z} & 0\rlap. \\};
	\arrows (11-) edge (-12) (12-) edge (-13) (13-) edge  node[above]{$B_\ast$} (-14) (14-) edge (-15);
\end{tikzmath}
Since $\beta$ is the composition of $B$ and the projection $H\Z\to H$, it shows that $H_\aast H\Z\cong\ker(\beta_\ast)$, and since $\beta$ is dual to $\tau_0$, this kernel is the $H_\aast$-submodule of $\scr A_\aast$ generated by the elements $\tau(E)\xi(R)$ with $\epsilon_0=0$. Denote by $\scr M_\Z\subset\scr M$ the set of those basis elements.

\begin{theorem}\label{thm:intsteenrod}
	The map
	\[\bigvee_{\zeta\in\scr M_\Z}\Sigma^{\abs\zeta}H\to H\wedge H\Z\]
	is an equivalence of $H$-modules.
\end{theorem}

\begin{proof}
	 In $\D(H)$, we have a commutative diagram
\begin{tikzmath}
	\diagram{
	\bigvee_{\zeta\in\scr M_\Z}\Sigma^{\abs\zeta} H & \bigvee_{\zeta\in\scr M}\Sigma^{\abs\zeta} H & \bigvee_{\zeta\in\scr M\minus \scr M_\Z}\Sigma^{\abs\zeta} H \\
	H\wedge H\Z & H\wedge H & H\wedge \Sigma^{1,0}H\Z \\
	};
	\arrows
	(21-) edge (-22) (22-) edge node[above]{$B$} (-23)
	(11-) edge (-12) (12-) edge (-13)
	(11) edge node[left]{$\alpha$} (21) (12) edge node[left]{$\simeq$} node[right]{\eqref{eqn:HH}} (22) (13) edge[dashed] node[right]{$\gamma$} (23);
\end{tikzmath}
    in which both rows are split cofiber sequences and $\alpha$ is to be proved an equivalence.
    First we show that the diagram can be completed by an arrow $\gamma$ as indicated. Let $\gamma$ be
	\begin{tikzmath}
		\diagram{\bigvee_{\zeta\in\scr M\minus \scr M_\Z}\Sigma^{\abs\zeta}H = \bigvee_{\zeta\in\scr M_\Z}\Sigma^{1,0}\Sigma^{\abs\zeta}H & H\wedge\Sigma^{1,0}H\Z\rlap, \\};
		\arrows (11-) edge node[above]{$\Sigma^{1,0}\alpha$} (-12);
	\end{tikzmath}
    where the equality is a reindexing. The commutativity of the second square is obvious. Applying $\spi_{\ast\ast}$ to this diagram, we deduce first that $\spi_{\ast\ast}(\alpha)$ is a monomorphism, whence that $\spi_{\ast\ast}(\gamma)$ is a monomorphism, and finally, using the five lemma, that $\spi_{\ast\ast}(\alpha)$ is an isomorphism.
\end{proof}

\section{The motivic cohomology of chromatic quotients of \texorpdfstring{$\MGL$}{MGL}}

In this section we compute the mod $\ell$ motivic cohomology of ``chromatic'' quotients of algebraic cobordism as modules over the motivic Steenrod algebra. The methods we use are elementary and work equally well to compute the ordinary mod $\ell$ cohomology of the analogous quotients of complex cobordism, such as connective Morava $K$-theory, at least if $\ell$ is odd (if $\ell=2$ the topological Steenrod algebra has a different structure and some modifications are required). The motivic computations require a little more care, however, mainly because the base scheme has plenty of nonzero cohomology groups (even if it is a field).

Throughout this section the base scheme $S$ is essentially smooth over a field.

\subsection{The Hurewicz map for \texorpdfstring{$\MGL$}{MGL}}
\label{sub:hurewicz}

Let $E$ be an oriented ring spectrum. We briefly review some standard computations from \cite[\S3–4]{Vezzosi:2001} and \cite[\S6]{Naumann:2009}. If $\BGL_r$ is the infinite Grassmannian of $r$-planes, we have
\begin{equation}\label{eqn:EBGL}
	E^\aast\BGL_{r}\cong E^\aast[[c_1,\dotsc,c_r]],
\end{equation}
where $c_i$ is the $i$th Chern class of the tautological vector bundle. This computation is obtained in the limit from the computation of the cohomology of the finite Grassmannian $\Gr(r,n)$, which is a free $E^\aast$-module of rank $\binom{n}{r}$. From \cite[Theorem~A.1]{Hu:2006} or \cite[Théorème 2.2]{Riou:2005}, we know that $\Sigma^\infty\Gr(r,n)_+$ is strongly dualizable in $\SH(S)$. If $X$ is the dual, then the canonical map
\[E_\aast X\tens_{E_\aast}E_\aast Y\to E_\aast(X\wedge Y)\]
is a natural transformation between homological functors of $Y$ that preserve direct sums, so it is an isomorphism for any cellular $Y$. It follows that the strong duality between $\Sigma^\infty\Gr(r,n)_+$ and $X$ induces a strong duality between $E_\aast\Gr(r,n)$ and $E_\aast X=E^{-\ast,-\ast}\Gr(r,n)$. In the limit we obtain canonical isomorphisms
\begin{gather*}
	E^\aast \BGL_{r}\cong \Hom_{E^\aast}(E_{-\ast,-\ast}\BGL_{r},E^\aast),\\
	E_\aast\BGL_{r}\cong\Hom_{E_\aast,c}(E^{-\ast,-\ast}\BGL_{r},E_\aast),
\end{gather*}
where $\Hom_{E_\aast,c}$ denotes continuous maps for the inverse limit topology on $E^\aast\BGL_{r}$. Taking further the limit as $r\to\infty$, we get duality isomorphisms for $\BGL$. There are Künneth formulas for finite products of such spaces. Now $\BGL$ has a multiplication $\BGL\times\BGL\to\BGL$ ``classifying'' the direct sum of vector bundles, which makes $E^\aast\BGL$ into a Hopf algebra over $E^\aast$. From the formula giving the total Chern class of a direct sum of vector bundles we obtain the formula for the comultiplication $\Delta$ on $E^\aast\BGL\cong E^\aast[[c_1,c_2,\dotsc]]$:
\[\Delta(c_n)=\sum_{i+j=n}c_i\tens c_j.\]
If $\beta_n\in E_\aast\BGL$ denotes the element which is dual to $c_1^n$ with respect to the monomial basis of $E^\aast\BGL$, then $\beta_0=1$ and it follows from purely algebraic considerations that the dual $E_\aast$-algebra $E_\aast\BGL$ is a polynomial algebra on the elements $\beta_n$ for $n\geq 1$ (see for example \cite[p.~176]{Madsen:1979}). Since the restriction map $E^\aast\BGL\to E^\aast \P^\infty$ simply kills all higher Chern classes, the $\beta_n$'s span $E_\aast\P^\infty\subset E_\aast\BGL$.

The multiplication $\MGL_r\wedge\MGL_s\to\MGL_{r+s}$ is compatible with the multiplication $\BGL_r\times\BGL_s\to\BGL_{r+s}$ under the Thom isomorphisms, and so the dual Thom isomorphism $E_\aast\BGL\cong E_\aast\MGL$ is an isomorphism of $E_\aast$-algebras. Thus, $E_\aast\MGL$ is also a polynomial algebra
\[E_\aast\MGL=E_\aast[b_1,b_2,\dotsc],\]
where $b_n\in \tilde E_{2n,n}\Sigma^{-2,-1}\MGL_1$ is dual to the image of $c_1^n$ by the Thom isomorphism \[E^\aast\P^\infty\cong \tilde E^{\aast}\Sigma^{-2,-1}\MGL_1.\]

The case $r=1$ of the computation~\eqref{eqn:EBGL} shows that we can associate to $E$ a unique graded formal group law $F_E$ over $E_{(2,1)\ast}$ such that, for any pair of line bundles $\scr L$ and $\scr M$ over $X\in\Sm/S$,
\[c_1(\scr L\tens\scr M)=F_E(c_1(\scr L),c_1(\scr M)).\]
The elements $b_n\in\MGL_\aast\MGL$ are then the coefficients of the power series defining the strict isomorphism between the two formal group laws coming from the two obvious orientations of $\MGL\wedge\MGL$. As the stack of formal group laws and strict isomorphisms is represented by a graded Hopf algebroid $(L,LB)$, where $L$ is the Lazard ring and $LB=L[b_1,b_2,\dotsc]$, we obtain a graded map of Hopf algebroids \[(L,LB)\to(\MGL_{(2,1)\ast},\MGL_{(2,1)\ast}\MGL)\] sending $b_n$ to $b_n\in \MGL_{2n,n}\MGL$ (see \cite[Corollary 6.7]{Naumann:2009}). We will often implicitly view elements of $L$ as elements of $\MGL_\aast$ through this map, and for $x\in L_n$ we simply write $\abs x$ for the bidegree $(2n,n)$.

Recall from~\S\ref{sub:motivic} that $H R$ is an oriented ring spectrum such that $H R_\aast$ carries the additive formal group law (since $[\scr L\tens\scr M]=[\scr L]+[\scr M]$ in the Picard group). It follows that we have a commutative square
\begin{tikzequation}\label{eqn:hurewicz}
	\def\colsep{2em}
	\diagram{L & R[b_1,b_2,\dotsc] \\ \MGL_\aast & HR_\aast\MGL \\};
	\arrows (11-) edge node[above]{$h_R$} (-12) (11) edge (21) (12) edge (22) (21-) edge (-22);
\end{tikzequation}
where the horizontal maps are induced by the right units of the respective Hopf algebroids. Explicitly, the map $h_R$ classifies the formal group law on $R[b_1,b_2,\dotsc]$ which is isomorphic to the additive formal group law via the exponential $\sum_{n\geq 0}b_nx^{n+1}$.

Let $I\subset L$ and $J\subset\Z[b_1,b_2,\dotsc]$ be the ideals generated by the elements of positive degree. By Lazard's theorem (\cite[Lemma 7.9]{Adams:1974}), $h_\Z$ induces an injective map $(I/I^2)_n\into (J/J^2)_n\cong\Z$ whose range is $\ell\Z$ if $n+1$ is a power of a prime number $\ell$ and $\Z$ otherwise. If $a_n\in L_n$ is an arbitrary lift of a generator of $(I/I^2)_n$, it follows easily that $L$ is a polynomial ring on the elements $a_n$, $n\geq 1$.

\begin{samepage}
\begin{definition}\label{def:typical}
	Let $\ell$ be a prime number and $r\geq 0$. An element $v\in L_{\ell^r-1}$ is called \emph{$\ell$-typical} if
	\begin{enumerate}
		\item\label{def:typical:1} $h_{\Z/\ell}(v)=0$;
		\item\label{def:typical:2} $h_{\Z/\ell^2}(v)\not\equiv 0$ modulo decomposables.
	\end{enumerate}
\end{definition}
\end{samepage}

For every $r\geq 0$, there is a canonical $\ell$-typical element in $L_{\ell^r-1}$, namely the coefficient of $x^{\ell^r}$ in the $\ell$-series of the universal formal group law. Indeed, since the formal group law classified by $h_{\Z/\ell}$ is isomorphic to the additive one, its $\ell$-series is zero, so condition \itemref{def:typical:1} holds. For $r\geq 1$, one can show that the image of that coefficient in $(I/I^2)_{\ell^r-1}$ generates a subgroup of index prime to $\ell$ (see for instance \cite[Lecture 13]{Lurie:2010}); since $h_\Z$ identifies $(I/I^2)_{\ell^r-1}$ with a subgroup of $(J/J^2)_{\ell^r-1}$ of index exactly $\ell$, condition \itemref{def:typical:2} holds.

\begin{remark}\label{rmk:lazard}
	Lazard's theorem shows in particular that $h_\Z\colon L\to\Z[b_1,b_2,\dotsc]$ is injective. It follows from the commutative square~\eqref{eqn:hurewicz} that $L\to \MGL_\aast$ is injective if $S$ is not empty. A consequence of the Hopkins–Morel equivalence is that the image of $L$ is precisely $\MGL_{(2,1)\ast}$ if $S$ is a field of characteristic zero (see Proposition~\ref{prop:piMGL}).
\end{remark}

\subsection{Regular quotients of \texorpdfstring{$\MGL$}{MGL}}
\label{sub:regular}

From now on we fix a prime number $\ell\neq\Char S$. We abbreviate $H\Z/\ell$ to $H$ and $h_{\Z/\ell}$ to $h$. By Lazard's theorem, $h(L)\subset\Z/\ell[b_1,b_2,\dotsc]$ is a polynomial subring $\Z/\ell[b_n'\suchthat n\neq \ell^r-1]$ where $b_n'\equiv b_n$ modulo decomposables. We choose once and for all a graded ring map
\[\pi\colon \Z/\ell[b_1,b_2,\dotsc]=\Z/\ell[b_1',b_2',\dotsc]\to h(L)\]
which is a retraction of the inclusion. For example, we could specify $\pi(b_{\ell^r-1})=0$ for all $r\geq 1$, but our arguments will work for any choice of $\pi$ and none seems particularly canonical.

\begin{theorem}\label{thm:HMGL}
	The coaction $\Delta\colon H_\aast\MGL\to\scr A_\aast\tens_{H_\aast}H_\aast\MGL$ factors through $\scr P_\aast\tens\Z/\ell[b_1,b_2,\dotsc]$
	and the composition
	\begin{tikzmath}
		\diagram{
			H_\aast\MGL & \scr P_\aast\tens\Z/\ell[b_1,b_2,\dotsc] & \scr P_\aast\tens h(L) \\
		};
		\arrows (11-) edge node[above]{$\Delta$} (-12) (12-) edge node[above]{$\id\tens\pi$} (-13);
	\end{tikzmath}
	is an isomorphism of left $\scr A_\aast$-comodule algebras.
\end{theorem}

Towards proving this theorem we explicitly compute the coaction $\Delta$ of $\scr A_\aast$ on $H_\aast\MGL$.
Since it is an $H_\aast$-algebra map, it suffices to compute $\Delta(b_n)$ for $n\geq 1$.
Consider the zero section
\[s\colon \P^\infty_+\to\MGL_1\]
as a map in $\H_\pt(S)$. In cohomology this map is the composition of the Thom isomorphism and multiplication by the top Chern class $c_1$. In homology, it therefore sends $\beta_n$ to $0$ if $n=0$ and to $b_{n-1}$ otherwise. Thus,
\begin{equation}\label{eqn:zerosection}
	\Delta(b_n)=\Delta(s_\ast(\beta_{n+1}))=(1\tens s_\ast)\Delta(\beta_{n+1}).
\end{equation}
The action of $\scr A^\aast$ on $c_1^n\in H^\aast\P^\infty=H^\aast[c_1]$ is determined by the Cartan formulas (Lemma~\ref{lem:cartan}). For degree reasons $Q_i$ acts trivially on elements in $H^{(2,1)\ast}\P^\infty$, and we get
\[P^R(c_1^n)=a_{n,R}c_1^{n+\abs R},\quad Q_i(c_1^n)=0,\]
where $\abs R=\sum_{i\geq 1} r_i(\ell^i-1)$ and $a_{n,R}$ is the multinomial coefficient given by
\[a_{n,R}=\binom{n}{n-\sum_{i\geq 1}r_i,r_1,r_2,\dotsc}\]
(understood to be $0$ if $\sum_{i\geq 1}r_i>n$).
Dualizing, we obtain
\[\Delta(\beta_n)=\sum_{m+\abs R=n}a_{m,R}\xi(R)\tens\beta_{m},\]
whence by~\eqref{eqn:zerosection},
\begin{equation}\label{eqn:HMGLcomod}
	\Delta(b_n)=\sum_{m+\abs R=n}a_{m+1,R}\xi(R)\tens b_{m}.
\end{equation}

\begin{lemma}\label{lem:HMGL1}
	The $H_\aast$-algebra map $f\colon H_\aast\MGL\to\scr P_\aast$ defined by
	\[f(b_n)=\begin{cases}
		\xi_r & \text{if $n=\ell^r-1$,}\\
		0 & \text{otherwise.}
	\end{cases}\]
	is a map of left $\scr A_\aast$-comodules.
\end{lemma}

\begin{proof}
	If $m$ is of the form $\ell^r-1$, then the coefficient $a_{m+1,R}$ vanishes mod $\ell$ unless $R=\ell^re_i$ for some $i\geq 0$, in which case $a_{m+1,R}=1$ and $m+\abs R=\ell^{r+i}-1$. Comparing~\eqref{eqn:HMGLcomod} with the formula for $\Delta(\xi_r)$ shows that $f$ is a comodule map.
\end{proof}

\begin{proof}[Proof of Theorem~\ref{thm:HMGL}]
	The formula~\eqref{eqn:HMGLcomod} shows that $\Delta$ factors through $\scr P_\aast\tens\Z/\ell[b_1,b_2,\dotsc]$. Let $g$ be the map to be proved an isomorphism. Note that it is a comodule algebra map since $\Delta$ is and $\pi$ is a ring map. Formula~\eqref{eqn:HMGLcomod} shows further that $g$ is extended from a map
	\[\tilde g\colon \Z/\ell[b_1,b_2,\dotsc]\to\Z/\ell[\xi_1,\xi_2,\dotsc]\tens h(L).\]
	If $n+1$ is not a power of $\ell$, we have
	\[g(b_n')\equiv 1\tens b_n'\]
	modulo decomposables by definition of $\pi$, whereas for $r\geq 0$ we have, by Lemma~\ref{lem:HMGL1},
	\[g(b_{\ell^r-1})\equiv \xi_r\tens 1\]
	modulo decomposables. These congruences show that $\tilde g$ is surjective. Now $\tilde g$ is a map between $\Z/\ell$-modules of the same finite dimension in each bidegree, so $\tilde g$ and hence $g$ are isomorphisms. 
\end{proof}

There is no reason to expect the isomorphism $g$ of Theorem~\ref{thm:HMGL} to be a map of $L$-modules (as $\Delta$ clearly is not), but it is easy to modify it so that it preserves the $L$-module structure as well. To do this consider the $H_\aast$-algebra map
\[\tilde f\colon\scr P_\aast\to H_\aast\MGL,\quad \tilde f(\xi_r)=g^{-1}(\xi_r\tens 1),\]
which is clearly a map of $\scr A_\aast$-comodule algebras.

\begin{corollary}\label{cor:HMGL}
	The map $\tilde f$ and the inclusion of $h(L)$ induce an isomorphism
	\[\scr P_\aast\tens h(L)\cong H_\aast\MGL\]
	of left $\scr A_\aast$-comodule algebras and of $L$-modules.
\end{corollary}

\begin{proof}
	We have $g^{-1}(\xi_r\tens 1)\equiv b_{\ell^r-1}$ modulo decomposables. It follows that the map is surjective and hence, as in the proof of Theorem~\ref{thm:HMGL}, an isomorphism.
\end{proof}

Now is a good time to recall the construction of general quotients of $\MGL$. Given an $\MGL$-module $E$ and a family $(x_i)_{i\in I}$ of homogeneous elements of $\pi_\aast\MGL$, the quotient $E/(x_i)_{i\in I}$ is defined by
\[
E/(x_i)_{i\in I}=E\wedge_\MGL\hocolim_{\{i_1,\dotsc,i_k\}\subset I}(\MGL/x_{i_1}\wedge_\MGL\dotsb\wedge_\MGL\MGL/x_{i_k}),
\]
where the homotopy colimit is taken over the filtered poset of finite subsets of $I$ and $\MGL/x$ is the cofiber of $x\colon \Sigma^{\abs x}\MGL\to\MGL$. It is clear that $E/(x_i)_{i\in I}$ is invariant under permutations of the indexing set $I$.

Let $x\in L$ be a homogeneous element such that $h(x)$ is nonzero. Then multiplication by $h(x)$ is injective and so there is a short exact sequence
\[0\to H_\aast\Sigma^{\abs x}\MGL\to H_\aast\MGL\to H_\aast(\MGL/x)\to 0.\]
It follows that $H_\aast(\MGL/x)\cong H_\aast[b_1,b_2,\dotsc]/h(x)$ and, by comparison with the isomorphism of Corollary~\ref{cor:HMGL}, we deduce that the map
\[\scr P_\aast\tens h(L)/h(x)\to H_\aast(\MGL/x)\]
induced by $\tilde f$ and the inclusion is an isomorphism of left $\scr A_\aast$-comodules and of $L$-modules. Let us say that a (possibly infinite) sequence of homogeneous elements of $L$ is \emph{$h$-regular} if its image by $h$ is a regular sequence. By induction we then obtain the following result.

\begin{lemma}\label{lem:regular}
	Let $x$ be an $h$-regular sequence of homogeneous elements of $L$. Then the maps $\tilde f$ and $h(L)/h(x)\into H_\aast(\MGL/x)$ induce an isomorphism
	\[\scr P_\aast\tens h(L)/h(x)\cong H_\aast (\MGL/x)\]
	of left $\scr A_\aast$-comodules and of $L$-modules.
\end{lemma}

We can now ``undo'' the modification of Corollary~\ref{cor:HMGL}:

\begin{theorem}\label{thm:regular}
	Let $x$ be an $h$-regular sequence of homogeneous elements of $L$. Then the coaction of $\scr A_\aast$ on $H_\aast(\MGL/x)$ and the map $\pi$ induce an isomorphism
	\[H_\aast (\MGL/x) \cong \scr P_\aast\tens h(L)/h(x)\]
	of left $\scr A_\aast$-comodules.
\end{theorem}

\begin{proof}
	Let $\tilde g$ be the isomorphism of Lemma~\ref{lem:regular} and let $g$ be the map to be proved an isomorphism. Then $g\tilde g(1\tens b)\equiv 1\tens b$ modulo decomposables and $g\tilde g(\xi_r\tens 1)=\xi_r\tens 1$, so $g\tilde g$ and hence $g$ are isomorphisms by the usual argument.
\end{proof}

If $x$ is a maximal $h$-regular sequence in $L$, \ie, an $h$-regular sequence which generates the maximal ideal in $h(L)$, then, by Theorem~\ref{thm:regular}, the coaction of $\scr A_\aast$ on $H_\aast(\MGL/x)$ and the projection $\Z/\ell[b_1,b_2,\dotsc]\to\Z/\ell$ induce an isomorphism
\[H_\aast(\MGL/x)\cong\scr P_\aast\]
of left $\scr A_\aast$-comodules. Note that this isomorphism does not depend on the choice of $\pi$ anymore and is therefore canonical.

As we noted in \S\ref{sub:duality}, $H\wedge\MGL$ is a psf $H$-module. Dualizing Theorem~\ref{thm:HMGL} and using Corollary~\ref{cor:P**dual}, we deduce that the map
\[\scr A^\aast/\scr A^\aast(Q_0,Q_1,\dotsc)\tens h(L)^\vee\to H^\aast\MGL,\quad [\phi]\tens m\mapsto\phi(m),\]
is an isomorphism of left $\scr A^\aast$-module coalgebras. Here the inclusion $h(L)^\vee\into H^\aast\MGL$ is dual to $\pi$.
If $x$ is an $h$-regular sequence in $L$, the computation of $H_\aast(\MGL/x)$ shows, with Lemma~\ref{lem:split}, that $H\wedge \MGL/x$ is a split direct summand of $H\wedge\MGL$, so it is also psf. By dualizing Theorem~\ref{thm:regular}, we obtain a computation of $H^\aast(\MGL/x)$. For example, if $x$ is a maximal $h$-regular sequence, we obtain that the map
\[\scr A^\aast/\scr A^\aast(Q_0,Q_1,\dotsc)\to H^\aast(\MGL/x),\quad [\phi]\mapsto \phi(\th),\]
where $\th\colon\MGL/x\to H$ is the lift of the Thom class, is an isomorphism of left $\scr A^\aast$-modules.

\subsection{Key lemmas}
\label{sec:key}

Recall that $\th\colon \MGL\to H$ is the universal Thom class.

\begin{lemma}\label{lem:PRtau}
	Let $R=(r_1,r_2,\dotsc)$. Then $P^R(\th)\in H^\aast\MGL$ is dual to $\prod_{i\geq 1}b_{\ell^i-1}^{r_i}\in H_\aast\MGL$.
\end{lemma}

\begin{proof}
	As $\th$ is dual to $1$, we must look for monomials $m\in \Z/\ell[b_1,b_2,\dotsc]$ such that $\Delta(m)$ has a term of the form $\xi(R)\tens 1$. By Lemma~\ref{lem:HMGL1}, such a term can only appear if $m$ is a monomial in $b_{\ell^i-1}$, and this monomial must be $\prod_{i\geq 1}b_{\ell^i-1}^{r_i}$.
\end{proof}

\begin{lemma}\label{lem:key2}
	Let $r\geq 0$ and $n=\ell^r-1$. Let $v\in L_n$ be an $\ell$-typical element, $\th'\colon \MGL/v\to H$ the unique lift of the universal Thom class, and $\delta$ the connecting morphism in the cofiber sequence
	\begin{tikzmath}
		\def\colsep{1.2em}
		\diagram{
		\Sigma^{2n,n}\MGL  & \MGL & \MGL/v & \Sigma^{2n+1,n}\MGL\rlap. \\
		};
		\arrows{
		(11-) edge node[above]{$v$} (-12) (12-) edge (-13) (13-) edge node[above]{$\delta$} (-14)
		};
	\end{tikzmath}
	 Then the square
\begin{tikzmath}
	\diagram{
	\MGL/v & \Sigma^{2n+1,n}\MGL \\
	H & \Sigma^{2n+1,n}H \\
	};
	\arrows
	(11-) edge node[above]{$\delta$} (-12)
	(11) edge node[left]{$\th'$} (21)
	(12) edge node[right]{$\Sigma^{2n+1,n}\th$} (22)
	(21-) edge node[below]{$Q_r$} (-22);
\end{tikzmath}
	commutes up to multiplication by an element of $\Z/\ell^\times$.
\end{lemma}

\begin{proof}
	We may clearly assume that $S$ is connected, so that $H^{0,0}\MGL\cong \Z/\ell$ with $\th$ corresponding to $1$.
	Since $H^{2n+1,n}\MGL=0$, $\delta^\ast\colon H^{0,0}\MGL\to H^{2n+1,n}(\MGL/v)$ is surjective, so it will suffice to show that $Q_r\th'\neq 0$. Recall that $Q_0=\beta$ and that, for $r\geq 1$, $Q_r=q_r\beta-\beta q_r$. In the latter case we have $\beta\th'=0$ for degree reasons. In all cases we must therefore show that $\beta q_r\th'\neq 0$. Consider the diagram of exact sequences
	\begin{tikzmath}
		\diagram{
		H^{2n,n}(\MGL/v,\Z/\ell^2) & H^{2n,n}(\MGL/v) & H^{2n+1,n}(\MGL/v) \\
		H^{2n,n}(\MGL,\Z/\ell^2) & H^{2n,n}\MGL & H^{2n+1,n}\MGL \\
		H^{0,0}(\MGL,\Z/\ell^2)\rlap. & & \\
		};
		\arrows
		(11-) edge (-12) (12-) edge node[above]{$\beta$} (-13)
		(21-) edge[->] (-22) (22-) edge node[above]{$\beta$} (-23)
		(11) edge (21) (12) edge (22) (13) edge (23)
		(21) edge node[right]{$v^\ast$} (31);
	\end{tikzmath}
		By Lemma~\ref{lem:PRtau}, $q_r\th$ is dual to $b_n$. Thus, an arbitrary lift of $q_r\th$ to $H^{2n,n}(\MGL,\Z/\ell^2)$  has the form $x+\ell y$ where $x$ is dual to $b_n$ and $y$ is any cohomology class, and we must show that $v^\ast(x+\ell y)\neq 0$. By duality,
		\[\langle v^\ast(x+\ell y),1\rangle =\langle x+\ell y,v_\ast(1)\rangle=\langle x+\ell y,h_{\Z/\ell^2}(v)\rangle.\]
		 Since $v$ is $\ell$-typical, $\langle \ell y,h_{\Z/\ell^2}(v)\rangle=0$ and $\langle x,h_{\Z/\ell^2}(v)\rangle\neq 0$, so $v^\ast(x+\ell y)\neq 0$.
\end{proof}

\begin{lemma}\label{lem:dim0}
	Assume that $S$ is the spectrum of a field. Let $v\in L$ be a homogeneous element and let $p\colon \MGL\to\MGL/v$ be the projection. If $x\in H^{2n,n}(\MGL/v)$ is such that $p^\ast(x)=0$, then $\beta(x)=0$.
\end{lemma}

\begin{proof}
	Since $p^\ast(x)=0$, $x=\delta^\ast(y)$ for some $y$, where $\delta\colon \MGL/v\to\Sigma^{1,0}\Sigma^{\abs v}\MGL$. Since $S$ is the spectrum of a field, we must have $y=\lambda z$ for some $\lambda\in H^{1,1}(S,\Z/\ell)$ and $z\in H^{(2,1)\ast}\MGL$. But then $\beta(y)=\beta(\lambda)z-\lambda\beta(z)=0$, whence $\beta(x)=0$.
\end{proof}

\begin{lemma}\label{lem:key2.5}
	Assume that $S$ is the spectrum of a field. Let $x\in H^{2n,n}\MGL$ and let $x'$ be a lift of $x$ in $H^{2n,n}(\MGL/\ell)$. Then the square
	\begin{tikzmath}
		\diagram{
		\MGL/\ell & \Sigma^{1,0}\MGL \\
		\Sigma^{2n,n}H & \Sigma^{2n+1,n}H \\
		};
		\arrows
		(11-) edge node[above]{$\delta$} (-12)
		(11) edge node[left]{$x'$} (21)
		(12) edge node[right]{$\Sigma^{1,0}x$} (22)
		(21-) edge node[below]{$\beta$} (-22);
	\end{tikzmath}
	is commutative.
\end{lemma}

\begin{proof}
	Since $\beta(x)=0$, $x$ lifts to $\hat x\in H^{2n,n}(\MGL,\Z/\ell^2)$ and it is clear that the square
	\begin{tikzmath}
		\diagram{
		\MGL & \MGL \\
		\Sigma^{2n,n}H & \Sigma^{2n,n}H\Z/\ell^2 \\
		};
		\arrows
		(11-) edge node[above]{$\ell$} (-12)
		(11) edge node[left]{$x$} (21)
		(12) edge node[right]{$\hat x$} (22)
		(21-) edge node[below]{$\ell$} (-22);
	\end{tikzmath}
	commutes, so there exists $y\colon \MGL/\ell\to \Sigma^{2n,n}H$ such that $p^\ast(y)=x$ and $\beta(y)=\delta^\ast(x)$. In particular, we have $p^\ast(x')=p^\ast(y)$. By Lemma~\ref{lem:dim0}, we obtain $\beta(x')=\beta(y)=\delta^\ast(x)$.
\end{proof}

\begin{lemma}\label{lem:key3}
	Let $i\geq 0$ and let $v\in L$ be a homogeneous element such that the coefficient of $b_{\ell^i-1}$ in $h_{\Z/\ell^2}(v)$ is zero. If $\th'\colon \MGL/v\to H$ lifts the universal Thom class, then $Q_i\th'=0$.
\end{lemma}

\begin{proof}
	Since $Q_i$ and $\th$ are compatible with changes of essentially smooth base schemes over fields, we may assume without loss of generality that $S$ is a field. We first consider the case $\abs v=0$. Then $v$ is an integer which must be divisible by $\ell$. If moreover $i=0$, then we assumed $v$ divisible by $\ell^2$ so that $\beta\th'=0$. If $i\geq 1$, we may assume that $v=\ell$ since $\th'$ factors through $\MGL/\ell$. By two applications of Lemma~\ref{lem:key2.5}, we have $\beta q_i\th'=\delta^\ast q_i\th=q_i\delta^\ast\th=q_i\beta\th'$, whence $Q_i\th'=q_i\beta\th'-\beta q_i\th'=0$.
	
	Suppose now that $\abs v\geq 1$. Then $\beta\th'=0$ for degree reasons, so we can assume $i\geq 1$ and we must show that $\beta q_i\th'=0$. Let $n=\ell^i-1$ and consider the diagram of exact sequences
\begin{tikzmath}
	\diagram{
	H^{2n,n}(\MGL/v,\Z/\ell^2) & H^{2n,n}(\MGL/v) & H^{2n+1,n}(\MGL/v) \\
	H^{2n,n}(\MGL,\Z/\ell^2) & H^{2n,n}\MGL & H^{2n+1,n}\MGL \\
	H^{2n,n}(\Sigma^{\abs{v}}\MGL,\Z/\ell^2)\rlap. & & \\
	};
	\arrows
	(11-) edge (-12) (12-) edge node[above]{$\beta$} (-13)
	(21-) edge[->] (-22) (22-) edge node[above]{$\beta$} (-23)
	(11) edge (21) (21) edge node[right]{$v^\ast$} (31)
	(12) edge[->] node[right]{$p^\ast$} (22)
	(13) edge[->] (23);
\end{tikzmath}
	By Lemma~\ref{lem:PRtau}, $q_i\th\in H^{2n,n}\MGL$ is dual to $b_n$. Let $x\in H^{2n,n}(\MGL,\Z/\ell^2)$ be dual to $b_n$, so that $x$ lifts $q_i\th$. For any $m\in H_\aast(\MGL,\Z/\ell^2)$, 
	\[\langle v^\ast(x),m\rangle=\langle x,v_\ast(m)\rangle=\langle x,h_{\Z/\ell^2}(v)m\rangle=0\]
	since $b_n$ is the only monomial with which $x$ pairs nontrivially. Thus, $v^\ast(x)=0$ and $x$ lifts to an element $\hat x\in H^{2n,n}(\MGL/v,\Z/\ell^2)$. Let $y$ be the image of $\hat x$ in $H^{2n,n}(\MGL/v)$. Then $p^\ast y=q_i\th=p^\ast q_i\th'$, and hence, by Lemma~\ref{lem:dim0}, $\beta q_i\th'=\beta y=0$.
\end{proof}

\subsection{Quotients of \texorpdfstring{$\BP$}{BP}}

\begin{lemma}\label{lem:modsplit}
	Let $E$ be an $\MGL$-module and let $x\in L$ be a homogeneous element such that $h(x)=0$. If $H\wedge E$ is psf, then $H\wedge E/x$ is psf.
\end{lemma}

\begin{proof}
	Since $h(x)=0$, we have a short exact sequence
	\[0\to H_\aast E\to H_\aast (E/x)\to H_\aast\Sigma^{1,0}\Sigma^{\abs x}E\to 0.\]
	This sequence splits in $H_\aast$-modules since the quotient is free, so we deduce from Lemma~\ref{lem:split} that
	\[H\wedge E/x\simeq (H\wedge E)\vee(H\wedge \Sigma^{1,0}\Sigma^{\abs x}E).\]
	Since $\abs x$ is of the form $(2n,n)$, it is clear that $H\wedge E/x$ is psf.
\end{proof}

\begin{lemma}\label{lem:projective}
	Let $E$ be an $\MGL$-module and let $x\in L$ be a homogeneous element such that $h(x)=0$. If $H_\aast E$ is projective over $H_\aast\MGL$, so is $H_\aast (E/x)$.
\end{lemma}

\begin{proof}
	The short exact sequence
	\[0\to H_\aast E\to H_\aast (E/x)\to H_\aast\Sigma^{1,0}\Sigma^{\abs x}E\to 0\]
	splits in $H_\aast\MGL$-modules.
\end{proof}

\begin{theorem}\label{thm:AmodQi}
	Let $M$ be a quotient of $\MGL$ by a maximal $h$-regular sequence, $I$ a set of nonnegative integers, and for each $i\in I$ let $v_i\in L_{\ell^i-1}$ be an $\ell$-typical element. Then there is an isomorphism of left $\scr A^\aast$-modules
	\[\scr A^\aast/\scr A^\aast(Q_i\suchthat i\notin I)\cong H^\aast(M/(v_i\suchthat i\in I))\]
	given by $[\phi]\mapsto \phi(\th)$, where $\th\colon M/(v_i\suchthat i\in I)\to H$ is the lift of the universal Thom class.
\end{theorem}

\begin{proof}
	We may clearly assume that $I$ is finite, and we proceed by induction on the size of $I$. If $I$ is empty, the theorem is true by Theorem~\ref{thm:regular}. Suppose it is true for $I$ and let $r\notin I$. Let $E=M/(v_i\suchthat i\in I)$. Since $E$ is a cellular $\MGL$-module and $H_\aast (\MGL/v_r)$ is flat over $H_\aast\MGL$ by Lemma~\ref{lem:projective}, the canonical map
	\[H_\aast E\tens_{H_\aast\MGL}H_\aast(\MGL/v_r)\to H_\aast(E\wedge_{\MGL}\MGL/v_r)=H_\aast(E/v_r)\]
	is an isomorphism. By Lemma~\ref{lem:modsplit}, all the $H$-modules $H\wedge\MGL$, $H\wedge E$, $H\wedge \MGL/v_r$, and $H\wedge E/v_r$ are psf. If we dualize this isomorphism and use Proposition~\ref{prop:psf}, we get an isomorphism
	\[H^\aast (E/v_r)\cong H^\aast E\cotens_{H^\aast\MGL}H^\aast(\MGL/v_r),\]
	where $\th$ on the left-hand side corresponds to $\th\tens\th$ on the right-hand side and the $\scr A^\aast$-module structure on the right-hand side is given by $\phi\cdot(x\tens y)=\Delta(\phi)(x\tens y)$. By the Cartan formula (Lemma~\ref{lem:cartan}),
	\[\Delta(Q_i)(\th\tens\th)=Q_i\th\tens\th+\th\tens Q_i\th +\sum_{j= 1}^i\psi_j(Q_{i-j}\th\tens Q_{i-j}\th)\]
	for some $\psi_j\in\scr A^\aast\tens_{H^\aast}\scr A^\aast$.
	By Lemma~\ref{lem:key3}, $Q_i\th\in H^\aast(\MGL/v_r)$ is zero when $i\neq r$. By induction hypothesis, $\Delta(Q_r)(\th\tens\th)=\th\tens Q_r\th$ and, if $i\notin I\cup\{r\}$, $\Delta(Q_i)(\th\tens\th)=0$. The latter shows that the map $[\phi]\mapsto \Delta(\phi)(\th\tens\th)$ is well-defined. We can thus form a diagram of short exact sequences of left $\scr A^\aast$-modules
\begin{tikzmath}
	\diagram{\scr A^\aast/\scr A^\aast(Q_i\suchthat i\notin I) & \scr A^\aast/\scr A^\aast(Q_i\suchthat i\notin I\cup\{r\}) & \scr A^\aast/\scr A^\aast(Q_i\suchthat i\notin I) \\
	H^\aast E\cotens H^\aast\MGL & H^\aast E\cotens H^\aast(\MGL/v_r) & H^\aast E\cotens H^\aast\MGL \\};
	\arrows
	(11-) edge[c->] node[above]{$\cdot Q_r$} (-12) (12-) edge[->>] (-13)
	(21-) edge[c->] node[above]{$1\cotens \delta^\ast$} (-22) (22-) edge[->>] (-23)
	(11) edge node[right]{$\cong$} (21)
	(12) edge (22)
	(13) edge node[right]{$\cong$} (23);
\end{tikzmath}
(where the cotensor products are over $H^\aast\MGL$).
	The right square commutes because the image of $1\in\scr A^\aast/\scr A^\aast(Q_i\suchthat i\notin I\cup\{r\})$ in the bottom right corner is $\th\tens\th$ either way. For the left square, the two images of $1$ are $\th\tens\delta^\ast\th$ and $\Delta(Q_r)(\th\tens\th)=\th\tens Q_r\th$. Lemma~\ref{lem:key2} then shows that the left square commutes up to multiplication by an element of $\Z/\ell^\times$, so the map between extensions is an isomorphism by the five lemma.
\end{proof}

\begin{remark}
	The title of this paragraph is justified by the fact that the motivic Brown–Peterson spectrum at $\ell$, which can be defined using the Cartier idempotent (as in \cite[\S5]{Vezzosi:2001}) or the motivic Landweber exact functor theorem (\cite[Theorem 8.7]{Naumann:2009}), is equivalent as an $\MGL$-module to $\MGL_{(\ell)}/x$ where $x$ is any regular sequence in $L$ that generates the vanishing ideal for $\ell$-typical formal group laws. This is in fact a consequence of the Hopkins–Morel equivalence (see Example~\ref{ex:landweber}). Since such a sequence $x$ is a maximal $h$-regular sequence, Theorem~\ref{thm:AmodQi} specializes \emph{a posteriori} to the computation of the mod $\ell$ cohomology of quotients of $\BP$.
\end{remark}

\begin{remark}
	If we forget the identification of $\scr A^\aast$ provided by Lemma~\ref{lem:steenrod}, all the results in this section remain true as long as we replace $\scr A^\aast$ by its subalgebra of standard Steenrod operations. In topology, the equivalence $M/(v_0,v_1,\dotsc)\simeq H\Z/\ell$ can be proved without knowledge of the Steenrod algebra, and so the topological version of Theorem~\ref{thm:AmodQi} (with $I=\{0,1,2,\dotsc\}$) gives a proof that the endomorphism algebra of the topological Eilenberg–Mac Lane spectrum $H\Z/\ell$ is generated by the reduced power operations and the Bockstein.
\end{remark}

\begin{example}
	As a counterexample to a conjecture one might make regarding the cohomology of more general quotients of $\MGL$ than those we considered, we observe that if $\ell=2$, $v\in L_3$ is $2$-typical, and $\hat x\in H^{4,2}(\MGL/v)$ is a lift of the dual $x$ to $b_2$, then $Q_1\hat x$ is nonzero. By~\eqref{eqn:HMGLcomod}, we have \[\Delta(b_3)=1\tens b_3+\xi_1^2\tens b_1+\xi_2\tens 1+\xi_1\tens b_2.\]
	By an analysis similar to those done in \S\ref{sec:key}, we deduce that for any lift $y$ of $q_1x$ to mod $4$ cohomology, $v^\ast(y)$ is nonzero, and hence that $Q_1\hat x$ is nonzero.
	This is only one instance of the following general phenomenon: if $\xi(R)\tens b_n$ is a term in $\Delta(b_{\ell^r-1})$ that does not correspond to a term in $\Delta(\xi_r)$, if $v\in L_{\ell^r-1}$ is $\ell$-typical, and if $\hat x\in H^\aast(\MGL/v)$ is a lift of the dual to $b_n$, then $\beta P^R(\hat x)$ is nonzero.
\end{example}

\section{The Hopkins–Morel equivalence}
\label{sec:towards_the_hopkins_morel_isomorphism}

In this section, $S$ is an essentially smooth scheme over a field (although in Lemmas~\ref{lem:HZ'effective}, \ref{lem:localhomology}, and \ref{lem:regular2} it can be arbitrary). Let $c$ denote the the characteristic exponent of $S$, \ie, $c=1$ if $\Char S=0$ and $c=\Char S$ otherwise.

\begin{definition}\label{def:generators}
	A set of generators $a_n\in L_n$, $n\geq 1$, will be called \emph{adequate} if, for every prime $\ell$ and every $r\geq 1$, $a_{\ell^r-1}$ is $\ell$-typical (Definition~\ref{def:typical}).
\end{definition}

Adequate sets of generators of $L$ exist: for example, the generators given in \cite[(7.5.1)]{Hazewinkel:1977} (where $m_n$ is the coefficient of $x^{n+1}$ in the logarithm of the formal group law classified by $h_\Z\colon L\into \Z[b_1,b_2,\dotsc]$) are manifestly adequate. We choose once and for all an adequate set of generators of $L$ and we write
\[\MH=\MGL/(a_1,a_2,\dotsc).\]
We remark that, for any morphism of base schemes $f\colon T\to S$, there is a canonical equivalence of oriented ring spectra $f^\ast\MGL_S\simeq\MGL_T$. In particular, the induced map $(\MGL_S)_\aast\to(\MGL_T)_\aast$ sends $a_n$ to $a_n$ and there is an induced equivalence $f^\ast\MH_S\simeq \MH_T$.

Recall that $H\Z$ has an orientation such that, if $\scr L$ is a line bundle over $X\in\Sm/S$, $c_1(\scr L)$ is the isomorphism class of $\scr L$ in the Picard group $\Pic(X)\cong H^{2,1}(X,\Z)$. It follows that the associated formal group law on $H\Z_\aast$ is additive and that there is a map $\MH\to H\Z$
factoring the universal Thom class $\th\colon\MGL\to H\Z$; this map is in fact unique for degree reasons.

\begin{lemma}\label{lem:HZ'effective}
	$\MH$ is connective.
\end{lemma}

\begin{proof}
    Follows from Corollary~\ref{cor:MGLeffective} since $\SH(S)_{\geq 0}$ is closed under homotopy colimits.
\end{proof}

\begin{lemma}\label{lem:HZeffective}
	$H\Z$ is connective.
\end{lemma}

\begin{proof}
	By Theorem~\ref{thm:EMLpb} and Lemma~\ref{lem:basechangetrunc}, we can assume that $S=\Spec k$ where $k$ is a perfect field.
	If $k\subset L$ is a finitely generated field extension, then by \cite[Lemma 3.9 and Theorem 3.6]{Mazza:2006} we have $\spi_{p,q}(H\Z)(\Spec L)=0$ for $p-q<0$. By Theorems~\ref{thm:fields} and~\ref{thm:t-structure}, this is equivalent to the connectivity of $H\Z$. 
\end{proof}

\begin{theorem}\label{thm:HZ_kk}
    The unit $\1\to H\Z$ induces an equivalence $(\1/\eta)_{\leq 0}\simeq H\Z_{\leq 0}$.
\end{theorem}

\begin{proof}
	By Theorem~\ref{thm:EMLpb} and Lemma~\ref{lem:basechangetrunc}, we can assume that $S=\Spec k$ where $k$ is a perfect field. By Theorem~\ref{thm:t-structure} and Lemma~\ref{lem:HZeffective}, it is necessary and sufficient to prove the exactness of the sequence
 \begin{tikzmath}
 	\def\colsep{.85em}
 	\diagram{\spi_{n-1,n-1}(\1) & \spi_{n,n}(\1) & \spi_{n,n}(H\Z) & 0 \\ };
 	\arrows (11-) edge node[above]{$\eta$} (-12) (12-) edge (-13) (13-) edge (-14);
 \end{tikzmath}
 for every $n\in\Z$.
    Furthermore, by Theorem~\ref{thm:fields}, it suffices to verify that the sequence is exact on the stalks at finitely generated field extensions $L$ of $k$. Write $\Spec L=\lim_\alpha X_\alpha$ where $X_\alpha\in\Sm/k$.
	 The Hopf map $\A^2\minus\{0\}\to\P^1$ defines a global section of the sheaf $\spi_{1,1}(\1)$ and in particular gives an element $\eta\in\spi_{1,1}(\1)(\Spec L)$. Since $\Hom_k(\Spec L,\G_m)=\colim_\alpha\Hom_k(X_\alpha,\G_m)$, any element $u\in L^\times$ defines a germ $[u]\in\spi_{-1,-1}(\1)(\Spec L)$. By \cite[Remark~6.42]{Morel:2012}, the graded ring $\bigoplus_{n\in\Z}\spi_{n,n}(\1)(\Spec L)$ is generated by the elements $\eta$ and $[u]$, and there is an exact sequence
\begin{tikzmath}
	\def\colsep{.85em}
	\diagram{\spi_{n-1,n-1}(\1)(\Spec L) & \spi_{n,n}(\1)(\Spec L) & K^M_{-n}(L) & 0\rlap,\\};
	\arrows (11-) edge node[above]{$\eta$} (-12) (12-) edge (-13) (13-) edge (-14);
\end{tikzmath}
    where the last map sends $[u]\in\spi_{-1,-1}(\1)(\Spec L)$ to the generator $\{u\}$ in the Milnor $K$-theory $K^M_{\ast}(L)$.
    On the other hand, combining \cite[Lemma 3.9]{Mazza:2006} with \cite[\S 5]{Mazza:2006}, we obtain an isomorphism of graded rings
    \[\lambda\colon K^M_{-\ast}(L)\to\bigoplus_{n\in\Z}\spi_{n,n}(H\Z)(\Spec L).\]
    It remains to prove that the triangle
\begin{tikzmath}
	\def\rowsep{2em}
	\diagram{\bigoplus_{n\in\Z}\spi_{n,n}(\1)(\Spec L) & \bigoplus_{n\in\Z}\spi_{n,n}(H\Z)(\Spec L) \\
	K^M_{-\ast}(L) & \\};
	\arrows
	(11-) edge (-12)
	([yshift=10pt] 11.south) edge (21)
	(21) edge node[below right]{$\lambda$} node[above left]{$\cong$} (12);
\end{tikzmath}
    is commutative. Since it is a triangle of graded rings, we can check its commutativity on the generators $\eta$ and $[u]$. The element $\eta$ maps to $0$ in both cases because it has positive degree. By inspection of the definition of $\lambda$, $\lambda\{u\}$ is the image of $u$ by the composition
	 \[L^\times=\colim_\alpha\Hom_k(X_\alpha,\G_m)\to\colim_\alpha[(X_\alpha)_+,\G_m]\to \colim_\alpha[(X_\alpha)_+,u_\tr\Z_\tr\G_m],\]
	 where the last map is induced by the unit $\G_m\to u_\tr\Z_\tr\G_m$. This is clearly also the image of $[u]$ by the unit $\1\to H\Z$.
\end{proof}

\begin{lemma}\label{lem:HZ<0}
    $\MGL_{\leq 0}\to H\Z_{\leq 0}$ is an equivalence.
\end{lemma}

\begin{proof}
	 Combine Theorems~\ref{thm:pikkMGL} and~\ref{thm:HZ_kk}, and the fact that $\MGL\to H\Z$ is a morphism of ring spectra.
\end{proof}

\begin{lemma}\label{lem:localhomology}
    Let $f$ be a map in $\SH(S)$. If $H\Q\wedge f$ and $H\Z/\ell\wedge f$ are equivalences for all primes $\ell$, then $H\Z\wedge f$ is an equivalence.
\end{lemma}

\begin{proof}
    Considering the cofiber $E$ of $f$, we are reduced to proving that if $H\Q\wedge E=0$ and $H\Z/\ell\wedge E=0$, then $H\Z\wedge E=0$. Since $\SH(S)$ is compactly generated, it suffices to show that $[X,H\Z\wedge E]=0$ for every compact $X\in \SH(S)$. By algebra, it suffices to prove that
\begin{gather*}
    [X,H\Z\wedge E]\tens\Q=0, \\
    [X,H\Z\wedge E]\tens\Z/\ell=0,\text{ and} \\
    \Tor^1([X,H\Z\wedge E],\Z/\ell)=0.
\end{gather*}
By Proposition~\ref{prop:Kcolim} \itemref{prop:Kcolim:2} and the compactness of $X$, we have $[X,H\Z\wedge E]\tens\Q =[X,H\Q\wedge E]$, which vanishes by assumption. The vanishing of the other two groups follows from the long exact sequence
\begin{tikzmath}
	\def\colsep{4.3em}
	\diagram{
	{[\Sigma^{1,0}X,{H\Z/\ell}\wedge E]} &[between origins] &[between origins] {[X,{H\Z}\wedge E]} &[between origins] &[between origins] {[X,{H\Z}\wedge E]} &[between origins] &[between origins] {[X,{H\Z/\ell}\wedge E]}\\
	& {\Tor^1([X,{H\Z}\wedge E],\Z/\ell)} & & & & {[X,{H\Z}\wedge E]\tens\Z/\ell} & \\
	};
	\arrows
	(11-) edge (-13) (13-) edge node[above]{$\ell$} (-15) (15-) edge (-17)
	(22) edge[<<-] (11) edge[c->] (13)
	(26) edge[<<-] (15) edge[c->] (17);
\end{tikzmath}
induced by the Bockstein cofiber sequence $H\Z\stackrel{\ell}{\to} H\Z\to H\Z/\ell$.
\end{proof}

Denote by $h\colon\MGL_\aast\to H\Q_\aast\MGL$ the rational Hurewicz map.

\begin{lemma}\label{lem:regular2}
	The sequence $(h(a_1),h(a_2),\dotsc)$ is regular in $H\Q_{\aast}\MGL$.
\end{lemma}

\begin{proof}
	By the calculus of oriented cohomology theories we have $H\Q_\aast\MGL\cong H\Q_\aast[b_1,b_2,\dotsc]$ and $h(a_n)\equiv u_nb_n$ modulo decomposables, for some $u_n\in\Z\minus \{0\}$. The lemma follows.
\end{proof}

\begin{lemma}\label{lem:HQ}
    $H\Q\wedge \MH\to H\Q\wedge H\Z$ is an equivalence.
\end{lemma}

\begin{proof}
Because $h(a_n)\in {H\Q}_{\ast\ast}\MGL$ maps to $0$ in both $H\Q_\aast \MH$ and $H\Q_\aast H\Z$, we have a commuting triangle
\begin{tikzmath}
	\def\colsep{3.5em}
	\diagram{
	&[between origins] {H\Q}_{\ast\ast}\MGL/(h(a_1),h(a_2),\dotsc) &[between origins] \\
	{H\Q}_{\ast\ast}\MH & & {H\Q}_{\ast\ast}H\Z\rlap.\\
	};
	\arrows (12) edge node[above left]{$\mu$} (21) edge node[above right]{$\nu$} (23) (21-) edge (-23);
\end{tikzmath}
    The map $\mu$ is an isomorphism by Lemma~\ref{lem:regular2}. By \cite[Corollary 10.3]{Naumann:2009}, $H\Q$ is the Landweber exact spectrum associated with the universal rational formal group law. In particular, $H\Q$ is cellular and the map $\MGL\to H\Q$ induces an isomorphism
\[{H\Z}_\aast \MGL\tens_{\Z[a_1,a_2,\dotsc]}\Q\cong {H\Z}_\aast H\Q.\]
This isomorphism can be identified with $\nu$ since ${H\Q}_\aast E\cong {H\Z}_\aast E\tens\Q$. 

This shows that $H\Q\wedge \MH\to H\Q\wedge H\Z$ is a $\pi_\aast$-isomorphism, whence an equivalence since both sides are cellular $H\Q$-modules (the right-hand side because $H\Q\wedge H\Z\simeq H\Q\wedge H\Q$).
\end{proof}

\begin{lemma}\label{lem:HZHZ}
    $H\Z\wedge \MH[1/c]\to H\Z\wedge H\Z[1/c]$
    is an equivalence.
\end{lemma}

\begin{proof}
    By Lemma~\ref{lem:localhomology}, it suffices to prove that
        \begin{gather*}
			  H\Q\wedge \MH[1/c]\to H\Q\wedge H\Z[1/c]\text{ and}\\
            H\Z/\ell\wedge \MH[1/c]\to H\Z/\ell\wedge H\Z[1/c]
        \end{gather*}
    are equivalences for every prime $\ell$. The former is taken care of by Lemma~\ref{lem:HQ}. In the latter, both sides are contractible if $\ell=c$. If $\ell\neq c$, then by Theorems~\ref{thm:intsteenrod} and~\ref{thm:AmodQi} (with $I=\{1,2,\dotsc\}$), $H\Z/\ell\wedge \MH\to H\Z/\ell\wedge H\Z$ is a $\pi_\aast$-isomorphism between cellular $H\Z/\ell$-modules, whence an equivalence (here we use the hypothesis that the generators $a_n$ are adequate in order to apply Theorem~\ref{thm:AmodQi}).
\end{proof}

\begin{lemma}\label{lem:main}
	Let $k$ be a field. Let $F\in\SH(k)$ be such that $H\Z\wedge F=0$ and let $X$ be a weak $\MGL$-module which is $r$-connective for some $r\in \Z$. Then $[F,X]=0$.
\end{lemma}

\begin{proof}
	By left completeness of the homotopy $t$-structure (Corollary~\ref{cor:nondeg}), there are fiber sequences of the form
    \begin{gather*}
        \prod_{n\in\Z}\Omega^{1,0}X_{\leq n}\to X\to\prod_{n\in\Z}X_{\leq n}\to\prod_{n\in\Z}X_{\leq n}\text{ and}\\
        \Omega^{1,0}X_{\leq n-1}\to \kappa_nX\to X_{\leq n}\to X_{\leq n-1}.
    \end{gather*}
    Since $X_{\leq n}=0$ if $n< r$, it suffices to show that $[\Sigma^{p,0}F, \kappa_nX]=0$ for every $p,n\in\Z$. Since $X$ is a weak $\MGL$-module, $\kappa_nX$ is a weak $\kappa_0\MGL$-module (see \S\ref{sub:tstructure}), so any $\Sigma^{p,0}F\to \kappa_nX$ can be factored as
\begin{tikzmath}
	\def\colsep{1.5em}
	\diagram{\Sigma^{p,0}F & \kappa_n X \\
	\kappa_0\MGL\wedge \Sigma^{p,0}F & \kappa_0\MGL\wedge \kappa_nX\rlap.\\};
	\arrows
	(11-) edge (-12) (21-) edge (-22) (11) edge (21) (22) edge (12);
\end{tikzmath}
Thus, it suffices to show that
\begin{equation}\label{eqn:smashnull}
\kappa_0\MGL\wedge F=0.
\end{equation}
By Corollary~\ref{cor:MGLeffective} and Lemma~\ref{lem:HZ<0}, $\kappa_0\MGL\simeq \MGL_{\leq0}\simeq H\Z_{\leq 0}$. By Lemma~\ref{lem:HZeffective} and \cite[Lemma 2.13]{Gutierrez:2012}, the canonical map $(H\Z\wedge H\Z)_{\leq 0}\to (H\Z_{\leq 0}\wedge H\Z)_{\leq 0}$ is an equivalence. Using the fact that the truncation $E\mapsto E_{\leq 0}$ is left adjoint to the inclusion $\SH(k)_{\leq 0}\subset\SH(k)$, it is easy to show that the composition
\[H\Z_{\leq 0}\to H\Z_{\leq 0}\wedge H\Z\to(H\Z_{\leq 0}\wedge H\Z)_{\leq 0}\simeq(H\Z\wedge H\Z)_{\leq 0}\to H\Z_{\leq 0}\]
is the identity, where the first map is induced by the unit $\1\to H\Z$ and the last one by the multiplication $H\Z\wedge H\Z\to H\Z$.
In particular, we get a factorization
\begin{tikzmath}
	\def\colsep{1.5em}
	\diagram{\kappa_0\MGL\wedge F & \kappa_0\MGL\wedge H\Z\wedge F \\
	& \kappa_0\MGL\wedge F\rlap. \\};
	\arrows
	(11-) edge (-12)
	(22) edge[<-] node[below left]{$\id$} (11) edge[<-] (12);
\end{tikzmath}
    Since $H\Z\wedge F=0$, this proves~\eqref{eqn:smashnull} and the lemma.
\end{proof}

\begin{theorem}\label{thm:main}
    $\MH[1/c]\to H\Z[1/c]$ is an equivalence.
\end{theorem}

\begin{proof}

Let $f\colon S\to\Spec k$ be essentially smooth, where $k$ is a field. Since $\th_S=f^\ast(\th_k)$ and $f^\ast(a_n)=a_n$, we may assume that $f$ is the identity.
Consider the fiber sequence
\[\Omega^{1,0}H\Z[1/c]\to F\to \MH[1/c]\to H\Z[1/c].\]
Then by Lemma~\ref{lem:HZHZ}, $H\Z\wedge F=0$.
Recall that $\MH$ is an $\MGL$-module by definition and that it is connective by Lemma~\ref{lem:HZ'effective}. By Lemma~\ref{lem:main}, we have
\begin{equation}\label{eqn:main1}
	[F,\MH[1/c]]=0.
\end{equation}
Similarly, $H\Z$ is a weak $\MGL$-module via the morphism of ring spectra $\th\colon\MGL\to H\Z$, and it is connective by Lemma~\ref{lem:HZeffective}. By Lemma~\ref{lem:main}, we have
\begin{equation}\label{eqn:main2}
	[F,\Omega^{1,0} H\Z[1/c]]=0.
\end{equation}
By~\eqref{eqn:main1}, the map $\Omega^{1,0}H\Z[1/c]\to F$ has a section, which is zero by~\eqref{eqn:main2}. Thus, $F=0$.
\end{proof}

\begin{remark}
	It follows from the results of~\S\ref{sub:convergence} that Theorem~\ref{thm:main} remains true if we drop the requirement that the generators $a_n$ be adequate.
\end{remark}

\section{Applications}

In this section we gather some consequences of Theorem~\ref{thm:main}. Throughout, the base scheme $S$ is essentially smooth over a field of characteristic exponent $c$. We denote by $L$ the Lazard ring which we regard as a graded ring with $\abs{a_n}=n$. Modules over $L$ will always be graded modules.

We note that if Theorem~\ref{thm:main} is true without inverting $c$, then so are all the results in this section.

\subsection{Cellularity of Eilenberg–Mac Lane spectra}

\begin{proposition}\label{prop:HZcell}
    For any $A\in\Sp(\s\Mod_{\Z[1/c]})$, the motivic Eilenberg–Mac Lane spectrum $HA$ is cellular.
\end{proposition}

\begin{proof}
    Since $\MGL$ is cellular, $\MGL/(a_1,a_2,\dotsc)[1/c]$ is also cellular. This proves the proposition for $A=\Z[1/c]$ by Theorem~\ref{thm:main}. The general case follows by Proposition~\ref{prop:Kcolim} \itemref{prop:Kcolim:2}.
\end{proof}

When $S$ is the spectrum of an algebraically closed field of characteristic zero, a different proof of the cellularity of $H\Z/2$ was given in \cite[Proposition 15]{Hu:2011}.

\subsection{The formal group law of algebraic cobordism}

\begin{proposition}\label{prop:piMGL}
	Suppose that $S$ is a field. Then the map $L[1/c]\to \MGL_{(2,1)\ast}[1/c]$ classifying the formal group law of $\MGL[1/c]$ is an isomorphism.
\end{proposition}

\begin{proof}
	We assume that $c=1$ to simplify the notations. We know from Lazard's theorem that the map is injective (see Remark~\ref{rmk:lazard}). For any $k\geq 0$, let $L(k)=L/(a_1,\dotsc,a_k)$ and let $\MGL(k)=\MGL/(a_1,\dotsc,a_k)$. We prove more generally that the induced map
	\begin{equation}\label{eqn:induction}
	L(k)_n\to\pi_{2n,n}\MGL(k)
	\end{equation}
	is surjective for all $n\in\Z$ and $k\geq 0$, and we proceed by induction on $n-k$. Since $\MGL(k)$ is connective, the map
	\[\pi_{2n,n}\MGL(k)\to\pi_{2n,n}\MGL(k+1)\]
	is an isomorphism when $n-k\leq 0$, by Corollary~\ref{cor:nondeg}. Taking the colimit as $k\to\infty$ and using Theorem~\ref{thm:main}, we obtain
	\[\pi_{2n,n}\MGL(k)\cong\pi_{2n,n}H\Z\]
	for $n-k\leq 0$, which proves that~\eqref{eqn:induction} is an isomorphism in this range since $\pi_{(2,1)\ast}H\Z$ carries the universal additive formal group law. If $n-k>0$, consider the commutative diagram with exact rows
	\begin{tikzmath}
		\diagram{L(k)_{n-k-1} & L(k)_n & L(k+1)_n \\
		\pi_{2(n-k-1),n-k-1}\MGL(k) & \pi_{2n,n}\MGL(k) & \pi_{2n,n}\MGL(k+1)\rlap. \\};
		\arrows (11-) edge[c->] node[above]{$a_{k+1}$} (-12) (12-) edge[->>] (-13)
		(21-) edge node[above]{$a_{k+1}$} (-22) (22-) edge (-23) (11) edge (21) (12) edge (22) (13) edge (23);
	\end{tikzmath}
	By induction hypothesis the left and right vertical maps are surjective. The five lemma then shows that~\eqref{eqn:induction} is surjective.
\end{proof}

\subsection{Slices of Landweber exact motivic spectra}

We now turn to some applications of the Hopkins–Morel equivalence to the slice filtration. Recall that $\SH^\eff(S)$ is the full subcategory of $\SH(S)$ generated under homotopy colimits and extensions by
\[\{\Sigma^{p,q} \Sigma^\infty X_+\suchthat X\in\Sm/S,\;p\in\Z,\;q\geq 0\}.\]
This is clearly a triangulated subcategory of $\SH(S)$. A spectrum $E\in\SH(S)$ is called \emph{$t$-effective} (or simply \emph{effective} if $t=0$) if $\Sigma^{0,-t}E\in\SH^\eff(S)$. The \emph{$t$-effective cover} of $E$ is the universal arrow $f_tE\to E$ from a $t$-effective spectrum to $E$, and the $t$th \emph{slice} $s_tE$ of $E$ is defined by the cofiber sequence
\[f_{t+1}E\to f_tE\to s_t E\to \Sigma^{1,0}f_{t+1}E.\]
Both functors $f_t\colon\SH(S)\to\SH(S)$ and $s_t\colon\SH(S)\to\SH(S)$ are triangulated and preserve homotopy colimits. By Remark~\ref{rmk:slices} and \cite[Theorem 5.2 (i)]{Gutierrez:2012}, $s_t$ has a canonical lift to a functor $\SH(S)\to\D(H\Z)$.

It is clear that the slice filtration is exhaustive in the sense that, for every $E\in\SH(S)$,
\begin{equation}\label{eqn:exhaustive}
E\simeq\hocolim_{t\to -\infty}f_tE.
\end{equation}
We say that $E$ is \emph{complete} if
\[
\holim_{t\to\infty} f_tE=0.
\]
Thus, by~\eqref{eqn:exhaustive}, slices detect equivalences between complete spectra. See \cite[Remark 2.1]{Voevodsky:2002} for an example of a spectrum which is not complete.

We refer to \cite{Naumann:2009} for the general theory of Landweber exact motivic spectra.

\begin{theorem}\label{thm:slicesLandweber}
	Let $M_\ast$ be a Landweber exact $L[1/c]$-module and $E\in\SH(S)$ the associated motivic spectrum. Then there is a unique equivalence of $H\Z$-modules $s_tE\simeq \Sigma^{2t,t}HM_t$ such that the diagram
	\begin{tikzmath}
		\diagram{\pi_{0,0}\MGL\tens M_t & \pi_{0,0}HM_t \\ \pi_{2t,t}E & \pi_{2t,t}s_tE \\};
		\arrows (11-) edge node[above]{$\th$} (-12) (11) edge (21) (21-) edge (-22) (12) edge node[right]{$\cong$} (22);
	\end{tikzmath}
	commutes.
\end{theorem}

\begin{proof}
	Suppose first that $c=1$ and $M_\ast=L$, so that $E=\MGL$. The assertion then follows from Theorem~\ref{thm:main} and \cite[Corollary 4.7]{Spitzweck:2010}, where the assumption that the base is a perfect field can be removed in view of Remark~\ref{rmk:slices}. In positive characteristic, it is easy to see that the proofs in \cite{Spitzweck:2010} are unaffected by the inversion of $c$ in Theorem~\ref{thm:main} and yield the statement of the theorem for $M_\ast=L[1/c]$. The result for general $M_\ast$ follows as in \cite[Theorem 6.1]{Spitzweck:2012}.
\end{proof}

For $M_\ast=L[1/c]$, we obtain in particular
\begin{equation}\label{eqn:sliceMGL}
	s_t\MGL[1/c]\simeq \Sigma^{2t,t}HL_t[1/c].
\end{equation}

\subsection{Slices of the motivic sphere spectrum}

Let $\mathbb L$ be the graded cosimplicial commutative ring associated with the Hopf algebroid $(L,LB)$ that classifies formal group laws and strict isomorphisms. For each $t\in\Z$, we can view its degree $t$ summand $\mathbb L_t$ as a chain complex (concentrated in nonpositive degrees) via the dual Dold–Kan correspondence, whence as an object in $\Sp(\s\Ab)$.

\begin{theorem}
	There is an equivalence of $H\Z$-modules $s_t\1[1/c]\simeq \Sigma^{2t,t}H\mathbb L_t[1/c]$ such that the diagram
	\begin{tikzmath}
		\diagram{\Sigma^{2t,t}H\mathbb L_t[1/c] & \Sigma^{2t,t}HL_t[1/c] \\ s_t\1[1/c] & s_t\MGL[1/c] \\};
		\arrows (11-) edge (-12) (11) edge node[left]{$\simeq$} (21) (12) edge node[left]{\eqref{eqn:sliceMGL}} node[right]{$\simeq$} (22) (21-) edge node[below]{\textnormal{unit}} (-22);
	\end{tikzmath}
	commutes.
\end{theorem}

\begin{proof}
	Follows from~\eqref{eqn:sliceMGL} as in \cite[\S8]{Levine:2013}.
\end{proof}

\subsection{Convergence of the slice spectral sequence}
\label{sub:convergence}

Consider a homological functor $F\colon\SH(S)\to\scr A$ where $\scr A$ is a bicomplete abelian category. The tower
\[\dotsb\to f_{t+1}\to f_t \to f_{t-1}\to\dotsb\]
in the triangulated category $\SH(S)$ gives rise in the usual way to a bigraded spectral sequence $\{F_r^{\aast}\}_{r\geq 1}$ with
\[
F_1^{s,t}(E)=F(s_t\Sigma^s E)\quad\text{and}\quad d_r\colon F_r^{s,t}\to F_r^{s+1,t+r}.
\]
We denote by $F_\infty^{\aast}$ the limit of this spectral sequence.
On the other hand, this tower induces a filtration of the functor $F$ by the subfunctors $f_tF$ given by
\[f_tF(E)=\Im(F(f_tE)\to F(E)),\]
and we denote by $s_tF$ the quotient $f_t F/f_{t+1}F$. We then have a canonical relation
\begin{equation}\label{eqn:foo}
s_tF(\Sigma^s E)\longleftrightarrow F_\infty^{s,t}(E),
\end{equation}
and we say that $E$ is \emph{convergent with respect to $F$} if it is an isomorphism for all $s,t\in\Z$. We say that $E$ is \emph{left bounded with respect to $F$} if for every $s\in\Z$, $F(f_t\Sigma^sE)=0$ for $t\gg 0$; this implies that \eqref{eqn:foo} is an epimorphism from the right-hand side to the left-hand side. If $F$ preserves sequential homotopy colimits and if sequential colimits are exact in $\scr A$, \eqref{eqn:foo} is always a monomorphism from the left-hand side to the right-hand side by virtue of~\eqref{eqn:exhaustive}, so that left boundedness implies convergence.
We simply say that a spectrum is \emph{convergent} (\resp{} \emph{left bounded}) if it is convergent (\resp{} left bounded) with respect to the functors $[\Sigma^{0,q}\Sigma^\infty X_+,\ph]$ for all $X\in\Sm/S$ and $q\in\Z$. Note that every left bounded spectrum is also complete.

\begin{lemma}\label{lem:sss}
	Let $k$ be a field and $E\in\SH(k)$. Suppose that there exists $n\in\Z$ such that $f_tE$ is $(t+n)$-connective for all $t\in\Z$. Then, for any essentially smooth morphism $f\colon S\to \Spec k$, $f^\ast E\in\SH(S)$ is left bounded and in particular complete and convergent.
\end{lemma}

\begin{proof}
	Let $X\in\Sm/S$ and $p,q\in\Z$. Suppose first that $f$ is the identity. Since $f_tE$ is $(t+n)$-connective, we have $[\Sigma^{p,q}\Sigma^\infty X_+,f_tE]=0$ as soon as $t>p-q-n+\dim X$ (Corollary~\ref{cor:nondeg}), so $E$ is left bounded. In general, let $f$ be the cofiltered limit of smooth maps $f_\alpha\colon S_\alpha\to\Spec k$ and let $X$ be the limit of smooth $S_\alpha$-schemes $X_\alpha$. Then by Lemma~\ref{lem:stablecontinuity}, we have
	\[
	[\Sigma^{p,q}\Sigma^\infty X_+,f^\ast(f_tE)]\cong\colim_\alpha [\Sigma^{p,q}\Sigma^\infty(X_\alpha)_+,f_\alpha^\ast(f_tE)]\cong \colim_\alpha [\Sigma^{p,q}\Sigma^\infty(X_\alpha)_+,f_tE]
	\]
	which is zero if $t>p-q-n+\essdim X$.
	
	It remains to show that $f^\ast(f_tE)\to f^\ast(E)$ is the $t$-effective cover of $f^\ast (E)$. Since $f^\ast(f_tE)$ is $t$-effective, it suffices to show that for any $X\in\Sm/S$, $p\in\Z$, and $q\geq t$, the map
	\[[\Sigma^{p,q}\Sigma^\infty X_+,f^\ast (f_tE)]\to [\Sigma^{p,q}\Sigma^\infty X_+, f^\ast(E)]\]
	is an isomorphism. If $f$ is smooth, this follows from the fact that the left adjoint $f_\sharp$ to $f^\ast$ preserves $t$-effective objects. In general, it follows from Lemma~\ref{lem:stablecontinuity}.
\end{proof}

\begin{lemma}\label{lem:fqMGL}
	$f_t\MGL[1/c]$ is $t$-connective.
\end{lemma}

\begin{proof}
	By Theorem~\ref{thm:main} and the proof of \cite[Theorem 4.6]{Spitzweck:2010}, $f_t\MGL[1/c]$ is the homotopy colimit of a diagram of $\MGL$-modules of the form $\Sigma^{2n,n}\MGL[1/c]$ with $n\geq t$. Thus, $f_t\MGL[1/c]$ is a homotopy colimit of $t$-connective spectra and hence is $t$-connective.
\end{proof}

\begin{lemma}\label{lem:fqLandweber}
	Let $M_\ast$ be a Landweber exact $L[1/c]$-module and $E\in\SH(S)$ the associated motivic spectrum. Then $f_tE$ is $t$-connective.
\end{lemma}

\begin{proof}
	Suppose first that $M_\ast$ is a flat $L$-module. It is then a filtered colimit of finite sums of shifts of $L[1/c]$, and $E$ is equivalent to the filtered homotopy colimit of a corresponding diagram in $\Mod_{\MGL}$. By Lemma~\ref{lem:fqMGL}, $f_t(\Sigma^{2n,n}\MGL[1/c])\simeq\Sigma^{2n,n}f_{t-n}\MGL[1/c]$ is $t$-connective for any $n\in\Z$. Since $f_t$ commutes with homotopy colimits, $f_t E$ is $t$-connective. In general, $E$ is a retract in $\SH(S)$ of $\MGL\wedge E$, so it suffices to show that $f_t(\MGL\wedge E)$ is $t$-connective. But $\MGL\wedge E$ is the spectrum associated with the Landweber exact left $L$-module $LB\tens_LM_\ast$ which is flat since it is the pullback of $M_\ast$ by $\Spec L\to\scr M_{\mathrm{fg}}^s$, where $\scr M_{\mathrm{fg}}^s$ is the stack represented by the Hopf algebroid $(L,LB)$.
\end{proof}

\begin{theorem}\label{thm:convergence}
	Let $M_\ast$ be a Landweber exact $L[1/c]$-module and $E\in\SH(S)$ the associated motivic spectrum. Then $E$ is left bounded and in particular complete and convergent.
\end{theorem}

\begin{proof}
	Since Landweber exact spectra are cartesian sections of $\SH(\ph)$ by definition, this follows from Lemmas~\ref{lem:sss} and~\ref{lem:fqLandweber}.
\end{proof}

\begin{example}\label{ex:landweber}
	Many interesting Landweber exact $L$-algebras are of the form $(L/I)[J^{-1}]$ where $I$ is a regular sequence of homogeneous elements and $J\subset L/I$ is a regular multiplicative subset. If $E$ is the Landweber exact motivic spectrum associated with $(L/I)[J^{-1}]$, there is a map
\[(\MGL/I)[J^{-1}]\to E\]
in $\D(\MGL)$. Assuming $c\in J$, we claim that this map is an equivalence. By Lemmas~\ref{lem:sss} and~\ref{lem:fqMGL}, the $\MGL$-module $(\MGL/I)[J^{-1}]$ is left bounded and hence complete, and so is $E$ by Theorem~\ref{thm:convergence}. Thus, it suffices to prove that this map is a slicewise equivalence, and this follows easily from Theorem~\ref{thm:slicesLandweber}. If $J_0\subset J$ is the subset of degree $0$ elements, we can prove in the same way that $f_0E\simeq(\MGL/I)[J_0^{-1}]$.
\end{example}

Combining Theorem~\ref{thm:convergence} with Theorem~\ref{thm:slicesLandweber}, we obtain for every $X\in\Sm/S$ a spectral sequence starting at $H^{\aast}(X,M_\ast)$ whose $\infty$-page is the associated graded of a complete filtration on $E^\aast(X)$:
\[H^{p+2t,q+t}(X,M_t)\Rightarrow E^{p,q}(X).\]
Note that $H^\aast(X,M_t)\cong H^\aast(X,\Z)\tens M_t$ since $M_t$ is torsion-free.
When $M_\ast=L[1/c]$, this spectral sequence takes the form
\begin{equation}\label{eqn:sssMGL}
	H^{p+2t,q+t}(X,\Z)\tens L_{t}[1/c]\Rightarrow \MGL^{p,q}(X)[1/c].
\end{equation}

In \cite{Levine:2007}, Levine and Morel define a multiplicative cohomology theory $\Omega^\ast(\ph)$ for smooth schemes over a field of characteristic zero and they make the following conjecture which can now be settled:

\begin{corollary}
	Let $k$ be a field of characteristic zero and $X\in\Sm/k$. There is a natural isomorphism of graded rings
	\[\Omega^\ast(X)\cong\MGL^{(2,1)\ast}(X).\]
\end{corollary}

\begin{proof}
	This is proved in \cite{Levine:2009} assuming the existence of the spectral sequence~\eqref{eqn:sssMGL}.
\end{proof}

\appendix
\section{Essentially smooth base change}
\label{sec:appendix}

In this appendix we show that the categories of motivic spaces, spaces with transfers, spectra, and spectra with transfers are ``continuous'' with respect to inverse limits of smooth morphisms of base schemes. Smooth morphisms and étale morphisms are always separated and of finite type.

\begin{definition}\label{def:esssmooth}
	Let $S$ be a base scheme. A morphism of schemes $T\to S$ is \emph{essentially smooth} if $T$ is a base scheme and if $T$ is a cofiltered limit $\lim_\alpha T_\alpha$ of smooth $S$-schemes where the transition maps $T_\beta\to T_\alpha$ are affine and dominant.
\end{definition}

The dominance condition is needed in the proof of Lemma~\ref{lem:scontinuity} \itemref{lem:scontinuity:2} below. If $X$ is smooth and quasi-projective over a field $k$ and $Z\subset X$ is a finite subset, then the semi-local schemes $\Spec \scr O_{X,Z}$, $\Spec \scr O_{X,Z}^\mathit{h}$, and $\Spec \scr O_{X,Z}^\mathit{sh}$ are examples of essentially smooth schemes over $k$.

With the notations of Definition~\ref{def:esssmooth}, if $U$ is any $T$-scheme of finite type, then by \cite[Théorème 8.8.2]{EGA4-3} it is the limit of a diagram of schemes of finite type $(U_\alpha)$ over the diagram $(T_\alpha)$. Moreover, if the morphism $U\to T$ is either
\begin{itemize}
	\item separated,
	\item smooth or étale,
	\item an open immersion or a closed immersion,
\end{itemize}
then we can choose each $U_\alpha\to T_\alpha$ to have the same property (this follows from \cite[Proposition 8.10.4]{EGA4-3}, \cite[Proposition 17.7.8]{EGA4-4}, and \cite[Proposition 8.6.3]{EGA4-3}, respectively). In particular, a composition of essentially smooth morphisms is essentially smooth. The following lemma shows that an essentially smooth scheme over a field is in fact essentially smooth over a finite field $\bb F_p$ or over $\Q$.

\begin{lemma}\label{lem:imperfectext}
	Let $k$ be a perfect field and $L$ a field extension of $k$. Then the morphism $\Spec L\to\Spec k$ is essentially smooth.
\end{lemma}

\begin{proof}
	We have $\Spec L =\lim_K\Spec K$ where $K$ ranges over all finitely generated extensions of $k$ contained in $L$. We may therefore assume that $L=k(x_1,\dotsc,x_n)$ for some $x_i\in L$. Since $k$ is perfect, $\Spec k[x_1,\dotsc,x_n]$ has a smooth dense open subset $U$ (\cite[Corollaire 17.15.13]{EGA4-4}). Then $\Spec L$ is the cofiltered limit of the nonempty affine open subschemes of $U$.
\end{proof}

From now on we fix a commutative ring $R$. We let $\H^s_\pt(S)$ (\resp{} $\H^s_\tr(S,R)$) denote the homotopy category of the category of pointed simplicial presheaves on $\Sm/S$ (\resp{} additive simplicial presheaves on $\Cor(S,R)$) with the projective model structure. Mapping spaces in the homotopy categories $\H^s_\pt(S)$ and $\H^s_\tr(S,R)$ will be denoted by $\Map^s(X,Y)$ to distinguish them from mapping spaces in $\H_\pt(S)$ and $\H_\tr(S,R)$, which we simply denote by $\Map(X,Y)$.

Let $\cat C(S)$ be any of the categories $\H^s_\pt(S)$, $\H^s_\tr(S,R)$, $\H_\pt(S)$, $\H_\tr(S,R)$, $\SH(S)$, and $\SH_\tr(S,R)$. In the terminology of \cite[\S1.1]{Cisinski:2012}, $\cat C(\ph)$ is then a complete monoidal $\Sm$-fibered category over the category of base schemes. In particular, a morphism of base schemes $f\colon T\to S$ induces a symmetric monoidal adjunction
\begin{tikzmath}
	\diagram{\cat C(S) & \cat C(T) \\};
	\arrows (11-) edge[vshift=\dbl] node[above=\dbl]{$f^\ast$} (-12) (-12) edge[vshift=\dbl] node[below=\dbl]{$f_\ast$} (11-);
\end{tikzmath}
where $f^\ast$ is induced by the base change functor $\Sm/S\to\Sm/T$ or $\Cor(S,R)\to\Cor(T,R)$, and if $f$ is smooth, it induces a further adjunction
\begin{tikzmath}
	\diagram{\cat C(T) & \cat C(S) \\};
	\arrows (11-) edge[vshift=\dbl] node[above=\dbl]{$f_\sharp$} (-12) (-12) edge[vshift=\dbl] node[below=\dbl]{$f^\ast$} (11-);
\end{tikzmath}
where $f_\sharp$ is induced by the forgetful functor $\Sm/T\to\Sm/S$ or $\Cor(T,R)\to\Cor(S,R)$. All this structure can in fact be defined at the level of model categories, and while we will not directly use any model structures, we will use homotopy limits and colimits. In other words, we consider $\cat C(S)$ as a derivator rather than just a homotopy category. The above adjunctions are then adjunctions of derivators in the sense that the left adjoints preserve homotopy colimits and the right adjoints preserve homotopy limits.

The (symmetric monoidal) adjunctions
\begin{tikzmath}
	\diagram{\H^s_\pt(\ph) & \H^s_\tr(\ph,R) \\ \H_\pt(\ph) & \H_\tr(\ph,R) \\ \SH(\ph) & \SH_\tr(\ph,R) \\};
	\arrows
	(11) edge[vshift=\dbl] node[above=\dbl]{$\L R_\tr$} (12) edge[<-,vshift=-\dbl] node[below=\dbl]{$u_\tr$} (12)
	(11) edge[vshift=-\dbl] node[left=\dbl]{$\L\id$} (21) edge[<-left hook,vshift=\dbl] node[right=\dbl]{$\R\id$} (21)
	(12) edge[vshift=-\dbl] node[left=\dbl]{$\L\id$} (22) edge[<-left hook,vshift=\dbl] node[right=\dbl]{$\R\id$} (22)
	(21) edge[vshift=\dbl] node[above=\dbl]{$\L R_\tr$} (22) edge[<-,vshift=-\dbl] node[below=\dbl]{$u_\tr$} (22)
	(31) edge[vshift=\dbl] node[above=\dbl]{$\L R_\tr$} (32) edge[<-,vshift=-\dbl] node[below=\dbl]{$u_\tr$} (32)
	(21) edge[vshift=-\dbl] node[left=\dbl]{$\Sigma^\infty$} (31) edge[<-,vshift=\dbl] node[right=\dbl]{$\R\Omega^\infty$} (31)
	(22) edge[vshift=-\dbl] node[left=\dbl]{$\Sigma^\infty_\tr$} (32) edge[<-,vshift=\dbl] node[right=\dbl]{$\R\Omega^\infty_\tr$} (32);
\end{tikzmath}
are compatible with the $\Sm$-fibered structures. This means that the left adjoint functors always commute with $f^\ast$, and, if $f$ is smooth, they also commute with $f_\sharp$. For the adjunctions $(\L R_\tr,u_\tr)$, this follows from the commutativity of the squares
\[
\begin{tikzpicture}
	\diagram{\Sm/S & \Cor(S,R) \\ \Sm/T & \Cor(T,R) \\};
	\arrows (11-) edge node[above]{$\Gamma$} (-12) (11) edge node[left]{$f^\ast$} (21) (21-) edge node[below]{$\Gamma$} (-22) (12) edge node[right]{$f^\ast$} (22);
\end{tikzpicture}
\qquad\text{and}\qquad
\begin{tikzpicture}
	\diagram{\Sm/T & \Cor(T,R) \\ \Sm/S & \Cor(S,R) \\};
	\arrows (11-) edge node[above]{$\Gamma$} (-12) (11) edge node[left]{$f_\sharp$} (21) (21-) edge node[below]{$\Gamma$} (-22) (12) edge node[right]{$f_\sharp$} (22);
\end{tikzpicture}
\]
(\cite[Lemmas 9.3.3 and 9.3.7]{Cisinski:2012}). For the vertical adjunctions this holds by definition of $f^\ast$ and of $f_\sharp$.

From now on we fix an essentially smooth morphism of base schemes $f\colon T\to S$, cofiltered limit of smooth morphisms $f_\alpha\colon T_\alpha\to S$ as in Definition~\ref{def:esssmooth}.

\begin{lemma}\label{lem:scontinuity}
	Let $K$ be a finite simplicial set and let $d\colon K\to \operatorname{N}(\Sm/T)$ be a diagram of smooth $T$-schemes, cofiltered limit of diagrams $d_\alpha\colon K\to \operatorname{N}(\Sm/T_\alpha)$. Let $X$ (\resp{} $X_\alpha$) be the homotopy colimit of $d$ in $\H_\pt^s(T)$ (\resp{} of $d_\alpha$ in $\H_\pt^s(T_\alpha)$).
	\begin{enumerate}
		\item\label{lem:scontinuity:1} For any $\scr F\in\H_\pt^s(S)$, the canonical map \[\hocolim_\alpha \Map^s(X_\alpha,f_\alpha^\ast\scr F)\to \Map^s(X,f^\ast\scr F)\] is an equivalence.
		\item\label{lem:scontinuity:2} For any $\scr F\in\H_\tr^s(S,R)$, the canonical map \[\hocolim_\alpha \Map^s(\L R_\tr X_\alpha,f_\alpha^\ast\scr F)\to \Map^s(\L R_\tr X,f^\ast\scr F)\] is an equivalence.
	\end{enumerate}
\end{lemma}

\begin{proof}
	Since filtered homotopy colimits commute with finite homotopy limits, we can assume that $K=\Delta^0$. Both sides preserve homotopy colimits in $\scr F$, so we may further assume that $\scr F=Y_+$ (\resp{} that $\scr F=\L R_\tr Y_+$) where $Y\in\Sm/S$. Then $f^\ast\scr F$ is represented by $Y\times_ST$ and the claim follows from \cite[Théorème 8.8.2]{EGA4-3} (\resp{} from \cite[Proposition 9.3.9]{Cisinski:2012}).
\end{proof}

A cartesian square
\begin{tikzmath}
	\diagram{W & V \\ U & X \\};
	\arrows (11-) edge (-12) (11) edge (21) (21-) edge node[above]{$i$} (-22) (12) edge node[right]{$p$} (22);
\end{tikzmath}
in $\Sm/S$ will be called a \emph{Nisnevich square} if $i$ is an open immersion, $p$ is étale, and $p$ induces an isomorphism $Z\times_XV\cong Z$, where $Z$ is the reduced complement of $i(U)$ in $X$ (by \cite[Proposition 17.5.7]{EGA4-4}, $Z\times_XV$ is always reduced and so it is the reduced complement of $p^{-1}(i(U))$ in $V$).

Recall that an object $\scr F$ in $\H^s_\pt(S)$ or $\H^s_\tr(S,R)$ is called \emph{$\A^1$-local} if, for every $X\in\Sm/S$, the projection $X\times\A^1\to X$ induces an equivalence $\scr F(X)\simeq \scr F(X\times\A^1)$, and it is \emph{Nisnevich-local} if it satisfies homotopical Nisnevich descent. Since $S$ is Noetherian and of finite Krull dimension, $\scr F$ is Nisnevich-local if and only if $\scr F(\emptyset)$ is contractible and, for every Nisnevich square $Q$, the square $\scr F(Q)$ is homotopy cartesian (this is a reformulation of \cite[Proposition 3.1.16]{Morel:1999}). The localization functors
\[
	\H^s_\pt(S)\to\H_\pt(S)\quad\text{and}\quad\H^s_\tr(S,R)\to \H_\tr(S,R)
\]
have fully faithful right adjoints identifying $\H_\pt(S)$ and $\H_\tr(S,R)$ with the full subcategories of $\A^1$- and Nisnevich-local objects in $\H_\pt^s(S)$ and $\H^s_\tr(S,R)$, respectively.

We now make the following observations.
\begin{itemize}
	\item Any trivial line bundle in $\Sm/T$ is the cofiltered limit of trivial line bundles in $\Sm/T_\alpha$.
	\item Any Nisnevich square in $\Sm/T$ is the cofiltered limit of Nisnevich squares in $\Sm/T_\alpha$.
\end{itemize}
The first one is obvious. Any Nisnevich square in $\Sm/T$ is a cofiltered limit of cartesian squares
\begin{tikzmath}
	\diagram{W_\alpha & V_\alpha \\ U_\alpha & X_\alpha \\};
	\arrows (11-) edge (-12) (11) edge (21) (21-) edge node[above]{$i_\alpha$} (-22) (12) edge node[right]{$p_\alpha$} (22);
\end{tikzmath}
in $\Sm/T_\alpha$, where $i_\alpha$ is an open immersion and $p_\alpha$ is étale. Let $Z_\alpha$ be the reduced complement of $i_\alpha(U_\alpha)$ in $X_\alpha$. It remains to show that $Z_\alpha\times_{X_\alpha}V_\alpha\to Z_\alpha$ is eventually an isomorphism. By \cite[Corollaire 8.8.2.5]{EGA4-3}, it suffices to show that $Z= \lim_\alpha Z_\alpha$ as closed subschemes of $X$. Now $\lim_\alpha Z_\alpha\cong Z_\alpha\times_{X_\alpha} X$ for large $\alpha$, and so $\lim_\alpha Z_\alpha$ is a closed subscheme of $X$ with the same support as $Z$. Moreover, it is reduced by \cite[Proposition 8.7.1]{EGA4-3}, so it coincides with $Z$.

\begin{lemma}\label{lem:A1Nislocal}
	The functors $f^\ast\colon\H^s_\pt(S)\to\H^s_\pt(T)$ and $f^\ast\colon\H^s_\tr(S,R)\to\H^s_\tr(T,R)$ preserve $\A^1$-local objects and Nisnevich-local objects.
\end{lemma}

\begin{proof}
	If $f$ is smooth this follows from the existence of the left adjoint $f_\sharp$ to $f^\ast$ and the observation that $f_\sharp$ sends trivial line bundles to trivial line bundles and Nisnevish squares to Nisnevich squares. Thus, each $f_\alpha^\ast$ preserves $\A^1$-local objects and Nisnevich-local objects. Since any trivial line bundle (\resp{} Nisnevich square) over $T$ is a cofiltered limit of trivial line bundles (\resp{} Nisnevich squares) over $T_\alpha$, Lemma~\ref{lem:scontinuity} shows that $f^\ast$ preserves $\A^1$-local objects and Nisnevich-local objects in general.
\end{proof}

\begin{lemma}\label{lem:continuity}
	Let $K$ be a finite simplicial set and let $d\colon K\to \operatorname{N}(\Sm/T)$ be a diagram of smooth $T$-schemes, cofiltered limit of diagrams $d_\alpha\colon K\to \operatorname{N}(\Sm/T_\alpha)$. Let $X$ (\resp{} $X_\alpha$) be the homotopy colimit of $d$ in $\H_\pt(T)$ (\resp{} of $d_\alpha$ in $\H_\pt(T_\alpha)$).
	\begin{enumerate}
		\item\label{lem:continuity:1} For any $\scr F\in\H_\pt(S)$, the canonical map \[\hocolim_\alpha \Map(X_\alpha,f_\alpha^\ast\scr F)\to \Map(X,f^\ast\scr F)\] is an equivalence.
		\item\label{lem:continuity:2} For any $\scr F\in\H_\tr(S,R)$, the canonical map \[\hocolim_\alpha \Map(\L R_\tr X_\alpha,f_\alpha^\ast\scr F)\to \Map(\L R_\tr X,f^\ast\scr F)\] is an equivalence.
	\end{enumerate}
\end{lemma}

\begin{proof}
	Combine Lemmas \ref{lem:scontinuity} and~\ref{lem:A1Nislocal}.
\end{proof}

It is now easy to deduce a stable version of Lemma~\ref{lem:continuity}. Recall that objects in $\SH(S)$ and $\SH_\tr(S,R)$ can be modeled by $\Omega$-spectra, \ie, sequences $(E_0,E_1,\dotsc)$ of $\A^1$- and Nisnevich-local objects $E_i$ with equivalences $E_i\simeq\Omega^{2,1}E_{i+1}$. If $E\in\SH(S)$ is represented by the $\Omega$-spectrum $(E_0,E_1,\dotsc)$ and $X\in\H_\pt(S)$, we have
\begin{equation}\label{eqn:Omega}
	[\Sigma^{p,q}\Sigma^\infty X,E]=[\Sigma^{p+2r,q+r}X,E_r]
\end{equation}
for any $r\geq 0$ such that $p+2r\geq q+r\geq 0$, and similarly if $E\in\SH_\tr(S,R)$ and $X\in\H_\tr(S,R)$.

If $f$ is smooth, then the existence of the left adjoint $f_\sharp$ to $f^\ast$ shows that $f^\ast$ commutes with the unstable bigraded loop functors $\Omega^{p,q}$. An easy application of Lemma~\ref{lem:continuity} then shows that this is still true for $f$ essentially smooth. Thus, the base change functors
\[f^\ast\colon\SH(S)\to\SH(T)\quad\text{and}\quad f^\ast\colon\SH_\tr(S,R)\to\SH_\tr(T,R)\]
can be described explicitly as sending an $\Omega$-spectrum $(E_0,E_1,\dotsc)$ to the $\Omega$-spectrum $(f^\ast E_0,f^\ast E_1,\dotsc)$. From~\eqref{eqn:Omega} and Lemma~\ref{lem:continuity} we obtain the following result.

\begin{samepage}
\begin{lemma}\label{lem:stablecontinuity}
	Let $X\in\Sm/T$ be a cofiltered limit of smooth $T_\alpha$-schemes $X_\alpha$ and let $p,q\in\Z$.
	\begin{enumerate}
		\item\label{lem:stablecontinuity:1} For any $E\in\SH(S)$, $[\Sigma^{p,q}\Sigma^\infty X_+,f^\ast E]\cong\colim_\alpha[\Sigma^{p,q}\Sigma^\infty(X_\alpha)_+,f_\alpha^\ast E]$.
		\item\label{lem:stablecontinuity:2} For any $E\in\SH_\tr(S,R)$, $[\L R_\tr\Sigma^{p,q}\Sigma^\infty X_+,f^\ast E]\cong\colim_\alpha[\L R_\tr\Sigma^{p,q}\Sigma^\infty(X_\alpha)_+,f_\alpha^\ast E]$.
	\end{enumerate}
\end{lemma}
\end{samepage}

\providecommand{\bysame}{\leavevmode\hbox to3em{\hrulefill}\thinspace}

\end{document}